\documentclass[12pt]{article}
\usepackage[margin=1in,vmargin=1in]{geometry}
\usepackage{graphicx}
\usepackage{hyperref}
\usepackage{caption}

\usepackage{url}
\usepackage{amsthm}
\usepackage{amssymb,amsmath}
\usepackage{amsmath}
\usepackage{fancyhdr}
\usepackage[final]{pdfpages}
\usepackage{authblk}
\usepackage[export]{adjustbox}
\usepackage{blindtext}

\theoremstyle{plain}
\newtheorem{theorem}{Theorem}[section]

\newtheorem{conjecture}[theorem]{Conjecture}
\newtheorem{remark}[theorem]{Remark}

\theoremstyle{definition}
\newtheorem{definition}[theorem]{Definition}

\theoremstyle{assumption}

\usepackage{float}
\usepackage[caption = false]{subfig}

\usepackage{stackengine}[2013-09-11]

\usepackage{scalerel}

\usepackage{commath}

 \stackMath

\title{
Resonant Tori, Transport Barriers, and Chaos \\
in a Vector Field with a Neimark-Sacker Bifurcation
\thanks{The second author was partially supported by NSF grant 
DMS-1813501.}}
\providecommand{\keywords}[1]{\textbf{Keywords:} #1}

\author[]{Emmanuel Fleurantin, J.D. Mireles James \vspace{-1ex}}
\affil[]{\small{Department of Mathematics, Florida Atlantic University, 777 Glades Rd, Boca Raton, FL 33431} \vspace{-2ex}}
\affil[]{\small{email: efleurantin2013@fau.edu, jmirelesjames@fau.edu}\vspace{-4ex}}
\date{}

\newcommand{\correction}[2]{#2}
\begin{document}
\maketitle
\begin{abstract}
    \noindent    
    We make a detailed numerical study of a three dimensional
     dissipative vector field
    derived from the normal form for a cusp-Hopf
    bifurcation.  The vector field 
     exhibits a Neimark-Sacker bifurcation giving
    rise to an attracting invariant torus.  Our 
    main goals are to (A) follow the torus via parameter continuation
    from its appearance to its disappearance, studying 
    its dynamics between these events, and to (B) 
    study the embeddings of the stable/unstable manifolds of the hyperbolic 
    equilibrium solutions over this parameter range, focusing on 
    their role as transport barriers and their participation in global 
    bifurcations.  Taken together the results
    highlight the main features of the global dynamics of the system.  
    
\end{abstract}

\begin{flushleft}
\keywords{Invariant tori, Neimark-Sacker bifurcation, parameterization method \\
\textbf{AMS Subject Classifications:} 34C45, 37G35, 37M05, 37C55}
\\
\end{flushleft}

\section{Introduction}\label{sec:intro}	
\noindent 
Interactions between equilibrium and  
oscillating states provide a basic mechanism
for generating complicated dynamics in nonlinear systems.  
Such interactions are the focus of the present investigation, 
where we study the 
global dynamics of a one parameter family of
three dimensional vector fields whose main features are
stable and saddle type equilibrium solutions and a   
periodic orbit with a complex conjugate pair of Floquet exponents.
The frequency of the periodic orbit together with the frequency of the 
complex exponent constitute two 
competing natural modes of oscillation.  Tension between these
internal frequencies gives rise to a number of interesting dynamical phenomena.
In particular the system admits a Neimark-Sacker bifurcation, where the real part of 
the complex conjugate Floquet exponents
crosses the imaginary axis as the parameter is changed \cite{MR0132256,MR2615427}.
The loss of stability of the periodic orbit
triggers the appearance of a smooth attracting invariant torus 
supporting quasiperiodic motions. 
Global bifurcations of the torus lead to resonant motions 
and eventually to the appearance of a chaotic attractor.

The local theory describing the 
appearance, evolution, and disappearance of invariant tori in 
dissipative multi-frequency systems is well developed
and we refer the reader to  the works of  
\cite{MR1005055,MR1115870,MR709899,
MR834186, MR906312,MR3095277,MR3062760,
MR3713933} on dissipative 
dynamics, the related work of 
\cite{MR783349,MR784530,MR932134,MR1285950}
on area and volume preserving systems, and 
to the numerical studies of 
\cite{MR709899,MR701669,
MR845031,MR2241302,MR3309008,MR3713932,
MR821035}.  Global questions about the dynamics of  
dissipative systems with attracting tori lead to 
difficult analytical and computational problems.
While many important theoretical questions have
by now been settled -- see for example
\cite{MR0462175,MR1733750,MR1237641,MR2471925}
and the references therein --
there remains much to be learned from careful qualitative 
studies of important special cases.

While many of the canonical examples of dynamical systems
theory come from specific physical or engineering applications, another 
source of compelling problems is to study the  
normal form of an interesting bifurcation.
Such systems caricature the universal features of
an entire class of problems, and
this is precisely the setting of the present paper. 
We study, from the numerical point of view,
a model derived from the normal form unfolding
the cusp-Hopf bifurcation.
This system, which is described in detail in Section \ref{sec:LangfordSystem},
was first introduced in \cite{MR821035} 
and is referred to hereafter as \textit{the Langford system}.
As already mentioned in the opening paragraph, a
main feature is that model undergoes a supercritical Neimark-Sacker bifurcation
resulting in the appearance of a smooth attracting invariant torus. 
We provide detailed computations of the torus, monitoring 
it as its dynamics change from quasi-periodic to resonant 
-- and as it changes from a $C^k$ to a $C^0$ invariant 
manifold -- before finally breaking up in a global bifurcation  
resulting in the appearance of a chaotic attractor.

In addition to undertaking a detailed description of the attracting 
invariant torus, the present work aims also to describe 
the dynamics nearby.  We are especially interested in any 
dramatic changes in the organization of the phase space as
the bifurcation parameter is varied.  Such changes may be 
triggered by either \textit{local} or \textit{global} bifurcations.  
More precisely we have the following distinction.

\begin{definition} We say that a bifurcation is {\it local }if it occurs due to a change in 
linear stability of an invariant object. \end{definition} 

\begin{definition} We say that a bifurcation is {\it global} if it is triggered by the formation 
of tangencies between invariant manifolds.\end{definition} 

\noindent In the present work we mainly observe
local bifurcations of equilibrium and periodic solutions --
and global bifurcations where the invariant manifolds do 
not intersect at all prior to, and intersect transversally
after the global bifurcation.

The discussion just presented makes it clear that the goals of the present work 
require careful examination of the embeddings of 
some hyperbolic invariant objects like stable/unstable
manifolds of equilibrium and periodic orbits. 
Much of the analysis is simplified 
by considering an appropriate surface of section, as this reduces
the invariant torus
and the stable/unstable manifolds of periodic orbits to one dimensional curves.
Embeddings of stable/unstable manifolds attached to equilibrium solutions
on the other hand are often difficult to characterize in 
a fixed section, and studying their structure is more delicate.
We employ the parameterization method of
\cite{Cabre1,Cabre2,Cabre3} to compute high order 
representations of the two dimensional local stable/unstable manifolds 
in the full three dimensional phase space.
The parameterization method is a functional 
analytic framework for studying invariant manifolds and in particular 
provides a natural notion of a-posteriori error analysis.  The local representations
 obtained using the parameterization method are
extended using standard adaptive numerical integration schemes.

The detailed numerical calculations performed in the main body of the 
paper provide insights into the dynamics of the system which are 
summarized in Section \ref{sec:conclusions}, 
and which give a coarse qualitative description of the global 
dynamics as a function of the bifurcation parameter.  
Since the Langford system is derived from a normal form, 
it is reasonable to expect qualitatively 
similar dynamics in an appropriately restricted region for any  
system undergoing the sequence of bifurcations unfolded 
by this vector field. 
Moreover, the approach of using the parameterization method
in conjunction with geometric analysis 
in Poincar\'{e} sections could be applied to the study of 
a wide variety of dynamical systems.

The remainder of the paper is organized as follows. In the next two subsections we first describe the three dimensional 
model under consideration, and then 
discuss briefly some related literature.  
In Section \ref{sec:parmMethod} we review the main ideas of the 
parameterization method for an equilibrium solution and apply them to 
the Langford system.
In Section \ref{sec:global} we
study the Neimark-Sacker bifurcation and the resulting attracting invariant
torus in an appropriate Poincar\'{e} section. We provide numerical evidence for 
a global bifurcation from a quasi-periodic torus to a resonant one, and 
 for a second global bifurcation which destroys the torus and appears to 
 create a chaotic attractor.
In Section \ref{sec:globalD} we study the invariant manifolds of the equilibrium solutions 
before and after the Neimark-Sacker bifurcation, with an emphasis on the 
omega limit sets of two dimensional unstable manifolds and on the role of the 
manifolds as separatrices.
We also study their role in further global bifurcations. 
We conclude the paper in Section \ref{sec:conclusions} with a 
summary of our observations about the global dynamics of the system and a
few further conclusions and observations.

\subsection{The Langford system} \label{sec:LangfordSystem}
We study the dynamical system generated by the
3D vector field 
\begin{equation}\label{eq:1}
f(x,y,z) = \begin{pmatrix} 
  (z - \beta)x - \delta y \\ 
  \delta x + (z - \beta)y \\
  \tau + \alpha z - \displaystyle\frac{z^3}{3} - (x^2 + y^2)(1+\varepsilon z) + \zeta z x^3 \\
\end{pmatrix},
\end{equation}
where $\varepsilon = 0.25$, $\tau = 0.6$, $\delta = 3.5$, $\beta = 0.7$, $\zeta = 0.1$, and with
$\alpha > 0$ treated as a bifurcation parameter.
The system was derived by Langford in \cite{MR821035} by
truncating to second order the 
normal form unfolding a simultaneous
Hopf/cusp bifurcation.  A third order term is then 
added to the vector field, breaking
the axial symmetry of the second order truncation.   This symmetry breaking 
is important for describing the dynamics following a generic bifurcation.  
Since the Hopf bifurcation creates a periodic orbit, and the cusp bifurcation 
creates three nearby equilibrium solutions, 
interesting interactions between these states are to be expected.

\begin{figure}[t!]
\centering
\subfloat{\includegraphics[width = 3in, height = 2.5in]{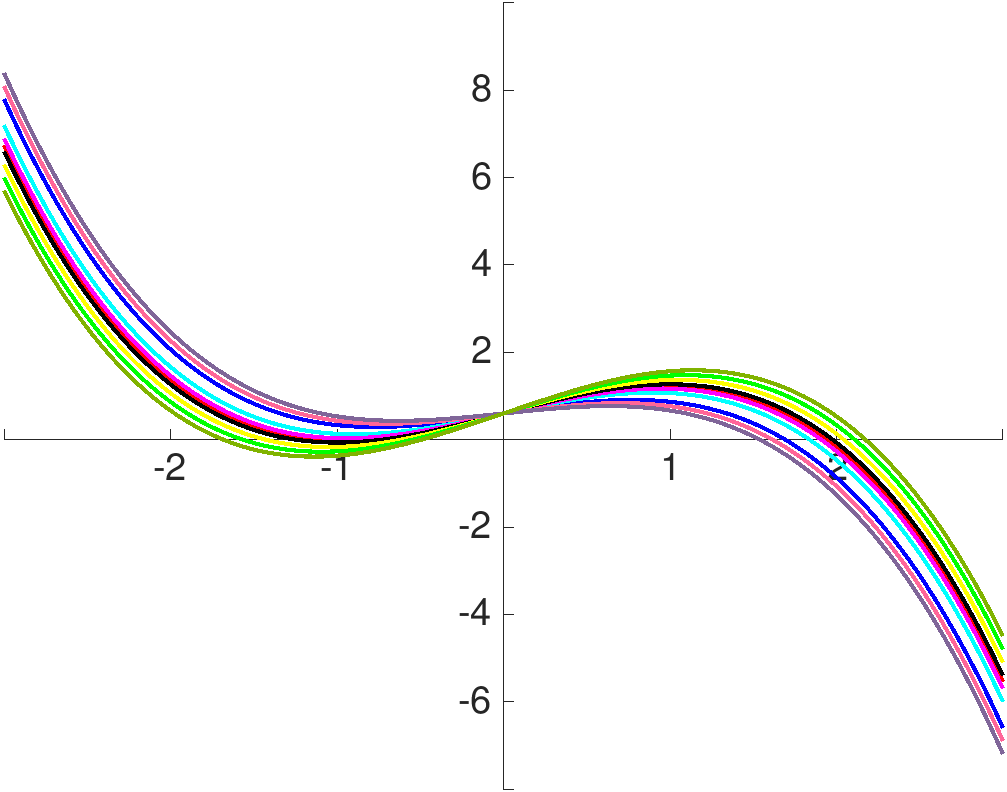}}
\caption{Graph of $g(z) = \tau + \alpha z - \displaystyle\frac{z^3}{3}$ for different parameter values of $\alpha$ and 
fixed $\tau = 0.6$.}
 \label{fig:gz}
\end{figure}

We begin with some elementary observations which inform the 
numerical study to follow.  
Note that the $z$-axis is an invariant sub-system as
$x = y = 0$ implies that $x' = y' = 0$. 
The dynamics on the $z$-axis are governed by the scalar differential 
equation 
\[
z' = \tau + \alpha z - \displaystyle\frac{z^3}{3} =: g(z).
\]
The function $g(z)$ is illustrated in Figure \ref{fig:gz}, and since $\tau > 0$, 
$g$ has one, two, or three zeros depending on the parameter $\alpha$.  
Moreover, equilibria of $f$ occur at $(0, 0, z_*)$ where $z_*$ is a zero of $g$.
Observe that for large positive $z$, $z' < 0$. While for large 
negative $z$, $z' > 0$.  That is, the field tends to diminish the $z$ value of a phase 
point whose $z$ value happens to be large. 

For all $\alpha \in \mathbb{R}$ \correction{comment:2}{Equation \eqref{eq:1} 
has at exactly one equilibrium solution with $x=y=0$ and $z>0$, which} we 
denote by $p_0 \in \mathbb{R}^3$.  This equilibrium has one stable eigenvalue, whose 
eigenvector coincides with the $z$-axis.
The remaining eigenvalues are complex conjugate unstable.  
At $\alpha \approx 0.9321697517861$ there is a saddle node bifurcation giving rise to 
a new pair of equilibrium points $p_1, p_2 \in \mathbb{R}^3$.  
These equilibria persist for all larger values of $\alpha$.
One of the equilibrium points 
appearing out of the saddle node bifurcation is fully stable, with three eigenvalues
having negative real parts, and we denote it by $p_2 \in \mathbb{R}^3$.
The other new equilibrium, which we denote by $p_1 \in \mathbb{R}^3$, is a saddle-focus with 
a complex conjugate pair of   
stable eigenvalues and one real unstable eigenvalue.  The unstable eigenvector again 
coincides with the $z$-axis.  Indeed, since the $z$-axis is invariant, 
the stable manifold of $p_0$ and the unstable manifold of $p_1$ 
\correction{comment:3}{coincide}, and are
contained in the $z$-axis.  This intersection is not transverse, and is rather forced by 
a rotational symmetry of the problem.

Now consider the plane $z=\beta$, and note that when the field is
projected onto this plane the nonlinear terms vanish from the first two 
components giving a pure rotation.
The plane is however not invariant, as $z'$ does not vanish there.
Nevertheless there is a periodic orbit $\gamma$ near the $z = \beta$ plane.
This periodic orbit, and the invariant $z$ axis organize the dynamics of 
the system.
The vector field along with the periodic
 orbit and the dynamics on the $z$-axis are illustrated
in the left frame of Figure \ref{fig:qualitity}.

As we will see below, the periodic orbit $\gamma$ has a pair of complex conjugate
Floquet exponents, hence solutions of the differential equation tend to 
circulate around $\gamma$.  
The orbit may be either attracting or repelling depending on the value of $\alpha$.
This circulation about the periodic orbit
is a dominant feature of the dynamics.

\begin{figure}[t!]
\subfloat[$\alpha=0.8$]{\includegraphics[width = 3in]{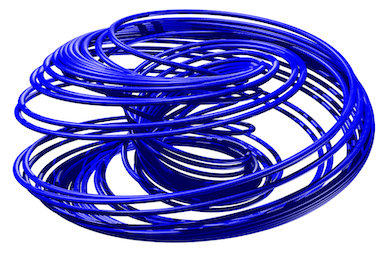}}\quad 
\subfloat[$\alpha=0.95$ ]{\includegraphics[width = 3in]{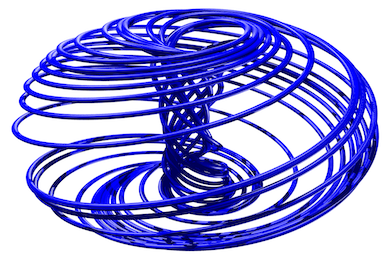}}\\
\caption{Direct simulation:  For many values of the bifurcation
parameter $\alpha$ the system 
appears to have an attractor with torus-like dynamics (product of two circles).
This is caused by circulation due to the complex conjugate Floquet exponents 
of the periodic orbit, and generates a kind of ``vortex''.}
\label{fig:attractors}
\end{figure}

Further insight into the dynamics is obtained by numerically 
integrating some trajectories (phase space sampling), 
as was done in the work of Langford \cite{MR821035}.  We provide, for 
the sake of completeness, the results of a few such simulations.  
The results illustrated in Figure \ref{fig:attractors} make clear
the typical behavior of the system,
and suggest the existence of a  ``torus-like'' attractor.  
Simulations were run for roughly one hundred time units.
The periodic orbit $\gamma$ runs through the 
center of the torus but is, as we will see, repelling for these parameter values. The  
saddle focus points $p_0$ and $p_1$ are at the top and bottom of the torus.

\begin{figure}[t!]
\centering
\subfloat{\includegraphics[width = 2.5in,valign=c]{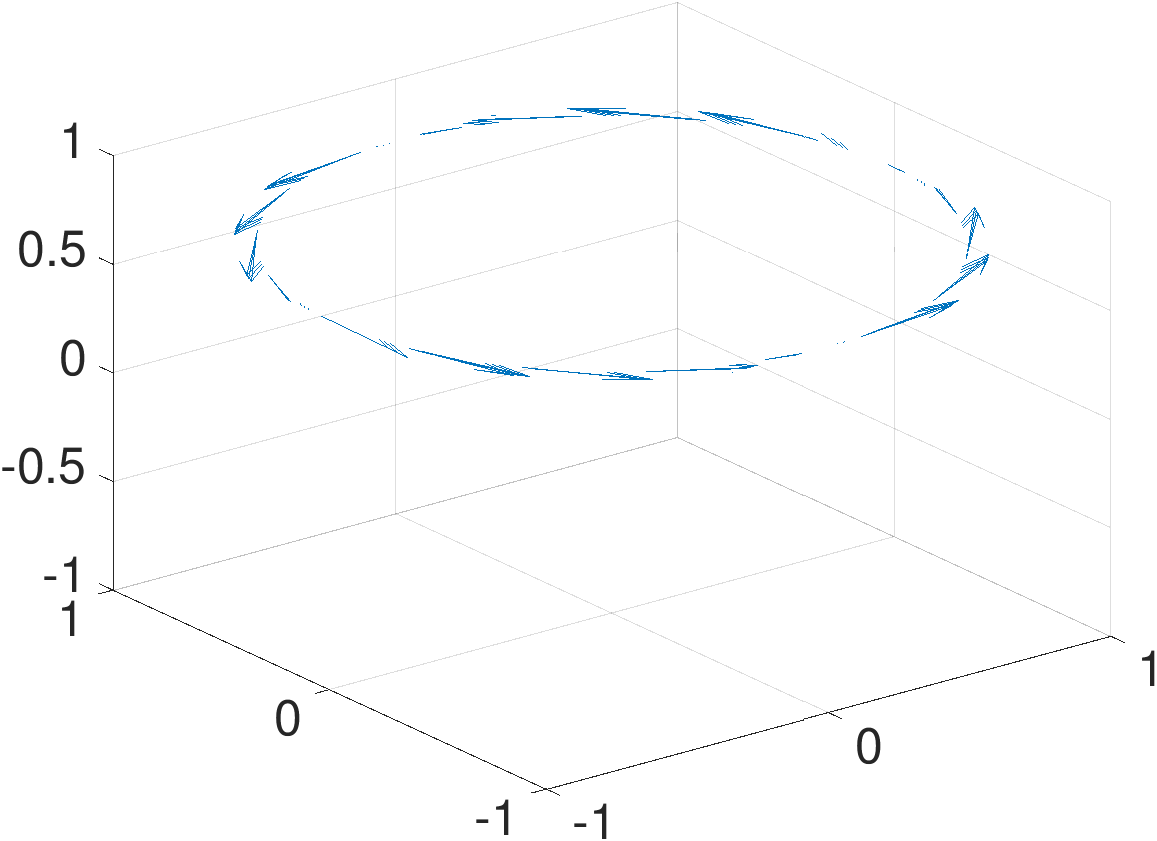}}\quad 
\subfloat{\includegraphics[width = 3.5in,valign=c]{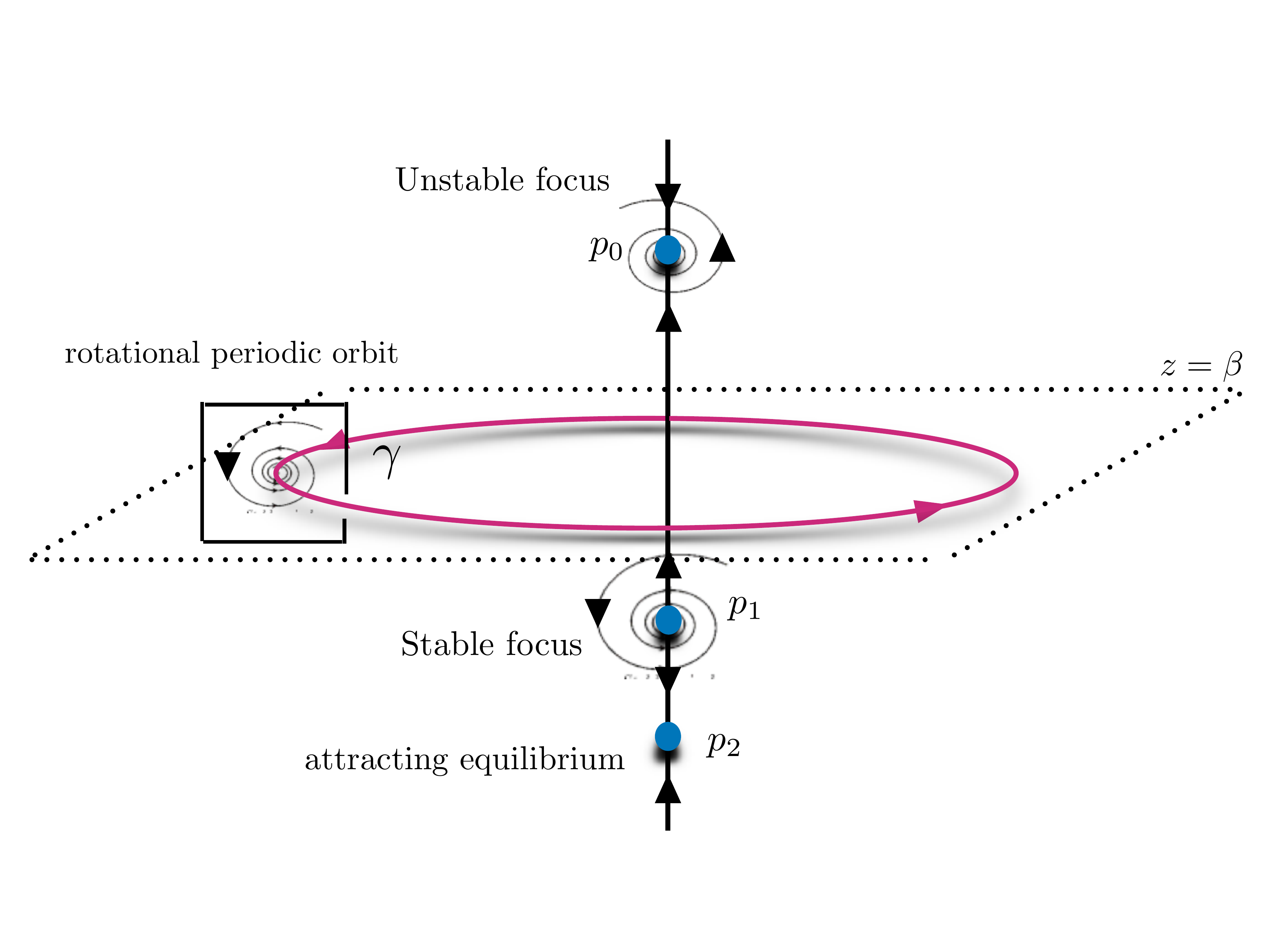}}\\
\caption{Phase space geography: the main features
of the system are the invariance of the $z$-axis, the rotation in the $z = \beta$
plane leading to a periodic orbit, and the unstable saddle focus at $p_0$.
The periodic orbit $\gamma$ is located near (but does not sit on) the $z = \beta$ plane.
The periodic orbit has complex conjugate Floquet multipliers which are 
stable for small $\alpha$ but which later cross the unit circle,
loosing stability in a Neimark-Sacker bifurcation.
For some $\alpha$ values there are an additional pair of equilibria
$p_1$-- stable focus and $p_2$ -- attracting point.  
This situation is illustrated in the schematic on the right.
Left: the phase portrait of the 
vector field along the periodic orbit.}
\label{fig:qualitity}
\end{figure}

\subsection{Some remarks on the literature}\label{ref:literature}

Roughly speaking, the dynamics described above 
suggests the system as a
toy model for dissipative vortex dynamics, or for a rotating viscus fluid.
There is a rich literature on the dynamics of vortex bubbles, and the 
interested reader might consult the works of 
\cite{MR1285950, Stone05, MR1704974, MR2481277, MR2259296}
for a more thorough discussion of the literature.
We remark that the torus bifurcations seen in the 
Langford system are similar to those seen in the piecewise linear 
electronic circuit of \cite{MR880159}, the commodity distribution 
model of \cite{MR1488520}, and the mechanical oscillators 
of \cite{MR3435117, MR3279518} to name only a few.
The appearance and destruction of invariant tori, 
as well as resonance phenomena and routes to chaos 
 are discussed much more generally in 
\cite{MR709899, Vadim} and the references found therein.

One further remark is in order.  The system given by 
Equation \eqref{eq:1} has been called  
\textit{the Aizawa system} by some researchers, and 
is the subject of some other recent work on
visualization.  For example researchers interested in 
computer animation \cite{aizawaVideo},
three dimensional printing \cite{bridges2018:491},
and even in graphical arts \cite{chaoticAtmsopheres}
have made interesting studies and use this name for the equations.  
This nomenclature seems to be a misnomer, as the equations 
do not appear in the works of Yoji Aizawa, and 
a more appropriate name for Equation \eqref{eq:1} would seem to be
the Langford system, due to the fact that -- as already mentioned 
above -- the system was proposed in \cite{MR821035}.

\section{Review of the parameterization method} \label{sec:parmMethod}
The parameterization method is a general functional analytic framework for studying invariant 
manifolds of discrete and continuous time dynamical systems, first developed in  
\cite{Cabre1, Cabre2, Cabre3} in the context of stable/unstable manifolds
attached to fixed points of nonlinear mappings on Banach spaces, and later 
extended in  \cite{Haro1, Haro2, MR2299977} for studying whiskered tori.
There is a thriving literature devoted to 
computational applications of the parameterization method, 
and the interested reader may want to consult
\cite{MR3713933, MR3309008, MR3713932, MR3118249, MR2551254, maximePOman, maximeJPMe, parmChristian, 
fastSlow, manifoldPaper1, jorgeMePerParm, jayChrisParmDDE, 
MR3783519, MR3797119, MR3705136}, 
though the list is far from being exhaustive.  A much more complete discussion is 
found in the book \cite{Cana}.

This section provides a practical overview of the parameterization method
with a strong emphasis on numerical aspects utilized in the sequel.  
The discussion focuses on analytic vector fields, and requisite formal series 
calculations are carried out for the specific example of Equation \eqref{eq:1}.
Since this material is not completely standard outside a certain 
circle of practitioners, it is included primarily so that the present work 
\correction{comment:4}{stands} alone for a broad readership.  
The reader either already familiar with or uninterested in these developments 
is encouraged to skip ahead to Section \ref{sec:global}, referring back to this section only 
as needed.

The Langford system admits equilibrium solutions with complex eigenvalues,
so that it is best to present the entire theory for complex vector fields.
Later we explain how to recover parameterizations of real invariant manifolds 
associated with complex conjugate eigenvalues of real vector fields. 
So, let  $f: \mathbb{C}^k \to \mathbb{C}^k$ be an analytic vector field 
and $\hat p \in  \mathbb{C}^k$ have $f(\hat p) = 0$ so that 
$x(t) = \hat p$ is an equilibrium solution of the differential equation $x' = f(x)$.
Assume for the sake of simplicity that $Df(\hat p)$ is diagonalizable over $\mathbb{C}$
having $k_s$ stable (and $k_u$ unstable) eigenvalues of multiplicity one.  
We do not necessarily assume that $k_s + k_u = k$, that is we do not rule out the 
possibility of some center directions at $\hat p$ (though this situation will not occur in the 
present work).

Label the stable eigenvalues as $\lambda_1^s, \cdots, \lambda_{k_s}^s$
and the unstable ones as $\lambda_1^u, \cdots, \lambda_{k_u}^u$ 
and order them according to the convention that 
\[\textrm{real}(\lambda_1^s) 
\leq \cdots \leq \textrm{real}(\lambda_{k_s}^s)< 0 < \textrm{real}(\lambda_1^u) 
\leq \cdots \leq \textrm{real}(\lambda_{k_u}^u).
\]
Since $Df(\hat p)$ is diagonalizable
there are linearly independent eigenvectors $\xi_1^u, \cdots, \xi_{k_u}^u \in \mathbb{C}^k$ and 
$\xi_1^s, \cdots, \xi_{k_s}^s \in  \mathbb{C}^k$ associated with the unstable and stable
eigenvalues respectively.

\begin{remark}
The assumption that $Df(\hat p)$ is diagonalizable is made only for the sake 
of convenience.  See 
\cite{Cabre1} for a much more general theoretical setup.  See also
\cite{parmChristian} for a complete description of the functional analytic 
set up and examples of the numerical implementation 
when there are repeated eigenvalues.  Nevertheless, the assumption holds
in the examples considered throughout the present work.  
\end{remark}

\begin{figure}[t!]
\centering
\includegraphics[width=0.6\textwidth]{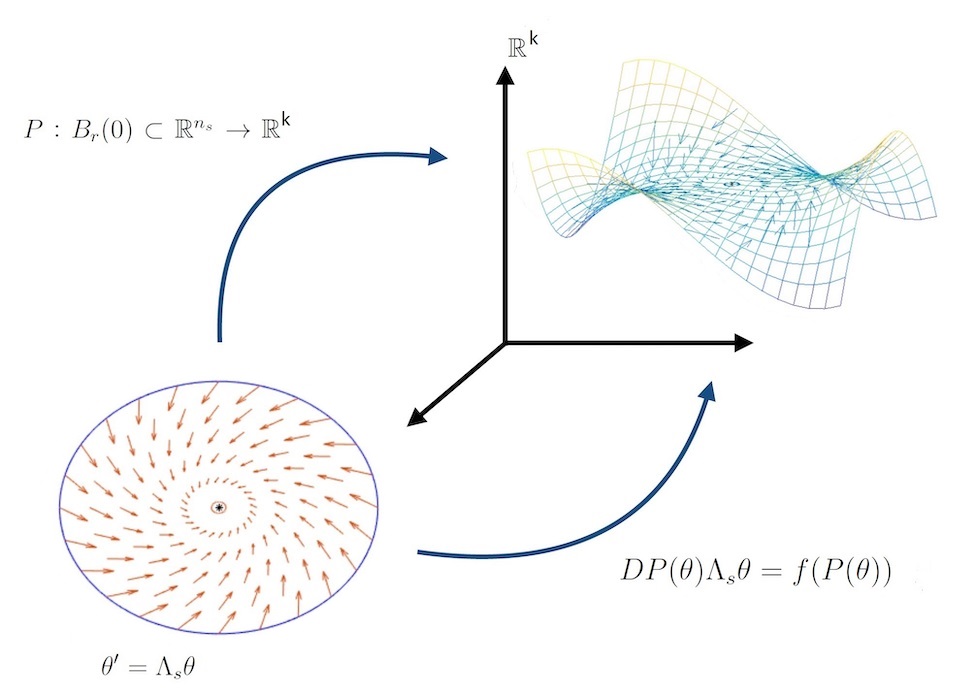}
\caption{Geometry of the parameterization method: the push forward of the linear 
vector field $\Lambda_s \theta$ by $P$ is equal to the given vector field $f$ 
restricted to the image of $P$.  Then the dynamics on the image of $P$ are conjugate to 
the linear dynamics generated by $\Lambda_s \theta$.} \label{fig:pushForward}
\end{figure}

\subsection{Invariance equation}
Given the setup introduced in the previous section we are
interested in computing an accurate representation of the $k_s$ dimensional 
local stable manifold attached to $\hat p$.
The parameterization method seeks 
a smooth surjective map $P$ 
satisfying the first order system of partial differential equations
\begin{equation}\label{eq:invEq}
\lambda^s_1 \theta_1\frac{\partial}{\partial \theta_1} P(\theta_1, \ldots, \theta_{k_s}) 
+ \ldots +\lambda^s_{k_s} \theta_{k_s} 
\frac{\partial}{\partial \theta_{k_s}} P(\theta_1, \ldots, \theta_{k_s}) 
= f(P(\theta_1, \ldots, \theta_{k_s})),
\end{equation} 
for $\theta = (\theta_1, \ldots, \theta_{k_s}) \in \mathbb{C}^{k_s}$, and 
subject to the first order constraints 
\begin{equation} \label{eq:firstOrder}
P(0, \ldots, 0) = \hat p, 
\quad \quad \quad \mbox{and} \quad \quad \quad
\frac{\partial}{\partial \theta_j}P = \xi^s_j, \quad \quad 1 \leq j \leq k_s.
\end{equation}
Equation \eqref{eq:invEq} is referred to as the invariance equation for $P$. 
A map $P$ solving Equation \eqref{eq:invEq}  subject to the first order
constraints of Equation \eqref{eq:firstOrder} is 
 a parameterization of the local stable manifold, as we explain below.
Making the obvious adjustments for the unstable eigenvalues/eigenvectors
leads to a parameterization method for the $k_u$ dimensional unstable manifold.

To explain the meaning of Equation \eqref{eq:invEq}, 
let
\[
\Lambda_s = \left(
\begin{array}{ccc}
\lambda_1 & \ldots & 0 \\
\vdots & \ddots & \vdots \\
0& \ldots & \lambda_{k_s}
\end{array}
\right),
\]
so that the invariance equation becomes
\[
DP(\theta) \Lambda_s \theta = f(P(\theta)).
\]
In the language of differential geometry, this equation says that the push forward 
by $P$ of the linear 
vector field $\Lambda_s \theta$ is equal to the vector field $f$ restricted
to the image of $P$.  Where the vector fields are equal  
they generate the same dynamics.  But the dynamics generated 
by $\Lambda_s \theta$
are completely understood: all orbits converge exponentially to the origin.  
It follows that all orbits on the image of $P$ converge  
to $\hat p$.  Since the image of $P$ is a smooth $k_s$ dimensional disk, 
it is a local stable manifold for $\hat p$. 
The situation is illustrated in Figure \ref{fig:pushForward}.

\begin{figure}[t!]
\centering
\includegraphics[width=0.7\textwidth]{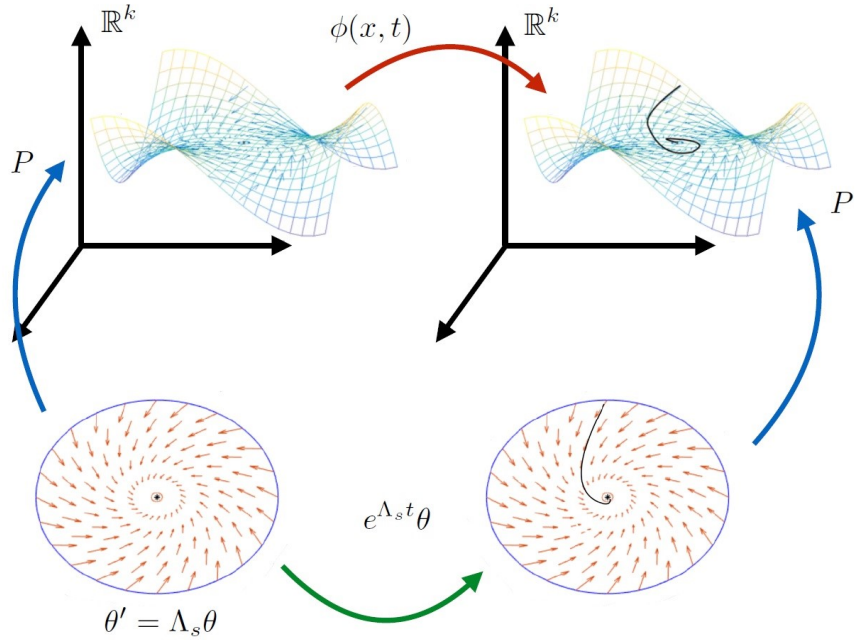}
\caption{Flow conjugacy: a mapping $P$ satisfying the invariance 
equation \eqref{eq:invEq} has that the diagram above commutes. } \label{fig:conjugacy}
\end{figure}

The observation is made more precise as follows.
Denote by $\phi \colon \mathbb{C}^k \times \mathbb{C} \to \mathbb{C}^k$
the flow generated by $f$.
The flow generated by $\Lambda_s \theta$ is given explicitly by 
\[
L(\theta, t) = e^{\Lambda_s t} \theta.
\]
One checks that $P$ satisfies Equation \eqref{eq:invEq}
if and only if
\begin{equation} \label{eq:flowConj}
\phi\left(P(\theta), t\right) = P\left(e^{\Lambda_s t} \theta\right), 
\end{equation}
for all $t \geq 0$.  This flow conjugacy is illustrated in 
Figure \ref{fig:conjugacy}.
Elementary proofs of these claims are found in any of the references
\cite{fastSlow, Cana, AMS}.  Moreover, replacing the stable by the unstable 
eigenvalues and eigenvectors in the discussion above and reversing time, 
the entire discussion carries through for the unstable manifold.

\begin{remark}[Real analytic vector fields and manifolds] \label{rem:real1} 
If $f \colon \mathbb{R}^k \to \mathbb{R}^k$ is a real analytic vector field 
with a real equilibrium $\hat p \in \mathbb{R}^k$ then the 
discussion above applies to an analytic extension of the vector field in a 
neighborhood of $\hat p$.  In this case any complex eigenvalues of $Df(\hat p)$
appear in complex conjugate pair, and the associated eigenvectors
can be taken complex conjugate.   
We look for a solution $P$ of Equation \eqref{eq:invEq}
taking real values on complex conjugate variables.  
This condition imposes a symmetry on the Taylor coefficients of the parameterization 
$P$, as illustrated explicitly in the examples below.
\end{remark}

\subsection{Formal series solution of Equation \eqref{eq:invEq} for the Langford system} 
\label{sec:formalSeries}
In this section we further restrict to the case of interest in the present work, 
where $\hat p \in \mathbb{C}^3$ and $\lambda_1, \lambda_2 \in \mathbb{C}$ 
are a pair of stable (or unstable) complex conjugate 
eigenvalues and $\lambda_3$ 
has the opposite stability.  Let $\xi_1, \xi_2 \in \mathbb{C}^3$ be an 
associated pair of linearly independent complex conjugate eigenvectors.
Since the field is analytic,  
we look for an analytic parameterization
\[
P(\theta_1, \theta_2) = \sum_{m=0}^{\infty} \sum_{n=0}^{\infty} \begin{pmatrix} 
  p_{mn}^1\\ 
  p_{mn}^2\\
  p_{mn}^3\\
\end{pmatrix} \theta_1^m \theta_2^n, 
\]
satisfying Equation \eqref{eq:invEq}, which in this case is reduced to 
\[
\lambda_1 \theta_1 \frac{\partial}{\partial \theta_1} P(\theta_1, \theta_2)
+ \lambda_2 \theta_2 \frac{\partial}{\partial \theta_2} P(\theta_1, \theta_2) 
= f(P(\theta_1, \theta_2)),
\]
where $f \colon \mathbb{C}^3 \to \mathbb{C}^3$ is the Langford vector field given in 
Equation \eqref{eq:1}.
Here $p_{mn}^j \in \mathbb{C}$ for all $j = 1,2,3$.
Imposing the linear constraints of Equation \eqref{eq:firstOrder} 
gives that $p_{00} = \hat{p}$, $p_{10} = \xi_1$ and $p_{01} = \xi_2$.

Now we would like to expand Equation \eqref{eq:invEq} in terms of the power series.
The left hand side of Equation \eqref{eq:invEq} is 
\[
\lambda_1 \theta_1 \frac{\partial}{\partial \theta_1} P(\theta_1, \theta_2)
+ \lambda_2 \theta_2 \frac{\partial}{\partial \theta_2} P(\theta_1, \theta_2) = 
\sum_{m=0}^\infty \sum_{n=0}^\infty (m \lambda_1 + n \lambda_2) p_{mn} \theta_1^m \theta_2^n,
\]
on the level of power series. 
To expand the right hand side we begin by writing 
\[
P_j(\theta_1, \theta_2) = \sum_{m=0}^\infty \sum_{n=0}^\infty p_{mn}^j \theta_1^m \theta_2^n,
\]
for $j = 1,2,3$ to denote the component power series.  
The field contains the nonlinear terms
$zx$, $zy$, $z^3$, $x^2 z$, $y^2 z$, and $z x^3$ (see again Equation \eqref{eq:1}). 
 Computing the power series 
for $f \circ P$ requires expanding these monomials of components of $P$, 
which is accomplished using Cauchy products.
For example the coefficients of $P_3 \cdot P_1$ are 
\[
(p^3 * p^1)_{mn} =  \sum_{j = 0}^m \sum_{k=0}^n p_{(m-j) (n-k)}^3 p_{jk}^1,
\]
while the coefficients of $P_3 \cdot P_1^3$ are 
\[
(p^3 * p^1 * p^1 * p^1)_{mn} = 
\sum_{i_1 = 0}^m \sum_{i_2 =0}^{i_1} \sum_{i_3=0}^{i_2} 
\sum_{k_1=0}^n \sum_{k_2 = 0}^{k_1} \sum_{k_3=0}^{k_2}
 p_{(m- i_1)(n-k_1)}^3 p_{(i_1- i_2)(k_1- k_2)}^1 p_{(i_2-i_3)(k_2 - k_3)}^1 p_{i_3 k_3}^1.
\]
Other products are similar.

\correction{comment:24}{Substituting} these power series expansions into the invariance equation \eqref{eq:invEq} 
and matching like powers of $\theta_1$ and $\theta_2$ leads to 
\begin{equation} \label{eq:scratchEq1}
(m \lambda_1 + n \lambda_1) 
\left[
\begin{array}{c}
p_{mn}^1 \\
p_{mn}^2 \\
p_{mn}^3
\end{array}
\right] = 
\end{equation}
\[
{\tiny
\left[
\begin{array}{c}
(p^3 * p^1)_{mn} - \beta p_{mn}^1 - \delta p_{mn}^2 \\
(p^3 * p^2)_{mn} - \beta p_{mn}^2 + \delta p_{mn}^1 \\
\alpha p_{mn}^3 - \displaystyle\frac{1}{3} (p^3 * p^3 * p^3)_{mn} - 
(p^1 * p^1)_{mn} - (p^2 * p^2)_{mn} - \varepsilon (p^1 * p^1 * p^3)_{mn}
- \varepsilon (p^2 * p^2 * p^3)_{mn} + \zeta (p^1 * p^3 * p^3 * p^3)_{mn}
\end{array}
\right] ,
}
\]
for $m + n \geq 2$.  To isolate terms of order $(m,n)$ 
consider that 
\begin{equation}\label{eq:homDeriv1}
(p^3 * p^1)_{mn} = p_{00}^3 p_{mn}^1 + p_{00}^1 p_{mn}^3 + (p^3 \hat * p^1)_{mn},
\end{equation}
where 
\[
(p^3 \hat * p^1)_{mn} =  \sum_{j = 0}^m \sum_{k=0}^n  \hat \delta_{jk}^{mn} p_{(m-j)(n-k)}^3 p_{jk}^1,
\]
and 
\[
\hat \delta_{jk}^{mn} = \begin{cases}
0 & \mbox{if } j = k = 0 \\
0 & \mbox{if } j = m \mbox{ and } k = n \\
1 & \mbox{otherwise}
\end{cases}.
\]
The point here is that $(p^3 \hat * q^1)_{mn}$ is precisely the sum left when terms 
containing $p_{mn}$ are extracted from the Cauchy product.

This expression is directly related to the derivative of $f$. 
To see this, let
\[
g(x, z) = xz, 
\]
and note that Equation \eqref{eq:homDeriv1} becomes
\[
(g \circ P)_{mn} = \nabla g(p_{00}^1, p_{00}^3) \left[
\begin{array}{c}
p_{mn}^1 \\
p_{mn}^3
\end{array}
\right]
 + (p^1 \hat * p^3)_{mn}. 
\]
Using this notation
the first component of Equation \eqref{eq:scratchEq1} is
\[
(m \lambda_1 + n \lambda_2) p_{mn}^1  =
\nabla g(p_{00}^1, p_{00}^3) \left[
\begin{array}{c}
p_{mn}^1 \\
p_{mn}^3
\end{array}
\right]
 + (p^1 \hat * p^3)_{mn}  - \beta p_{mn}^1  - \delta p_{mn}^2.
\]
Isolating terms of order $(m,n)$ on the left and lower order terms 
on the right gives
\[
\nabla g(p_{00}^1, p_{00}^3) \left[
\begin{array}{c}
p_{mn}^1 \\
p_{mn}^3
\end{array}
\right]
 - \beta p_{mn}^1  -  \delta p_{mn}^2
- (m \lambda_1 + n \lambda_2) p_{mn}^1
 = -   (p^1 \hat * p^3)_{mn},
\]
which is linear in $p_{mn}^1$.
Comparing the right hand side in the equation above 
with the vector field $f$, and recalling that 
$\hat p = p_{00}$, we see that  
\[
\nabla g(p_{00}^1, p_{00}^3) \left[
\begin{array}{c}
p_{mn}^1 \\
p_{mn}^3
\end{array}
\right]
 - \beta p_{mn}^1  -  \delta p_{mn}^2
 = \nabla f_1(\hat p)  \left[
\begin{array}{c}
p_{mn}^1 \\
p_{mn}^3
\end{array}
\right].
\]
Combining the equation above with 
a nearly identical computation for the 
second component, and a somewhat lengthier
computation for the third component, and 
noting that 
\[
Df(\hat p) = 
\left[
\begin{array}{c}
\nabla f_1(\hat p) \\
\nabla f_2(\hat p) \\
\nabla f_3 (\hat p)
\end{array}
\right],
\]
we obtain the expansion
\[
(f \circ P)_{mn} = 
Df(\hat p) p_{mn} + 
\]
\[
{\tiny
\left[
\begin{array}{c}
(p^3 \hat * p^1)_{mn}  \\
(p^3 \hat * p^2)_{mn}  \\
 - \displaystyle\frac{1}{3} (p^3 \hat * p^3 \hat * p^3)_{mn} - 
(p^1 \hat * p^1)_{mn} - (p^2 \hat * p^2)_{mn} - \varepsilon (p^1 \hat * p^1 \hat * p^3)_{mn}
- \varepsilon (p^2 \hat * p^2 \hat * p^3)_mn + \zeta (p^1 \hat * p^3 \hat * p^3 \hat * p^3)_{mn}
\end{array}
\right].
}
\]
\correction{comment:24}{Substituting} this expansion into Equation \eqref{eq:scratchEq1}
gives
\[
(m \lambda_1 + n \lambda_2) p_{mn} = 
Df(\hat p) p_{mn} + 
\]
\[
{\tiny
\left[
\begin{array}{c}
(p^3 \hat * p^1)_{mn}  \\
(p^3 \hat * p^2)_{mn}  \\
 - \displaystyle\frac{1}{3} (p^3 \hat * p^3 \hat * p^3)_{mn} - 
(p^1 \hat * p^1)_{mn} - (p^2 \hat * p^2)_{mn} - \varepsilon (p^1 \hat * p^1 \hat * p^3)_{mn}
- \varepsilon (p^2 \hat * p^2 \hat * p^3)_mn + \zeta (p^1 \hat * p^3 \hat * p^3 \hat * p^3)_{mn}
\end{array}
\right],
}
\] 
and by isolating terms of order $(m,n)$ on the left
we obtain the 
linear \textit{homological equations}  
\begin{equation} \label{eq:homEq} 
[Df(\hat p) - (m \lambda_1 + n \lambda_2) \mbox{Id}] p_{mn} = s_{mn}, 
\end{equation}
for $p_{mn}$, where 
\[
s_{mn} = 
\left(
\begin{array}{c}
s_{mn}^1 \\
s_{mn}^2 \\
s_{mn}^3
\end{array}
\right),
\]
with 
\[
s_{mn}^1 =   - (p^3 \hat * p^1)_{mn},
\]
\[
s_{mn}^2 =  - (p^3 \hat * p^2)_{mn},
\]
and 
\[
s_{mn}^3 = 
\]
\[
\frac{1}{3} (p^3 \hat * p^3 \hat * p^3)_{mn} +
(p^1 \hat * p^1)_{mn} + (p^2 \hat * p^2)_{mn}
+ \varepsilon (p^1 \hat * p^1 \hat * p^3) + \varepsilon (p^2 \hat * p^2 \hat * p^3)
- \zeta (p^1 \hat * p^1 \hat * p^1 \hat * p^3)_{mn}.
\]
We make the following observations:
\begin{itemize}
\item While our derivation of Equation \eqref{eq:homEq} is particular to the Langford system
of Equation \eqref{eq:1}, 
we remark that the form of the homological equations is always the same.  Only the 
right hand side depends on the particular nonlinearity of the given system. 
\item The matrix acting on $p_{mn}$ is
the characteristic matrix for the differential at $\hat p$.  Then the 
equation is uniquely solvable at order $(m,n)$ if $m \lambda_1 + n \lambda_2$
is not an eigenvalue.    
\item Since $\lambda_3$ has the opposite stability of $\lambda_1, \lambda_2$,
we obtain the \textit{non-resonance condition}
\[
m \lambda_1 + n \lambda_2 \neq \lambda_j, \quad \quad \quad j = 1,2.
\]
If the non-resonance conditions are satisfied for all $m, n \in \mathbb{N}$
with $m + n \geq 2$, then the formal series solution of 
Equation \eqref{eq:invEq} is formally well defined to all orders.  
\item If $\lambda_2 = \overline{\lambda_1}$, that is if we consider the complex conjugate case, 
then there is no possibility of a resonance and we can compute
the power series coefficients of the parameterization 
to any desired finite order. 
\item When $\lambda_1, \lambda_2$ are 
complex conjugates, the coefficients of $P$ 
have the symmetry $\overline{p_{nm}} = p_{mn}$ 
for all $m + n \geq 2$.  
This is seen by taking complex conjugates of both sides of the 
homological equation, and using the fact that $Df(\hat p)$ is a 
real matrix.  

Since $\hat p$ is real, choosing complex conjugate
eigenvectors $\xi_2 = \overline{\xi_1}$ enforces the symmetry to all orders.
The power series solution $P$ has complex coefficients, but we obtain the real image of $P$ by taking 
complex conjugate variables.  That is, we define for the real parameters $s_1, s_2$ the function
\[
\hat{P}(s_1, s_2) = P(s_1 + i s_2, s_1 - i s_2),
\]
which parameterizes the real stable/unstable manifold.
\end{itemize}

\subsection{Numerical considerations}
The homological equations derived in the previous section 
allow us to recursively compute the power series
coefficients of the stable/unstable 
manifold parameterization $P$ to any desired order $m + n = N$.
The coefficients are uniquely determined up to the 
choice of the scaling of the eigenvectors.  
In practical applications we have to decide
how to answer the following questions:
\begin{itemize}
\item To what order $N$ will we compute the approximate parameterization?
\item What scale to choose for the eigenvectors?
\item On what domain do we to restrict the polynomial $P^N$?
\end{itemize}

In practice we proceed as follows.  First we choose a convenient 
value for $N$, based on how long we want to let the computations run.
Then, we always restrict $P$ to the unit disk for the sake of numerical 
stability.  Finally, we choose the eigenvector scaling so that the 
last coefficients, the coefficients of order $N$, are smaller than 
some prescribed tolerance.  A good empirical rule of thumb is that 
the truncation error is roughly the same magnitude
as the $N$-th order coefficients.

In practice we can prescribe the size of the $N$-th order terms
as soon as we know the exponential decay rate of the coefficients.  
In the next section we describe the relationship between the 
scale of the eigenvectors and the exponential decay rate.

\subsubsection{Rescaling the eigenvectors}
In Section \ref{sec:formalSeries} we saw that the power series coefficients 
of the parameterization are uniquely determined up to the choice of the 
eigenvector.  Since the eigenvectors are unique up to the choice of length,
we have that the length determines uniquely the coefficients.  
In fact the effect of rescaling the eigenvectors is made completely 
explicit as follows.  The material in this section is discussed in greater 
detail in \cite{maximeJPMe}.

Suppose that 
\[
P(\theta_1, \theta_2) = \sum_{m=0}^\infty \sum_{n=0}^\infty p_{mn} \theta_1^m \theta_2^n,
\] 
is the formal solution of Equation \eqref{eq:invEq},
with
\[
p_{00} = \hat p, \quad \quad \quad 
p_{10} = \xi_1, 
\quad \quad \quad \mbox{and} \quad \quad \quad 
p_{01} = \xi_2,
\]
where $\| \xi_1 \| = \| \xi_2 \| = 1$.  Assuming that $P$ is bounded and analytic 
on the complex poly-disk with radii $R_1, R_2 > 0$, there is a $C > 0$
so that 
\[
| p_{mn} |  \leq \frac{C}{R_1^m R_2^n},
\] 
by the Cauchy estimates.  

Now choose non-zero $s_1, s_2 \in \mathbb{R}$ 
and define the rescaled eigenvectors
\[
\eta_1 = s_1 \xi_1, 
\quad \quad \quad \mbox{and} \quad \quad \quad 
\eta_2 = s_2 \xi_2.
\]
The new parameterization associated with the rescaled eigenvectors is given by 
\[
Q(\theta_1, \theta_2) = \sum_{m=0}^\infty \sum_{n = 0}^\infty  
q_{mn} \theta_1^m \theta_2^n,
\]
where 
\begin{equation} \label{eq:coeffRescale}
q_{mn} = s_1^m s_2^n p_{mn}.
\end{equation}
See \cite{maximeJPMe} for a proof of this identity, and also the discussion in \cite{Cabre1,Cabre3}.
%
The coefficients of the rescaled parameterization 
have the new exponential decay rate given by 
\begin{align*}
|q_{mn}| & = \left| s_1^m s_2^m p_{mn} \right| \\
&\leq s_1^m s_2^n \frac{C}{R_1^m R_2^n} \\
&\leq \frac{C}{\left( \frac{R_1}{s_1}\right)^m \left( \frac{R_2}{s_2}\right)^n}.
\end{align*}

These observations lead to a practical algorithm.  
First compute the parameterization $P$ with an arbitrary 
choice of eigenvector scaling (for example scaled to have length one). 
Then solve the homological equations to some order $N_0$
using this scaling, and compute $C$, $R_1$ and $R_2$ using 
an exponential best fit.  
Suppose that $\varepsilon_0 > 0$ is the desired tolerance, that is 
the desired size of the order $N \geq N_0$ coefficients.  Then  
we choose $s_1$ and $s_2$ so that
\[
\frac{C}{\left( \frac{R_1}{s_1}\right)^N \left( \frac{R_2}{s_2}\right)^N} \leq \varepsilon_0.
\]
Finally we recompute the coefficients $q_{mn}$ for $2 \leq m + n \leq N$.
The rescaled coefficients could be computed from the old coefficients
using the formula of Equation \eqref{eq:coeffRescale}.  In practice
however better results are obtained by recomputing the coefficients $q_{mn}$ from scratch
via the homological equations.

We remark that in the case of complex conjugate eigenvalues we 
want the eigenvectors to be complex conjugates.  Assuming that
$\xi_2 = \overline{\xi}_1$ we take $s_1 = s_2 \in \mathbb{R}$
so that $\eta_2 = \overline{\eta}_1$. 
Also note that by choosing our domain to be the unit poly-disk, we have that
$R_1 = R_2 = 1$, further simplifying the analysis.

\subsubsection{A-posteriori error}
Once we have chosen the polynomial order $N$ and the scaling of the 
eigenvectors, that is once we have uniquely specified our parameterization 
to order $N$, we would like a convenient measure of the truncation 
error.  As mentioned above, a good heuristic indicator is that the 
error is roughly the size of the highest order coefficients 
(assuming we take the unit disk as the domain of our approximate 
parameterization).  In this section we discuss a more quantitative 
indicator.  

We remark that there exist methods 
of a-posteriori error analysis for the parameterization method, 
which -- when taken to their logical conclusion -- lead to 
mathematically rigorous computer assisted error bounds on the 
truncation errors.  
The interested reader will find fuller discussion and more
references to the literature in 
\cite{Cabre3, maximeJPMe, parmChristian, Cana, AMS, MR2821596}
and discussion of related techniques in 
\cite{MR3022075, MR3281845, MR2644324}.

The analysis in the present work is qualitative and we don't
require the full power of mathematically rigorous error bounds.  
Instead we employ an error indicator inspired by the fact 
that the parameterization satisfies the flow invariance 
property given in Equation \eqref{eq:flowConj}.  We 
choose $T {\not =} 0$, and a partition of the interval $[0, 2\pi]$
into $K$ angles, $\alpha_j = 2 \pi j/(K+1)$,
for $0 \leq j \leq K$.  
Since we are interested in the case of complex conjugate 
eigenvalues $\lambda, \overline \lambda \in \mathbb{C}$, we 
define complex conjugate parameters 
\[
\theta_j = (\theta_1^j + i\theta_2^j, \theta_1^j - i\theta_2^j)
= (\cos(\alpha_j) + i \sin(\alpha_j), \cos(\alpha_j) - i \sin(\alpha_j)),
\] 
and the linear mapping
\[
e^{\Lambda T} = \left(
\begin{array}{cc}
e^{\lambda T} & 0 \\
0 & e^{\overline \lambda T}
\end{array}
\right).
\]
which maps complex conjugate inputs to complex conjugate 
outputs.  The a-posteriori indicator is 
\[
\mbox{Error}_{\mbox{{\tiny conj}}}\left(N, T\right) =
\max_{0 \leq j \leq K} \left\|
\phi(P^N(\theta_j), T)
- P^N(e^{\Lambda T} \theta_j)
\right\|. 
\] 
Here $T > 0$ if the complex conjugate eigenvalues 
$\lambda, \overline \lambda$ are stable
and $ T < 0$ if they are unstable.
In practice the flow map $\phi(x, t)$ will be evaluated using 
a numerical integration scheme, and the accuracy of the
indicator is limited by the accuracy of the integrator.

\begin{figure}[t!]
\centering
\includegraphics[width=0.9 \textwidth]{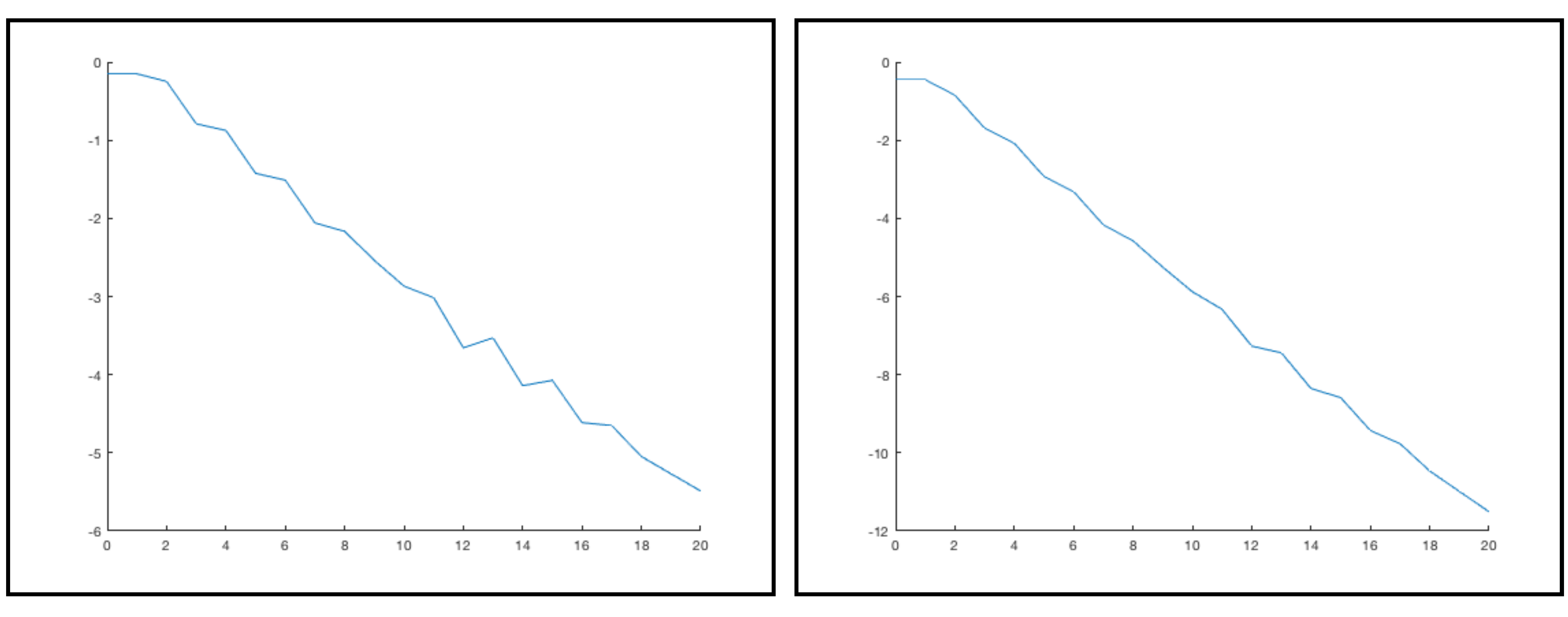}
\caption{\textbf{Rescaling the eigenvector and coefficient decay}: 
the left frame shows the coefficient decay when the eigenvectors 
are scaled to unit length.  The right frame is with scaling one half.
Both figures plot coefficient magnitude $\max_{i + j = n} \log(|p_{ij}|)$ (vertical axis) versus 
polynomial order $n$ (horizontal axis).  When the eigenvector is scaled 
to unit length we see that the order 20 coefficient are on the order of 
$10^{-6}$, which is small but far from machine epsilon.  We should either
increase the order of the polynomial or decrease the scale of the eigenvector.
Indeed, when the scale is decreased to one half we see that the last coefficients
have magnitude on the order of a few thousand multiples of machine epsilon.
} \label{fig:decay}
\end{figure}

\begin{figure}[h!]
\centering
\includegraphics[width=0.9 \textwidth]{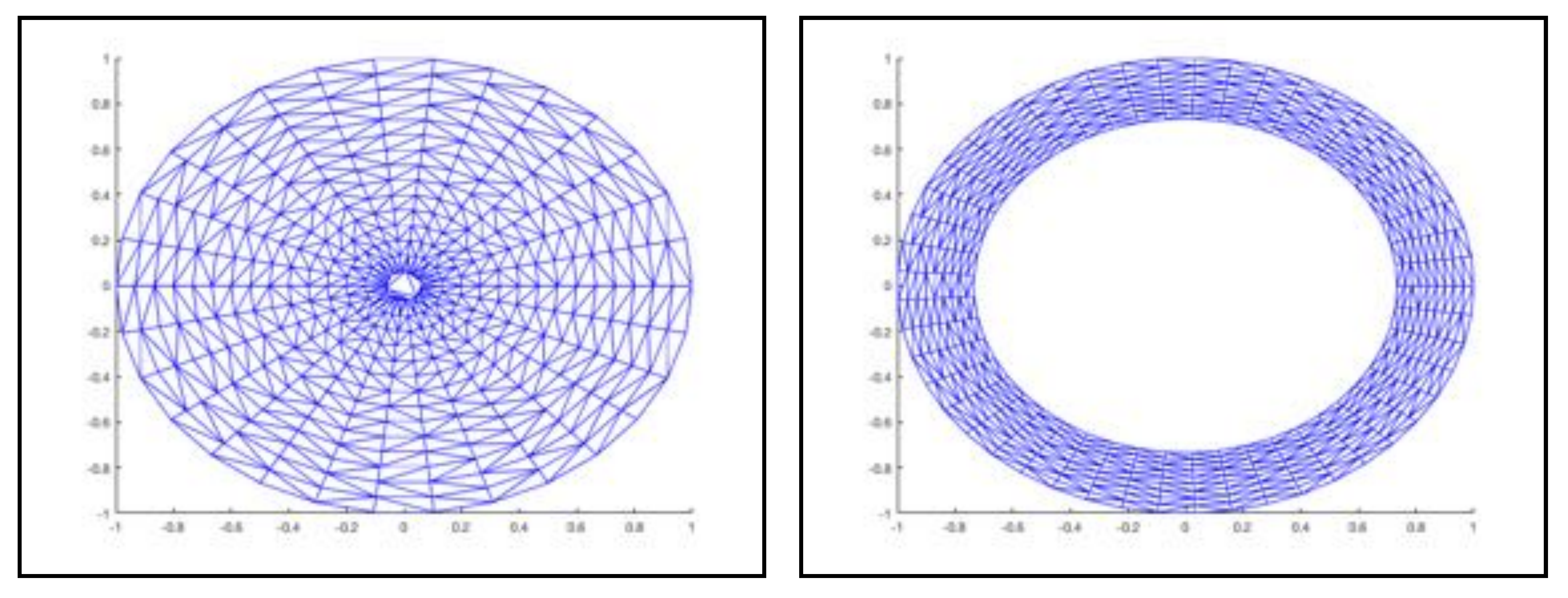}
\caption{\textbf{Triangulating the local invariant manifold and fundamental domain}: 
for the local parameterization we subdivide the unit disk -- fifteen radial subdivisions
by 30 angular subdivisions (left).  Since the domain is simply connected, the triangulation 
is computed using Delauney's algorithm (built into MATLAB).
For a fundamental domain we take the unit circle 
as the outer boundary, and the circle of radius $|e^{-\lambda_u \tau}|$ as the 
inner boundary of an annulus.  We take ten radial subdivisions and fifty angular subdivisions.
We compute a Delauney triangulation, but this ``fills in the hole'' of the annulus.  This 
is fixed by removing triangles with a long side from the triangulation and results in the 
mesh illustrated in the right.
} \label{fig:triangulations}
\end{figure}

\subsubsection{A numerical example} \label{sec:numericalExampleParm}
As an example of the performance of the method, consider the parameterization of the 
two dimensional unstable manifold of the Langford system (Equation \eqref{eq:1}) 
at the equilibirum $p_0$, computed to order $N = 20$.
Figure \ref{fig:decay} illustrates the effect of the choice of the eigenvector scaling on 
the decay rate of the Taylor coefficient.  
We remark that the magnitude of the last Taylor coefficient computed is a good 
heuristic indicator of the size of the truncation error.  
For example if we choose eigenvectors scaled to length one, we obtain the 
decay rate illustrated in the left frame of Figure \ref{fig:decay}, and we see that the 
norm of the largest coefficient of order twenty is about $10^{-6}$.  
On the other hand if we rescale to eigenvector to have length 
$1/2$ then the coefficients decay as in the  right frame of Figure \ref{fig:decay},
and the largest norm of any coefficient of order twenty is now about $10^{-12}$.

\begin{figure}[t!]
\centering
\includegraphics[width=1 \textwidth]{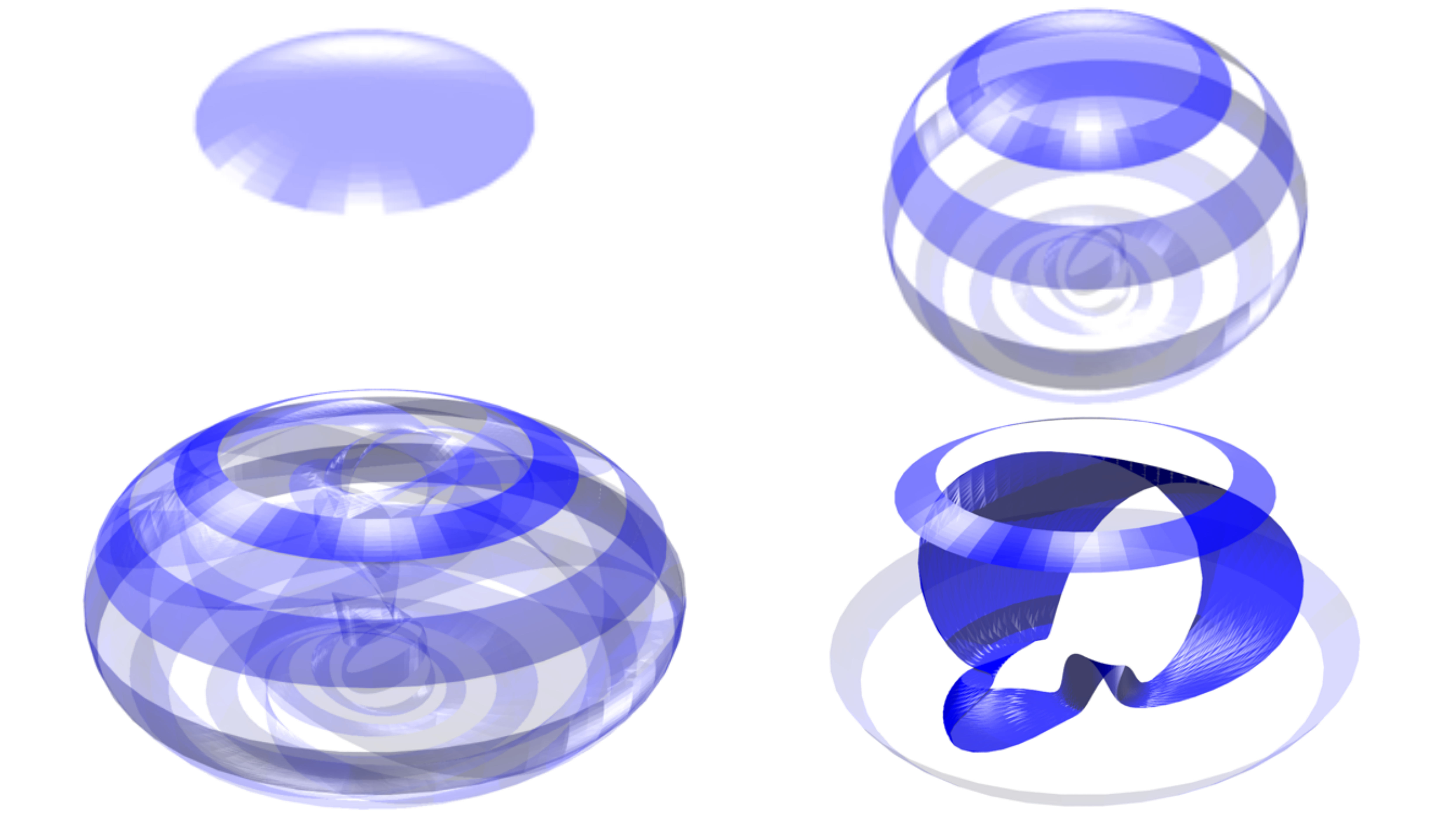}
\caption{\textbf{Growing the unstable manifold one fundamental domain at a time}: 
(Top left) the initial local unstable manifold obtained using the parameterization 
method. (Top right) the initial local manifold parameterization as well as the 
first, third, fifth, seventh, ninth, eleventh, thirteenth, and fifteenth iterate
of the fundamental domain. (Bottom left) the first through thirty third iterates
of the fundamental domain.
(Bottom right) the sixtieth iterate of the fundamental domain, and we see that the 
image is substantially folded.  The first and sixth iterates are shown as well to provide the 
overall shape of the bubble.  
In the bottom frames the initial parameterized local unstable manifold is not shown.} \label{fig:buildBubble}
\end{figure}

To visualize the parameterized local manifold we evaluate the polynomial 
approximation on the unit disk.  First we take a Delaunay triangulation
of the unit disk as illustrated in the left frame of Figure \ref{fig:triangulations}. 
This triangulation of the unit disk is pushed forward to the phase space
$\mathbb{R}^3$ by the polynomial parameterization, resulting in a 
triangulation of the two dimensional local unstable manifold as 
illustrated in the top left frame of Figure \ref{fig:buildBubble}.

To ``grow'' a larger representation of the unstable manifold 
we choose a fundamental domain, for example by taking 
$\tau = 0.25$ and considering the annulus in parameter space formed by the 
boundary of the unit disk and by the circle of radius
$R = \left| e^{\lambda_u \tau} \right| \approx 0.733$.
We mesh this annulus using  
$100$ angular subdivisions and $40$ radial subdivisions,    
as illustrated in the right frame of Figure \ref{fig:triangulations}.
We lift this fundamental domain to the phase space and 
repeatedly apply the time $\tau = 0.25$ map via numerical integration of the 
vertices of the triangulation.  
We refine the mesh whenever any side of a triangle in phase space  
gets too large.   In the present work we measure ``too large'' just by looking at the 
resulting picture.

The top right, and bottom frames of Figure \ref{fig:buildBubble} illustrate the results of iterating a 
triangulation of a fundamental domain for the local unstable manifold at $p_0$,
and we see that the ``bubble'' grows in a quite regular way.  
However, by the time we take $60$ iterates the embedding of the initial annulus 
is becoming quite complicated.

\subsection{Numerical approximation methods for stable/unstable manifolds} \label{rem:invManLit}
The literature devoted to numerical approximation of stable/unstable manifolds
is substantial, and we take a moment to reframe the techniques 
just discussed in this light.
A classic general reference is the work of \cite{1990mmcm.conf..285S}.
The essential remark is that computational methods for studying stable/unstable 
manifolds decompose naturally into two independent tasks: 
\begin{itemize}
\item \textbf{Step 1:} Calculate an approximation of the local invariant manifold.
\item \textbf{Step 2:} Advect the local approximation, ``growing'' the representation of the manifold.
\end{itemize}
A natural approach to Step 1 is to approximate the local manifold to first order, 
reducing the problem to linear algebra.  That is, by computing the eigenvalues/eigenvectors
of the differential at the equilibrium we can approximate the local stable/unstable manifolds
by the stable/unstable eigenspaces.
Step 2 is in general much more difficult, due to the fact that 
nonlinearities cause the manifold
to grow in a highly nonuniform way.  For this reason, much work focuses
on the development of powerful methods for Step 2.  
We refer the interested reader to the works of  
\cite{MR1713086,MR1870261,MR2114735,MR2834454,
MR1391509,MR1391509,MR2338026,MR2179490,MR2989589}, 
and also to the survey paper of 
\cite{MR2136745} for much fuller discussion of the topic.  
The woks just cited develop sophisticated adaptive subdivision 
schemes to control the accuracy and complexity of the advection 
problem, growing the stable/unstable manifolds in a uniform way.  

Another way to fight the nonuniformity encountered at Step 2 is to employ a higher order 
approximation scheme, and hence to compute a larger portion of the 
stable/unstable manifold at Step 1. The idea is that  a 
manifold approximation holding in a large neighborhood of the equilibrium
reduces the dramatic expansion which results from integrating a very small 
polygonal manifold patch until it describes a large portion of the manifold.

The parameterization method as discussed in the present section accomplishes 
this.  Indeed, deriving the homological equations for the system facilitates 
the implementation of  programs which compute the Taylor coefficients 
of the local parameterization to 
any desired order.  
We refer back to the calculations discussed in Section \ref{sec:numericalExampleParm}
where we saw the parameterized local manifold grow quite uniformly after 
the initial high order computation.   
The error from computing the local invariant manifold from 
Step 1 can be estimated even in a 
large neighborhood of the equilibrium using the a-posteriori indicator.

Of course, even when the parameterization method is used in Step 1, 
we have to employ advection schemes to see a larger portion of the 
manifold.  The parameterization method simply provides a high order
approach to Step 1:  it does not eliminate the need for Step 2. In fact 
any of the Step 2 schemes mentioned above could be used in conjunction 
with a high order approximation computed at Step 1 using the parameterization 
method.  
Far from being competitors, the various techniques complement one another.
See \cite{MR2835474,MR3021639,shaneNumericalPaper} for examples of calculations which 
combine the parameterization method in Step 1 with adaptive 
advection schemes in Step 2.

\section{The invariant torus} \label{sec:global}
The first of our two main goals is to study the appearance of the smooth 
attracting invariant torus, the major changes in its dynamics as the bifurcation 
parameter increases -- including interestingly enough the loss of differentiability --
and finally the disappearance of the torus in a global bifurcation resulting in the 
appearance of a new chaotic attractor.  The discussion takes place in a 
Poincar\'{e} section, where periodic orbits are reduced to collections of points,
their stable/unstable manifolds are reduced to curves, and invariant tori 
are reduced to invariant circles.

\subsection{Neimark-Sacker bifurcation in the return map}
We begin by studying the dynamics near the periodic orbit $\gamma$
as the bifurcation parameter $\alpha$ varies.  To this end 
we fix as a surface of section the half plane $\Sigma$ given by $x = 0, y > 0$ (with $z$ free)
and consider the first return map $R \colon \Sigma \to \Sigma$, which is 
well defined in a (possibly quite large) neighborhood of the periodic orbit $\gamma$.
In the discussion that follows all fixed points and $k$-cycles of $R$ are computed 
using standard Newton schemes, and derivatives of the Poincar\'{e} map are computed 
by integrating the variational equations of the flow.

We first observe that 
for $0 < \alpha \leq 0.65$ the first return map has an attracting 
fixed point $p_* \in \Sigma$
corresponding to the attracting periodic orbit $\gamma$
discussed in the introduction. 
At $\alpha_1 \approx 0.697144898322973$ the fixed point looses stability, triggering a
super-critical Neimark-Sacker bifurcation (see \cite{MR699057} for precise definitions).  
This results in the appearance of a smooth attracting invariant circle $\Gamma$
near $p_*$ in $\Sigma$, which is of course an invariant Torus $\mathcal{T}$ for the flow. 
The bifurcation value is computed using a Newton scheme for an appropriate 
augmented system where the parameter $\alpha$ is treated as one of the unknowns, 
hence the bifurcation parameter is known to roughly machine precision.
Such techniques are discussed at length in the classic works of 
\cite{doedel_paris,MR910499,MR1314079}.

\begin{figure}[t!]
\centering
\subfloat{\includegraphics[width = 2.4in]{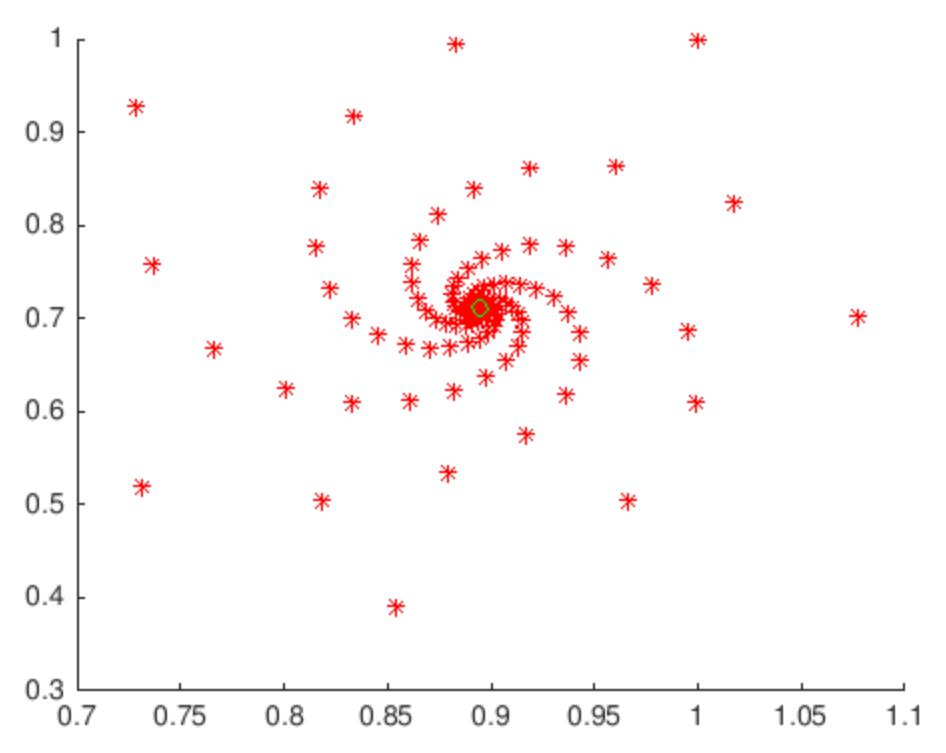}}\quad
\subfloat{\includegraphics[width = 2.4in]{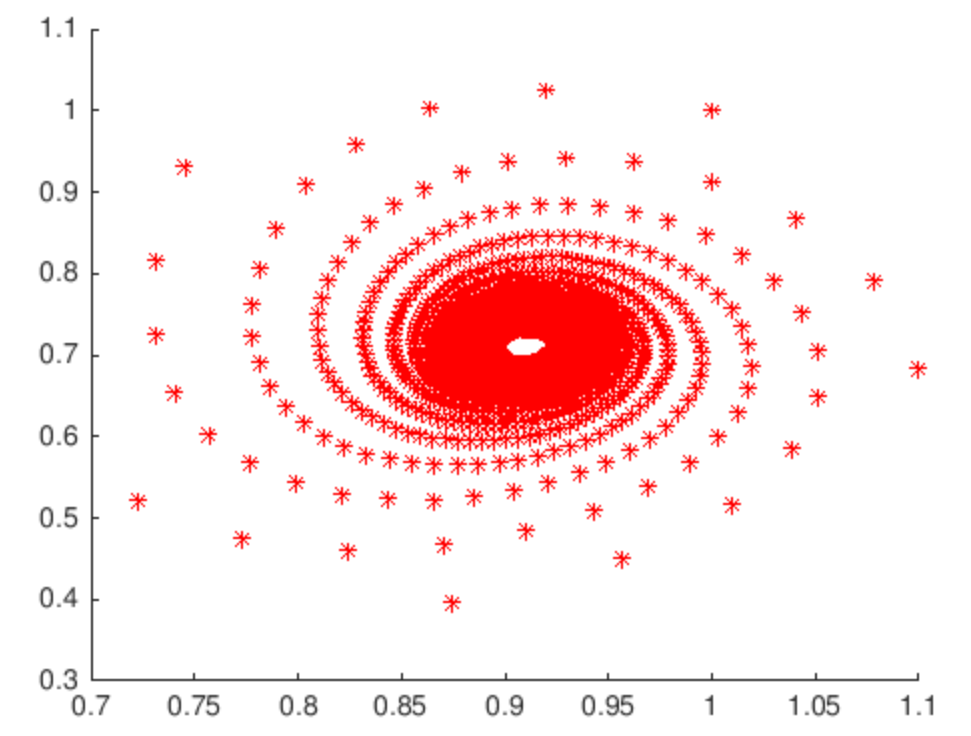}}\\
\subfloat{\includegraphics[width = 2.4in]{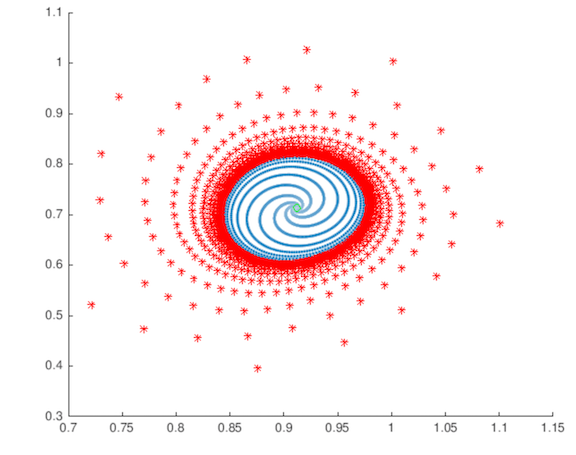}} \quad
\subfloat{\includegraphics[width = 2.4in]{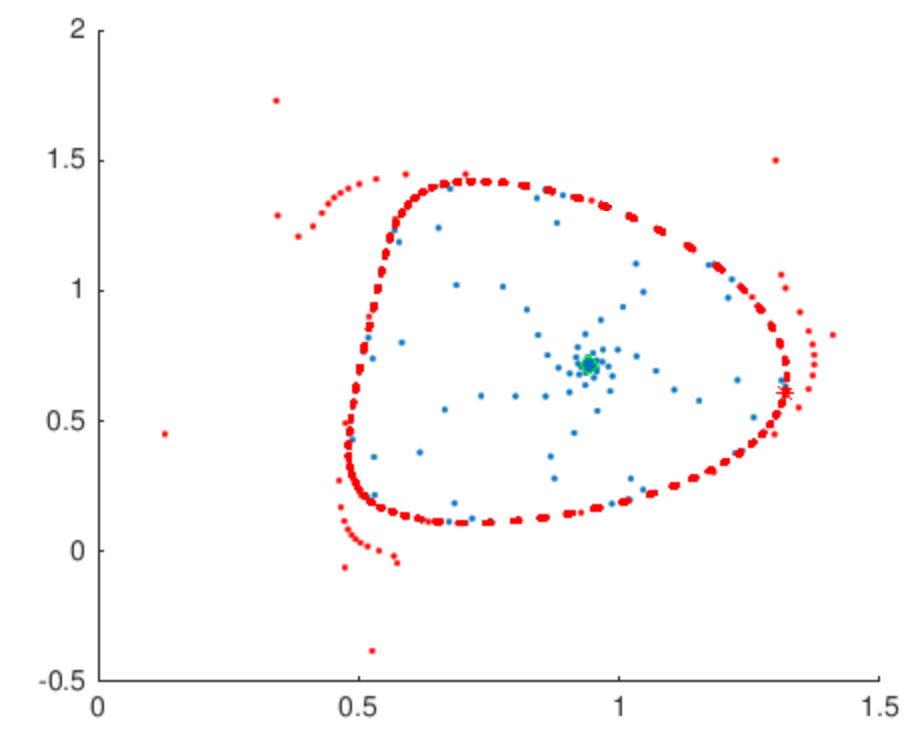}}
\caption{\textbf{Poincar\'{e} sections: attracting fixed point}.
(Top left) attracting fixed point in the Poincare section for $\alpha = 0.65$.
(Top right)  Neimark-Sacker bifurcation for $\alpha = \alpha_1$.
(Bottom left) repelling fixed point of the Poincar\'{e} map and 
attracting invariant circle for $\alpha=0.7$.
(Bottom right) repelling fixed point and attracting invariant circle for $\alpha=0.8$.
In the bottom frames, blue points represent orbits diverging from the repelling 
fixed point and converging to the attracting invariant circle from inside.  In all frames
red points represent orbits converging to the attractor from the outside.
This circle itself is located by iterating the Poincar\'{e} map sufficiently long.
 }
\label{fig:hopf}
\end{figure}

The dynamics in the section just at and just after the bifurcation are illustrated in Figure \ref{fig:hopf}.
For $\alpha > \alpha_1$ the fixed point $p_*$ is repelling and the invariant circle is attracting.
The general theory of the  Neimark-Sacker bifurcation \cite{MR699057} dictates that for 
small enough $\epsilon > 0$, the invariant circle at 
$\alpha = \alpha_1 + \epsilon$ is smoothly conjugate to an irrational rotation.  
The four frames of Figure \ref{fig:hopf} illustrate the initially 
attracting fixed point (top left frame), the Neimark-Sacker bifurcation (top right frame), and the attracting 
invariant circle surrounding the now repelling fixed point where the size of the circle gets larger as 
$\alpha$ increases (bottom two frames).


\begin{figure}[t!]
\centering
\includegraphics[width = 3in]{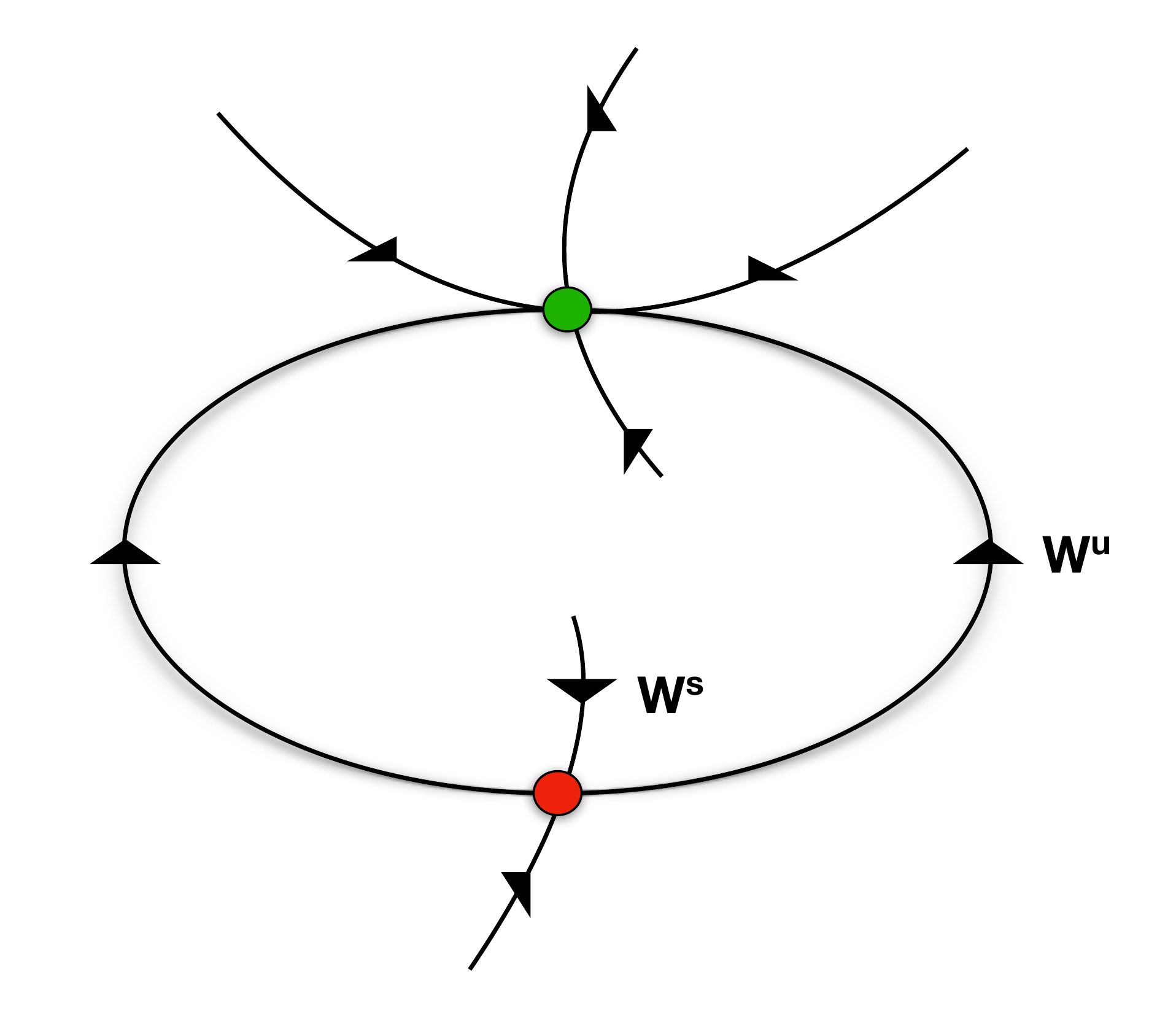}
\caption{\textbf{Schematic of a resonant torus:} the green dot is a stable cycle and the red dot 
a saddle cycle.  Black curves are stable/unstable manifolds.  The unstable manifold 
of the saddle cycle accumulates at the 
stable cycle, forming an invariant topological circle. The smoothness of the circle depends on the 
eigenvalues at the stable cycle, and if the eigenvalues at the stable cycle 
are complex conjugate the circle cannot 
be globally differentiable.  }
\label{fig:resSchematic}
\end{figure}

\subsection{Resonant tori}\label{sec:resonantTori}
When  $\alpha > 0$ is large enough the dynamics on the invariant circle $\Gamma$ 
change in a fundamental way, as we now discuss. 
We say that $q$ is a period $k$ point for $R \colon \mathbb{R}^2 \to \mathbb{R}^2$
if $q, R(q), R^2(q), \ldots, R^{k-1}(q)$ is a collection of $k$ distinct points 
having $R^k(q) = q$.  We say that the point set 
\[
\mathcal{Q} = \mbox{orbit}(q) = q \cup R(q) \cup R^2(q) \cup \ldots \cup R^{k-1}(q) \subset \mathbb{R}^2,
\]
is a $k$-cycle for $R$.  Notions like stability and stable/unstable manifolds 
of $k$-cycles are defined in 
the obvious way after observing that $q$ and each of its iterates
are fixed points of the composition map $R^k$.
See for example \cite{jorgeMePerParm} for precise definitions and references 
to the literature.
The following notion is critical in the discussion to follow.

\begin{definition}
Let $R \colon \mathbb{R}^2 \to \mathbb{R}^2$ be a smooth map of the plane and 
$\Gamma \subset \mathbb{R}^2$ be a topological circle invariant under $R$.
We say that $\Gamma$ is a \textit{simple resonant invariant circle} if there 
is an attracting $k$-cycle $\mathcal{Q}_1$ and a saddle 
$k$-cycle $\mathcal{Q}_2$ so that 
\[
\Gamma = \{\mathcal{Q}_1 \} \cup \{\mathcal{Q}_2 \} \cup W^u(\mathcal{Q}_2).
\]
\end{definition}

\noindent The situation is that 
 the one dimensional unstable manifold of the saddle cycle is completely absorbed into the 
basin of attraction of the stable cycle, in such a way that a circle is formed.  
In this case the dynamics on the invariant circle are conjugate 
to a gradient system. 
Observe that  the unstable manifold of the saddle cycle is smooth (analytic if the 
map is), even if -- as we will see below -- the regularity of the invariant circle is 
another matter completly.
The situation is illustrated in Figure \ref{fig:resSchematic} for the 
simple case of a one-cycle.

\begin{remark}
If the decomposition of $\Gamma$ requires multiple stable and saddle cycles
of different periods we say that we have a 
compound resonant invariant circle. However we do not encounter this situation in the present study,
and for this reason we usually drop the term ``simple'' and say simply that we have a resonant 
invariant circle.  
\end{remark}

\begin{remark}
Suppose that $R$ is a Poincar\'{e} map for a 3-dimensional smooth vector field $f$,
and that $R$ has a simple resonant invariant circle $\Gamma$.
Then the flow generated by $f$ has a resonant invariant torus given by 
\begin{equation*}
\mathcal{T} :=\left\{ \phi _{t}(v):v\in \Gamma ,t\in \mathbb{R}\right\}.
\end{equation*}%
\end{remark}

We further remark that the global regularity of a resonant invariant circle (or torus) is determined
by the linearization of $R$ at $q_1$ or any of its iterates.  
So for example if $DR(q_1)$ has real distinct stable eigenvalues then the 
resulting invariant circle is finitely differentiable, with regularity determined by the ratio of these 
eigenvalues.   If on the other hand $DR(q_1)$ has complex conjugate eigenvalues then the 
torus in phase space is only $C^0$, as the unstable manifold of $\mathcal{Q}_2$ is forced to approach 
$\mathcal{Q}_1$ in a spiraling fashion and the resulting curve cannot be differentiable or even 
Lipschitz at $q_1$ or any of its iterates.

\subsection{Resonant tori in the Langford system}
The formation of a resonant invariant torus in the Langford system
of Equation \eqref{eq:1} involves a 
global bifurcation which 
can be observed in the Poincar\'{e} section, as we now describe.  
We begin with the observation that at 
$\alpha = 0.82$ there is an attracting  $3$-cycle, 
which we denote by $\mathcal{Q}_1$, and which lies outside 
the invariant circle $\Gamma$.  The basin of attraction of the $3$-cycle 
is fairly small, as there is a saddle 
type $3$-cycle nearby, which we denote by $\mathcal{Q}_2$.

\begin{figure}[t!]
\centering
\subfloat{\includegraphics[width = 3in]{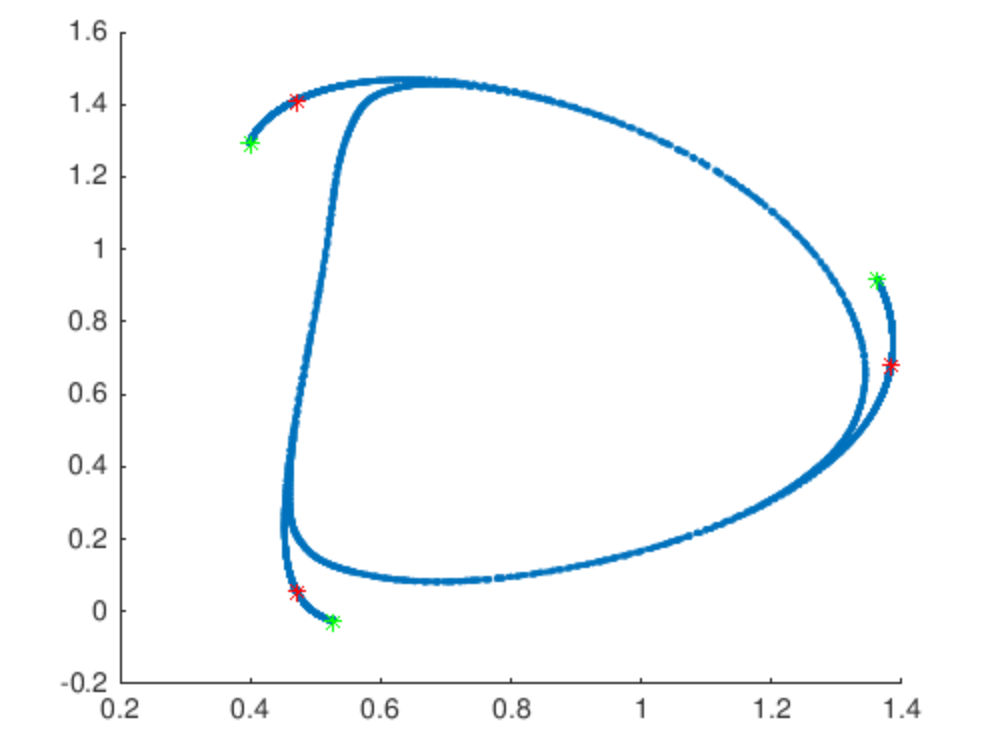}}\quad
\subfloat{\includegraphics[width = 3in]{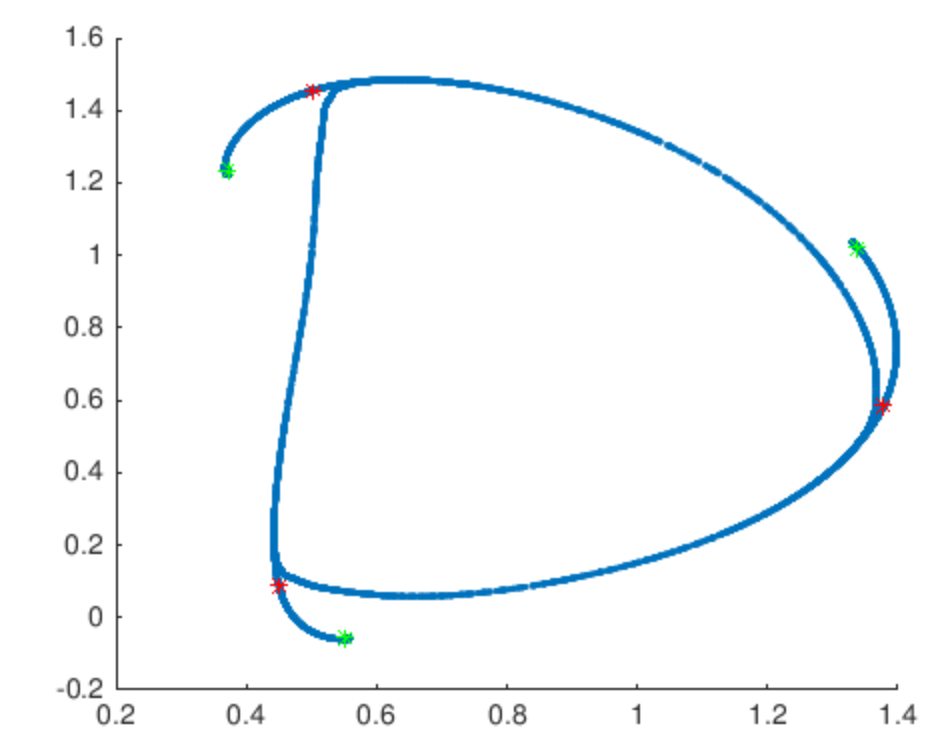}}
\caption{\textbf{Unstable manifold of the saddle  $3$-cycle accumulating 
to both the invariant circle and the stable $3$-cycle:}
The three red dots illustrate the saddle $3$-cycle $\mathcal{Q}_2$
while the three green points illustrate the attracting 
$3$-cycle $\mathcal{Q}_1$.  The blue curve represents the unstable manifold of 
$\mathcal{Q}_2$.  In both cases
half the unstable manifold accumulates on $\Gamma$, and half
accumulates on $\mathcal{Q}_1$.  
In the left frame $(\alpha = 0.82)$ the saddle is far from the invariant circle but in the 
right frame $(\alpha = 0.8224)$ it has moved much closer in anticipation of the 
coming global bifurcation.
 }
\label{fig:preResTorus}
\end{figure}

For parameter values near $\alpha = 0.82$, the unstable manifold of $\mathcal{Q}_2$
has the following behavior:  half of $W^u(\mathcal{Q}_2)$ accumulates on 
the attracting invariant circle $\Gamma$ while the other half 
accumulates to the attracting $3$-cycle $\mathcal{Q}_1$. 
Things remain much the same for 
nearby parameter values, for example at $\alpha = 0.8224$, with the 
caveat that the saddle $3$-cycle $\mathcal{Q}_2$ has moved closer to $\Gamma$.
The situation is illustrated in Figure \ref{fig:preResTorus}.
The stable manifold of $\mathcal{Q}_2$ appears to form a separatrix 
between the basins of attraction of $\Gamma$ and $\mathcal{Q}_1$.  
See for example the left frame of Figure \ref{fig:preResTorus}.

\begin{figure}
\begin{center}
\includegraphics[width=0.425\textwidth]{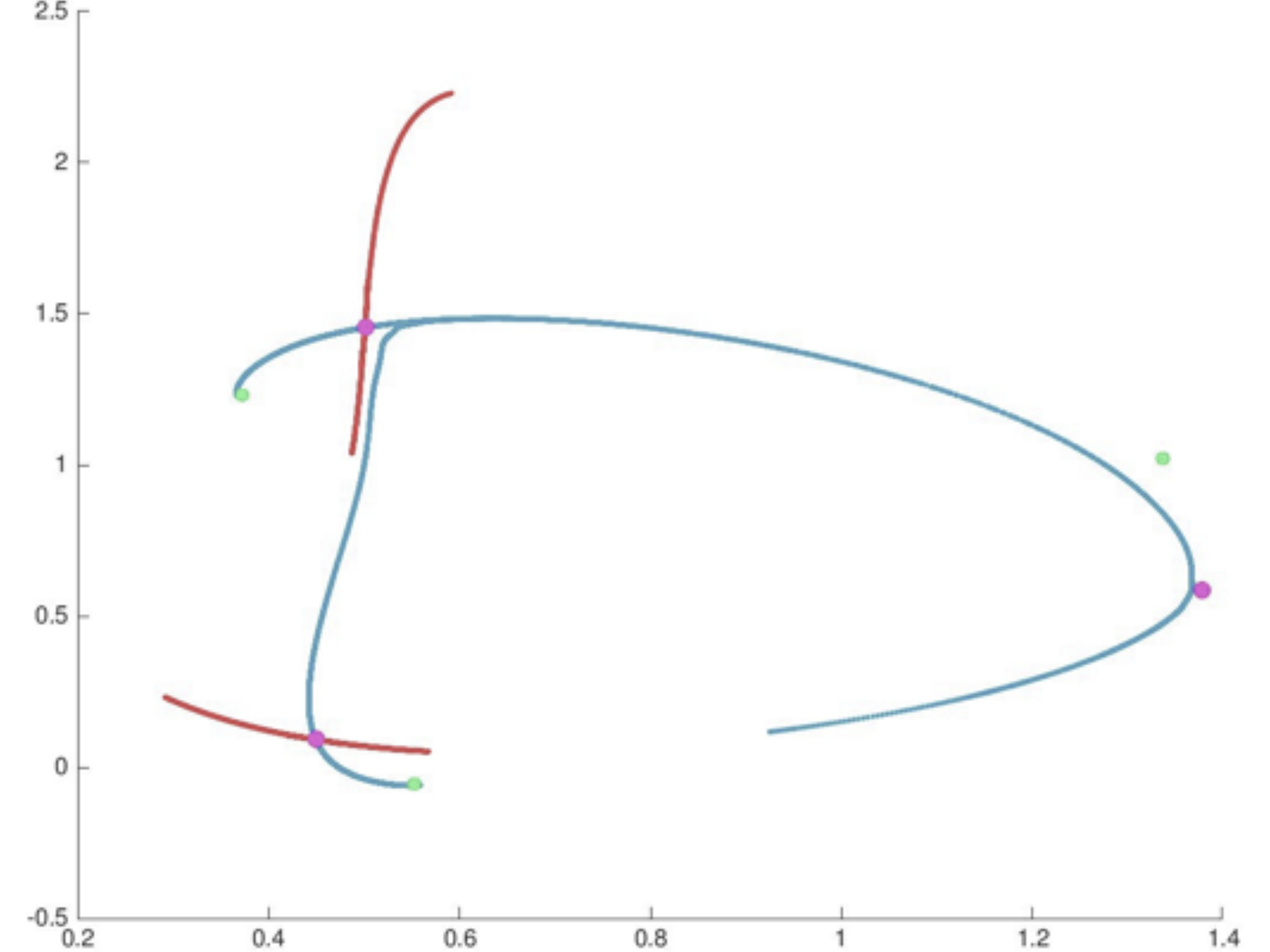}\quad
\includegraphics[width=0.425\textwidth]{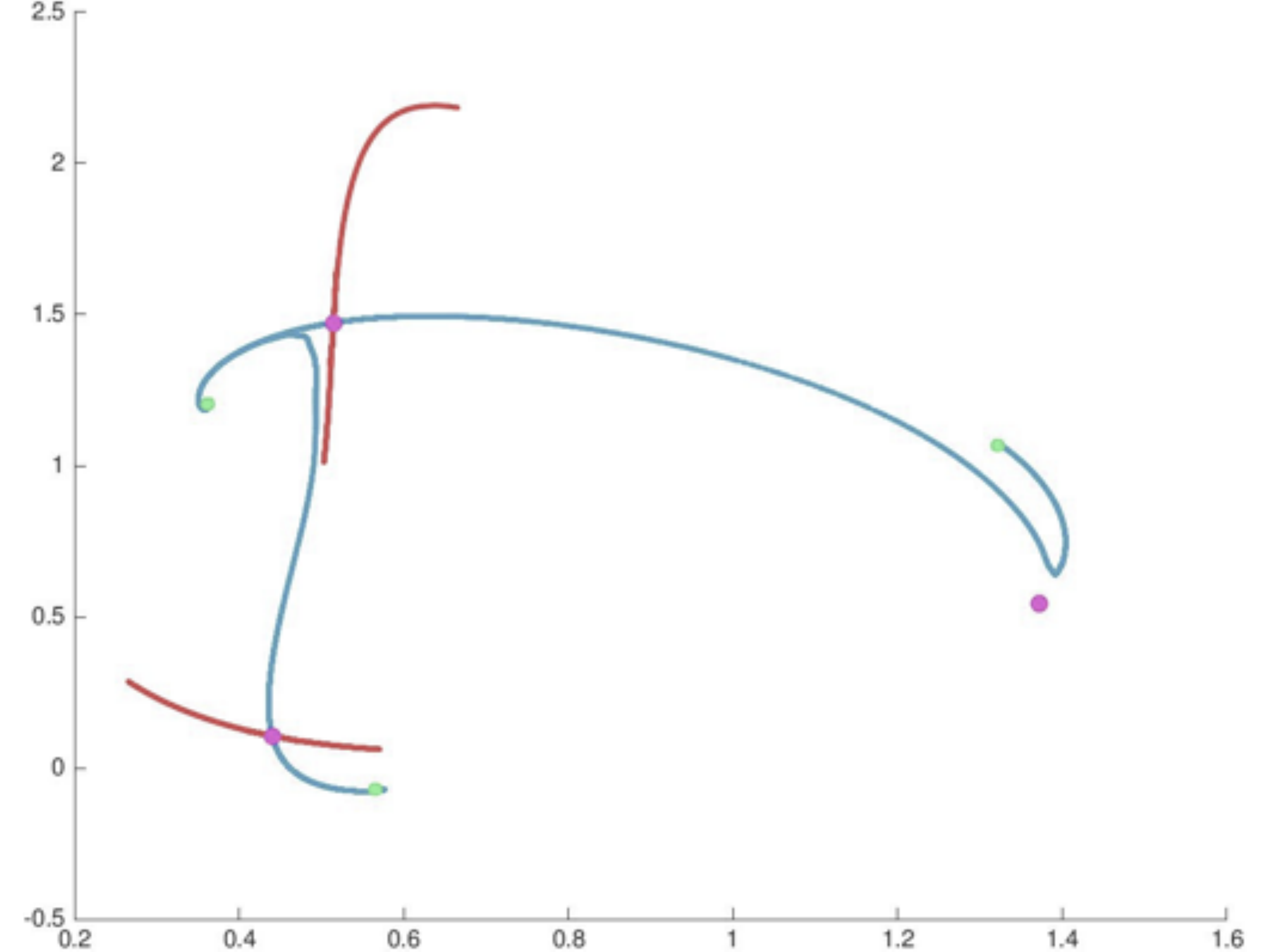}
\end{center}
\caption{\textbf{The resonant torus: before and after.} $W^s(\mathcal{Q}_2)$ 
is red and $W^u(\mathcal{Q}_2)$ is blue. The stable cycle is marked with three green points
and the saddle cycle marked by magenta.
 (Left) at $\alpha = 0.822$ note the top left magenta point.  The 
 left side of its unstable manifold goes to the attracting orbit (green point) 
 while its right side wraps around the attracting invariant circle. (Right) at $\alpha = 0.826$ 
 the bifurcation has occurred and the invariant circle is resonant,
 now comprised of the two $3$-cycles and the unstable manifold.
Looking again at the top left magenta point,
the left side of $W^u(\mathcal{Q}_2)$ still accumulates to the top left green point in the attracting orbit, 
the right side now loops back and is ``captured'' by the top right green point.  Hence both 
sides of the unstable manifold now accumulate to the attracting cycle.  }
\label{fig:bifurcation1}
\end{figure}

\begin{figure}[h!]
\centering
\subfloat{\includegraphics[width = 2.5in]{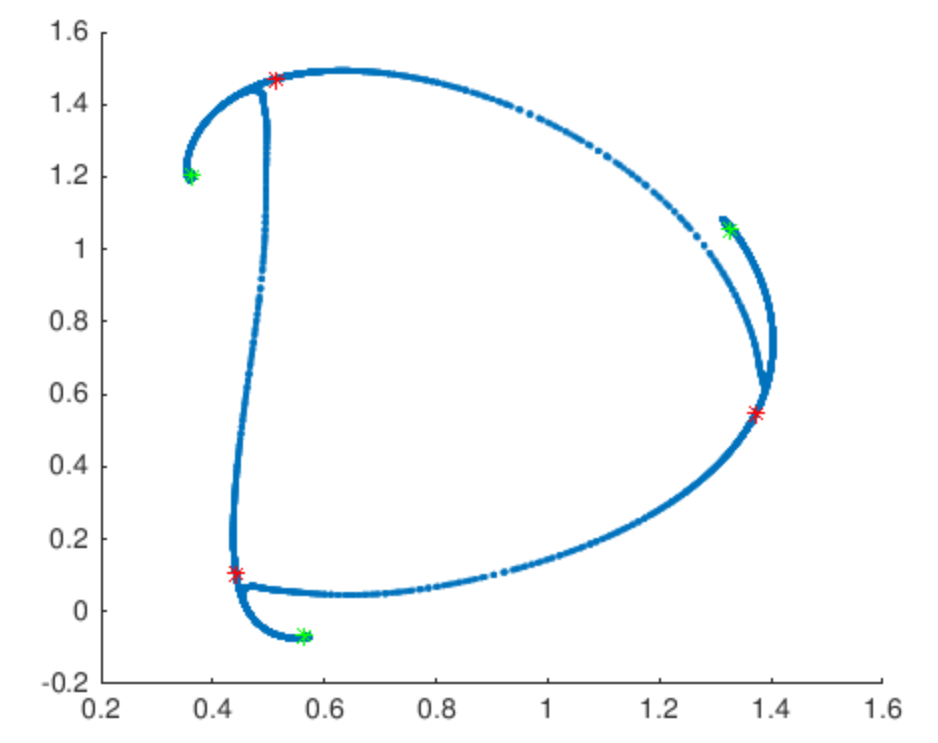}}\quad
\subfloat{\includegraphics[width = 2.5in]{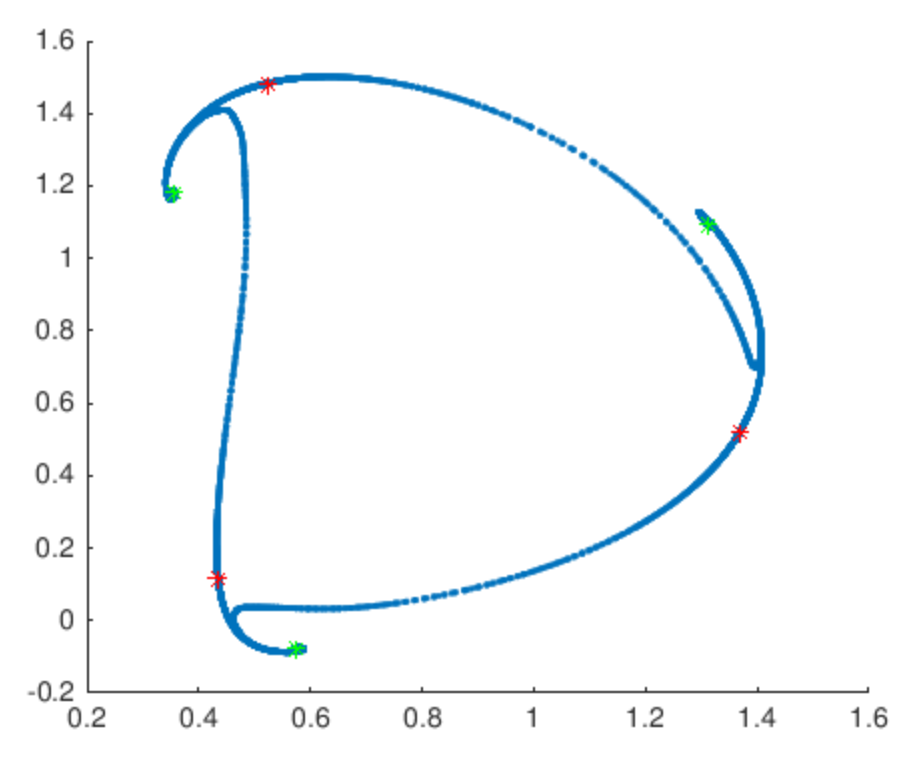}}
\caption{\textbf{Resonant invariant circles in the Poincar\'{e} section.}
Colors have the same meaning as in Figure \ref{fig:preResTorus}.
(Left) at $\alpha = 0.825$ the saddle $3$-cycle has collided 
with the invariant circle.
(Right) a larger value of $\alpha = 0.83$,  
 the resonant invariant circle $\Gamma$ is becoming less regular 
 and lacking differentiability.
 }
\label{fig:resTorusBorn}
\end{figure}

For some parameter value $0.8224 < \alpha_2 < 0.825$ there appears to be
a global bifurcation where $\mathcal{Q}_2$ collides with the invariant circle 
$\Gamma$.  At this point $W^u(\mathcal{Q}_2)$, rather than accumulating on the 
invariant circle $\Gamma$,
\textit{has become $\Gamma$}.   Both halves of the unstable manifold
accumulate at $\mathcal{Q}_1$, which is now inside the invariant circle
as well.  See Figure \ref{fig:bifurcation1} for an illustration of 
the phase space configuration just before and just after the global bifurcation.  
The situation remains for parameter values $\alpha > \alpha_2$
as illustrated in Figure \ref{fig:resTorusBorn}.

Let $q_1$ and $q_2$ be points on the $3$ cycles 
$\mathcal{Q}_1$ and $\mathcal{Q}_2$ respectively.  By numerical calculation we 
find that  the eigenvalues of $DR^3(q_1)$ are complex conjugate stable.
Then the resonant invariant circle appearing in this bifurcation is only $C^0$.
This is a dramatic change, as for $\alpha < \alpha_2$ the simulations indicate
that the torus is smooth (at least finitely differentiable).  The global 
bifurcation just described gives a vivid natural example of a low regularity invariant 
manifold for a smooth (in fact analytic) vector field.    

\begin{figure}
\subfloat{\includegraphics[width = 2in]{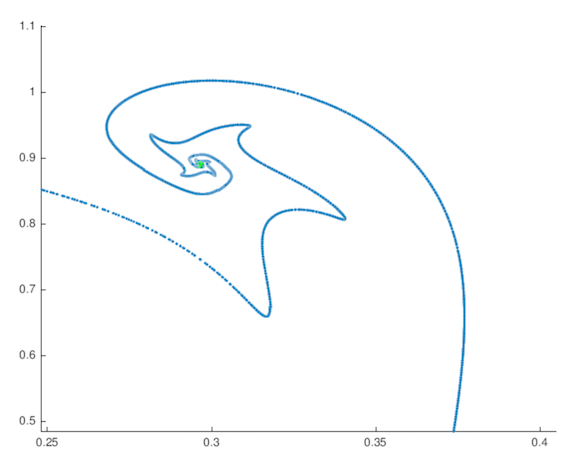}}\quad
\subfloat{\includegraphics[width = 2in]{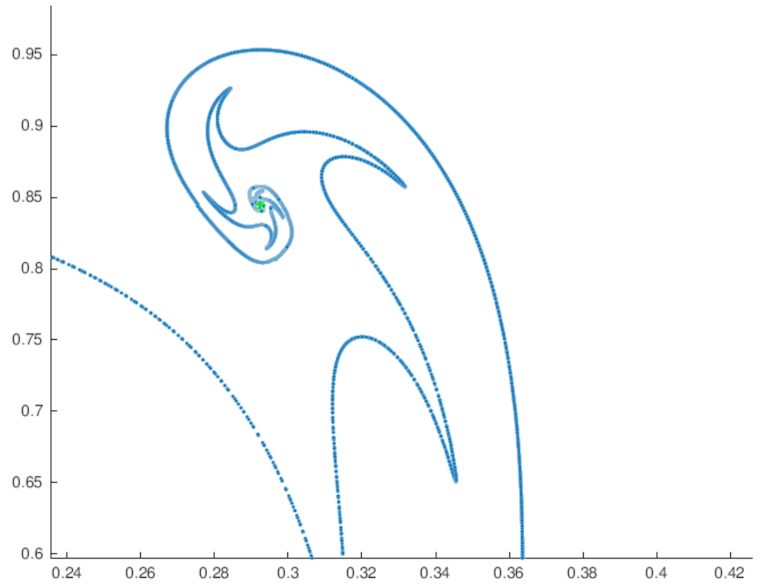}} \quad
\subfloat{\includegraphics[width = 2in]{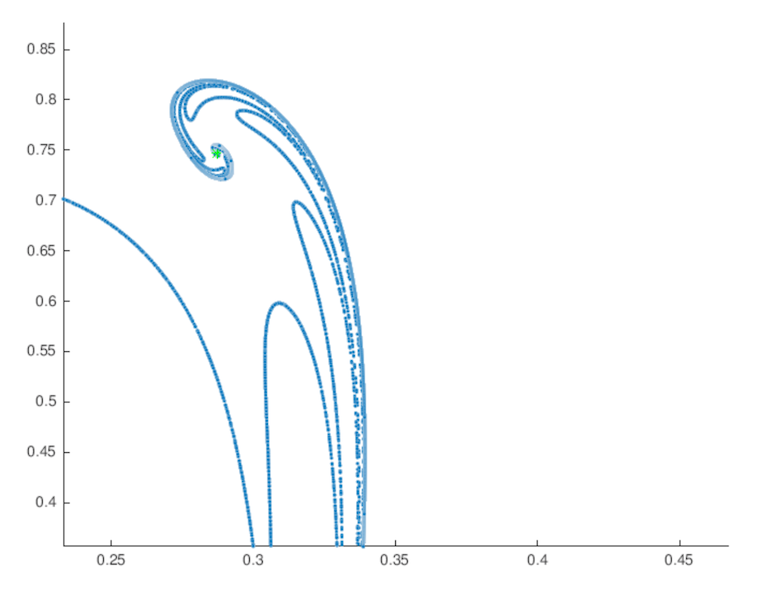}}
\caption{\textbf{Resonant torus near the attracting $3$-cycle:}
Closeups on an attracting period $3$ point for three different values of $\alpha$
larger than $\alpha_2$.  It is clear that the invariant circle is becoming increasingly irregular, 
developing sharp cusp-like edges in its embedding.    
 }
\label{fig:spirialWaves}
\end{figure}

\begin{figure}[t!]
\begin{center}
\includegraphics[width=0.75\textwidth]{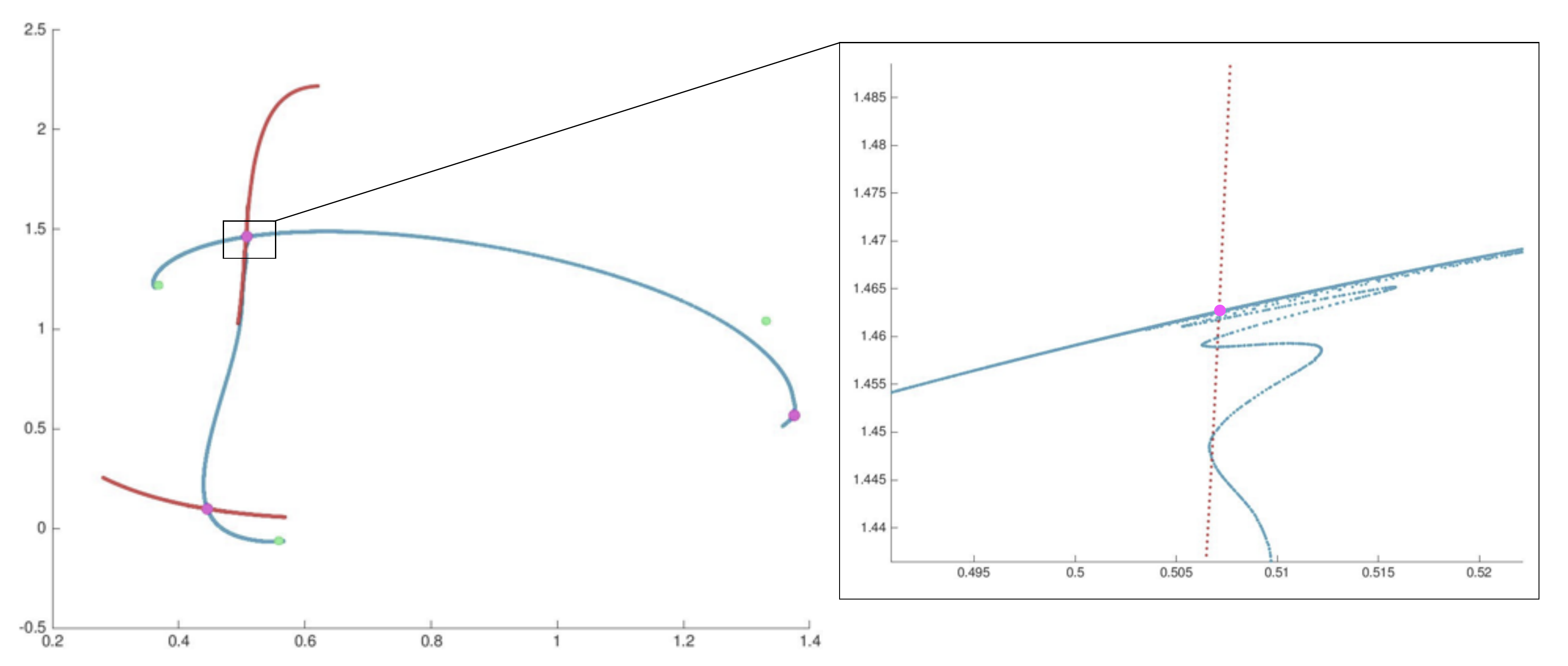}
\end{center}
\caption{\textbf{Transient chaos:}  At $\alpha = 0.8225$ there are transverse intersections of 
$W^u(\mathcal{Q}_2)$  and $W^s(\mathcal{Q}_2)$ indicating the presence of 
Smale horseshoes and thus chaotic dynamics
near the attracting resonant invariant torus $\mathcal{T}$ in phase space.
 Note that, in contrast, the dynamics on the 
attractor are very simple -- that is, the chaos is transient. This image tells us 
that, since $W^u(\mathcal{Q}_2)$  and $W^s(\mathcal{Q}_2)$ do not intersect at $\alpha = 0.8224$
(see Figure \ref{fig:bifurcation1}), the global 
bifurcation to a resonant torus occurs for $0.8224 < \alpha_2 < 0.8225$. }
\label{fig:birfucation2}
\end{figure}

\subsection{Transient chaotic motions}
Increasing the bifurcation parameter $\alpha$ past the global bifurcation at
$\alpha_2$ shows that the embedding of the attracting $C^0$ resonant 
invariant circle $\Gamma$ appears to get even ``wilder''.  See the
three frames of Figure \ref{fig:spirialWaves}.   The blue curve 
illustrates the unstable manifold of $\mathcal{Q}_2$ and 
in the left two frames is contained in the invariant circle $\Gamma$, 
indicating that the circle is loosing regularity.

A more quantitative discussion about dynamical complexity in the system begins by
observing  that just before the global bifurcation at $\alpha_2$, as $\mathcal{Q}_2$ 
is approaching $\Gamma$,
there is the appearance of chaotic dynamics in the system.  To see this 
we observe that in the top left corner of the left hand frame of Figure \ref{fig:bifurcation1},
 the unstable manifold of $\mathcal{Q}_2$ (vertical blue curve) is to the right side of the stable manifold 
 (vertical red curve).  In the top left corner of the right hand frame in 
the same figure we see that the situation is reversed: the unstable manifold
now being to the left of the stable one (again vertical blue and red respectively).
Since the curves move continuously this suggest that there should be a 
range of parameters where they intersect.  

Figure \ref{fig:birfucation2} shows that this is indeed the case. At $\alpha = 0.8225$
we can see that  the resonant torus in phase space has not formed yet, as the unstable manifold (blue
curve) does not accumulate at the stable $3$-cycle (green point).  Here, on close 
inspection we see that $W^s(\mathcal{Q}_2)$ and $W^u(\mathcal{Q}_2)$ do intersect, in fact transversally.
Then there is a Smale horse shoe and hence chaotic dynamics in a neighborhood of $\Gamma$
\cite{MR0228014}.
Note however that the invariant circle $\Gamma$ is still attracting and that the dynamics on $\Gamma$
are simple before and after the global bifurcation at $\alpha_2$.  This suggests that   
the chaotic motions are only transient, in the sense that the horse shoe is not 
in the attractor.

\subsection{Destruction of the torus and appearance of a chaotic attractor}
By further increasing the bifurcation parameter we eventually observe the destruction of 
the invariant torus, as we now describe.  We begin by 
recalling a classical result concerning the disappearance of an invariant circle.  
We refer to the discussion in \cite{Vadim} for the details of the proof and 
generalizations to higher dimensions.  See also 
\cite{MR516994} and 
\cite{MR709899,MR701669,MR845031}.
The set up is as follows.

Suppose that a one parameter family of smooth discrete time dynamical systems 
$R_\alpha \colon \mathbb{R}^2 \to \mathbb{R}^2$
has at $\alpha = \alpha_0$ a resonant $C^0$ invariant circle 
$\Gamma \subset \mathbb{R}^2$ formed by the closure of the unstable manifold of a 
saddle cycle $\mathcal{Q}_2$ accumulating to a stable cycle $\mathcal{Q}_1$ as in 
Figure \ref{fig:resSchematic}.
We then have a resonance region in the sense of \cite{Vadim} and, while it may be obvious
it is nevertheless important to note that, the resonant torus $\Gamma$ 
is robust
under small perturbations including small changes in $\alpha$.  This is because the saddle 
cycle, the attracting cycle, and any local unstable manifold of the saddle are structurally stable
objects.  Since a local unstable manifold is all that is required to reach the basin of 
attraction of $\mathcal{Q}_1$,
we have that the resonant torus is robust. 
It follows that there is a one parameter family of attracting invariant circles 
$\Gamma(\alpha)$ for $\alpha$ near $\alpha_0$.  Indeed the tori vary continuously in $\alpha$,
again see  \cite{Vadim}.

\begin{figure}[t!]
\centering
\subfloat{\includegraphics[width = 2.5in, height=2.1in, valign=c]{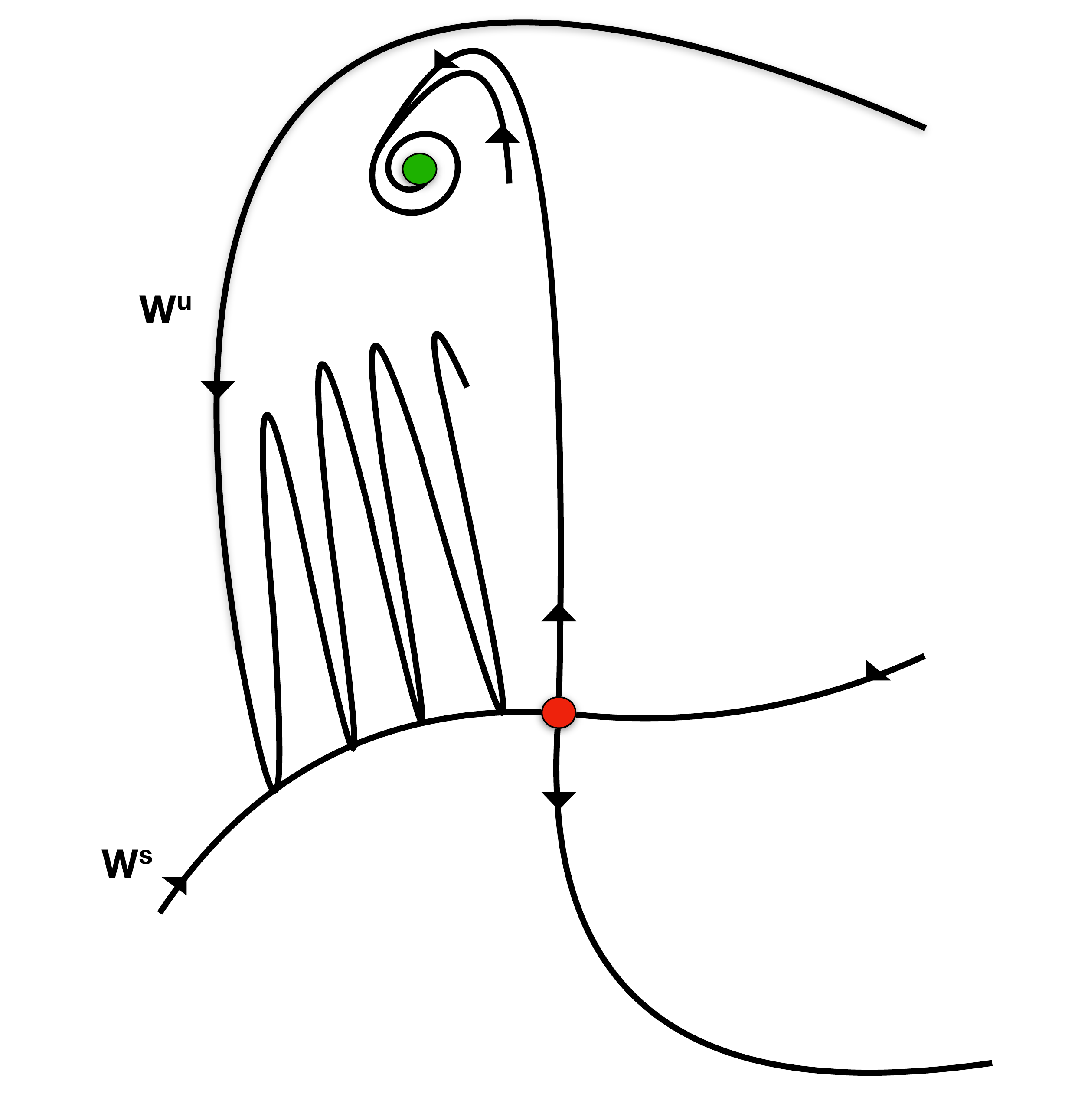}}\quad
\subfloat{\includegraphics[width = 2.5in, height=2.1in, valign=c]{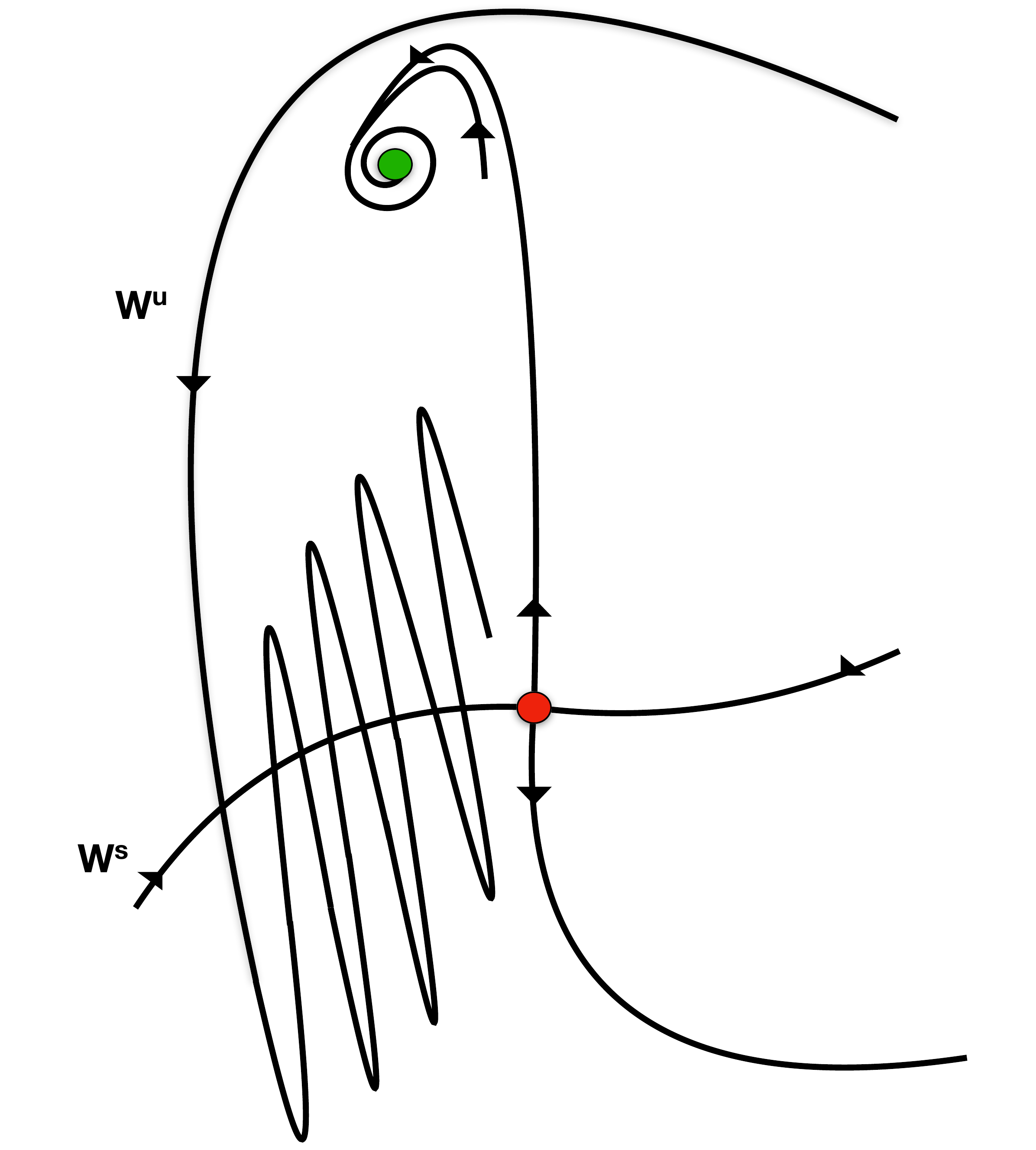}} 
\caption{\textbf{Schematic of the tangency:}
one mechanism for the destruction of a $C^0$ resonant invariant circle (see Figure
\ref{fig:resSchematic}) is the formation of a homoclinic tangency (we stress that many 
tangencies may appear at the same time
\cite{MR699057}).  Suppose that before 
the tangency $W^u$ of the saddle cycle is absorbed into the basin of attraction of a stable cycle (green point).
Once a tangency forms $W^u$ of the saddle cycle must also accumulate in a $C^1$ fashion on
$W_{\mbox{{\tiny loc}}}^u$ near the saddle (red point) -- while still accumulating at 
the stable cycle -- and the resonant torus is destroyed.   
The bifurcation is discussed in greater detail in  \cite{Vadim}. 
 }
\label{fig:Tangency}
\end{figure}

Now, suppose that at some parameter $\hat{\alpha}> \alpha_0$  the invariant torus no longer exists. 
Then, by the least upper bound property of $\mathbb{R}$ there is an $\alpha^*$ so that $\Gamma$ 
is a continuous attracting invariant circle on the interval $[\alpha_0, \alpha^*)$, but that for $\alpha > \alpha^*$
this fails to be true.
The {\it theorem on torus breakdown} \cite{Vadim} asserts the following 
three possible mechanisms for the destruction of $\Gamma$: (i) loss of cycle stability-- 
i.e. a local bifurcation at $\mathcal{Q}_1$, 
(ii) occurrence of a tangency bifurcation of the stable and saddle cycles on $\Gamma$, or
(iii) formation of a homoclinic tangency between $W^u(\mathcal{Q}_2)$ and $W^s(\mathcal{Q}_2)$.  
Mechanism (iii) is illustrated schematically in Figure \ref{fig:Tangency}.

Of course the theorem gives only a trichotomy.  It does not say which alternative actually occurs
in a given example, and it is with this in mind that we investigate the fate of the invariant torus in the 
Langford system.  The situation is illustrated in the two frames of Figure \ref{fig:Tangency}, where we observe 
the formation of a homoclinic tangle for $W^{u,s}(\mathcal{Q}_2)$.  Since these manifolds do not intersect in 
the frame on the left, and do intersect in the frame on the right, they must develop a tangency at some 
point $0.92 < \alpha_3 < 0.93$.  While the closure of the unstable manifold (blue) on the left is still an
attracting invariant circle, in the right frame $W^u(\mathcal{Q}_2)$ is no longer contained in the attractor
suggesting that the torus was destroyed in the homoclinic tangency -- that is that we have in alternative (iii).

\begin{figure}
\centering
\subfloat{\includegraphics[width = 6.5in, height=2.5in, valign=c]{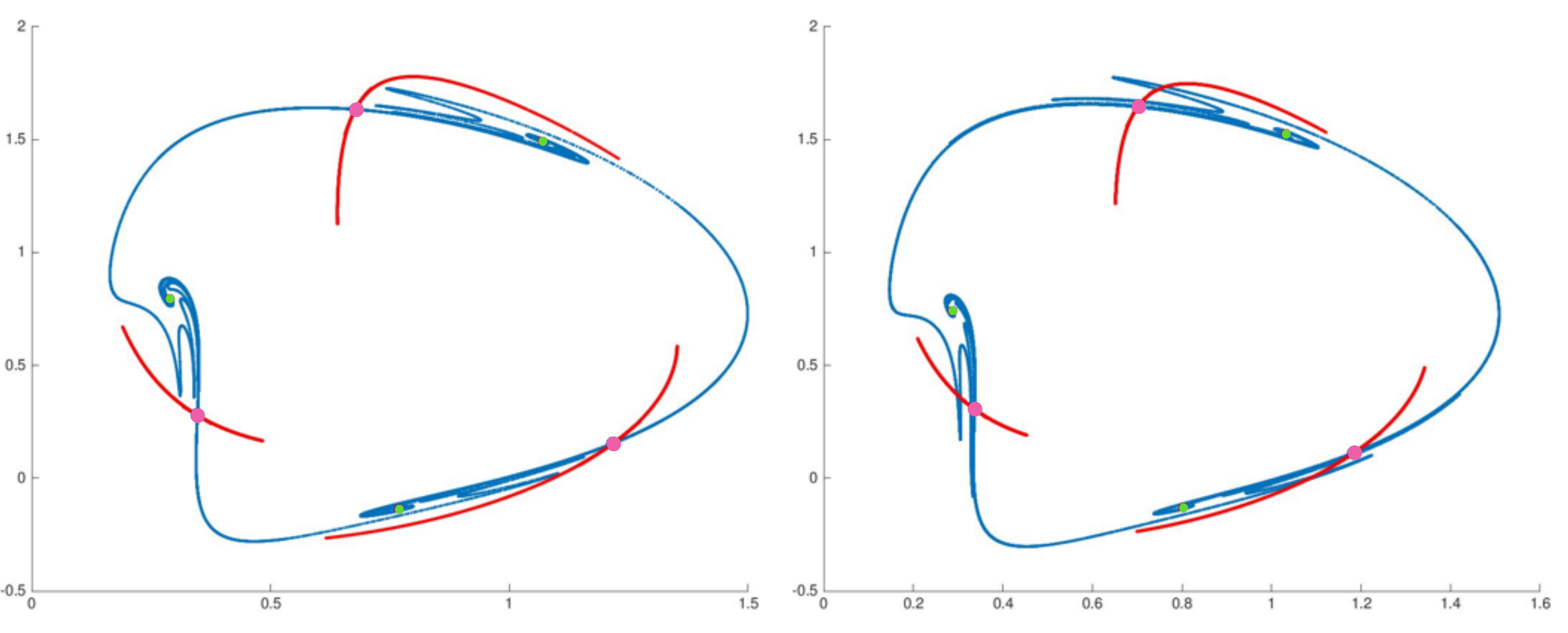}}
\caption{\textbf{Global bifurcation and destruction of the invariant circle:}
(Left) for $\alpha=0.92$ the unstable manifold of $\mathcal{Q}_2$ (blue) still accumulates
 (albeit in a complicated way)
to the stable three cycle $\mathcal{Q}_1$, and the attractor is still a resonant invariant circle.  
That being said, we can see sharp turns developing in the embedding of the unstable manifold near
the stable manifold (red curve) of $\mathcal{Q}_2$. 
(Right) at $\alpha=0.93$, these sharp turns in the unstable manifold 
embedding have moved across the stable manifold resulting in points of transverse intersection 
and hence a Smale horse shoe nearby.  
 }
\label{fig:Tangency}
\end{figure}

\begin{figure}
\centering
\subfloat{\includegraphics[width = 5in]{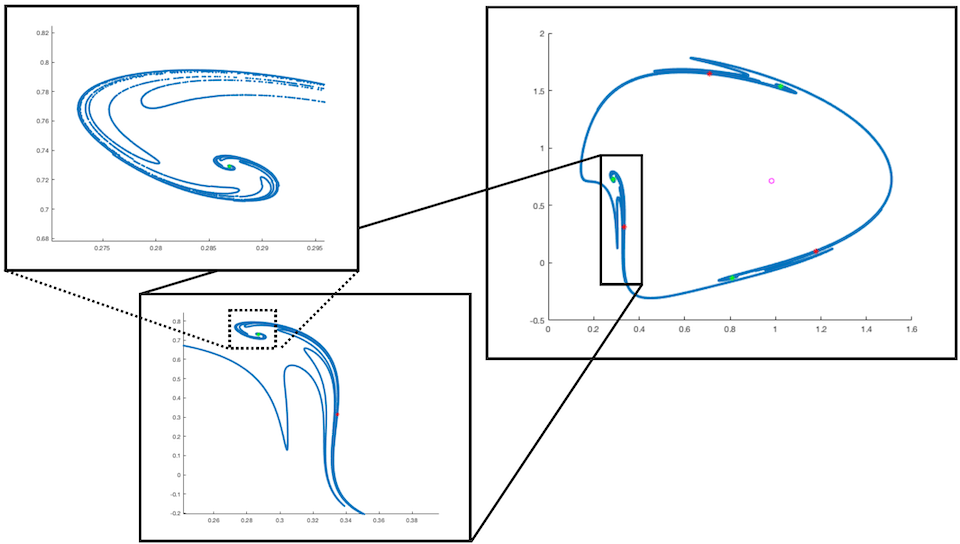}}
\caption{\textbf{Fractal structure of $W^u(\mathcal{Q}_2)$:}
at $\alpha=0.9321$, after the formation of the homoclinic tangency, 
the structure of the unstable manifold is much more complicated.  
}
\label{fig:horseShoe}
\end{figure}

\begin{figure}
\centering
\subfloat{\includegraphics[width = 2.1in]{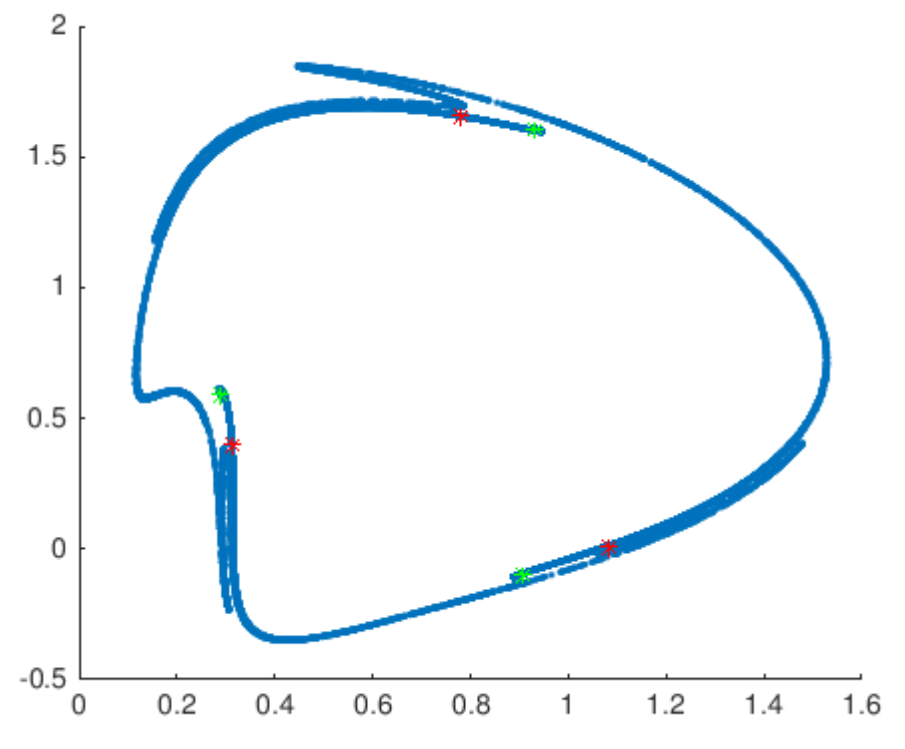}}
\subfloat{\includegraphics[width = 2.1in]{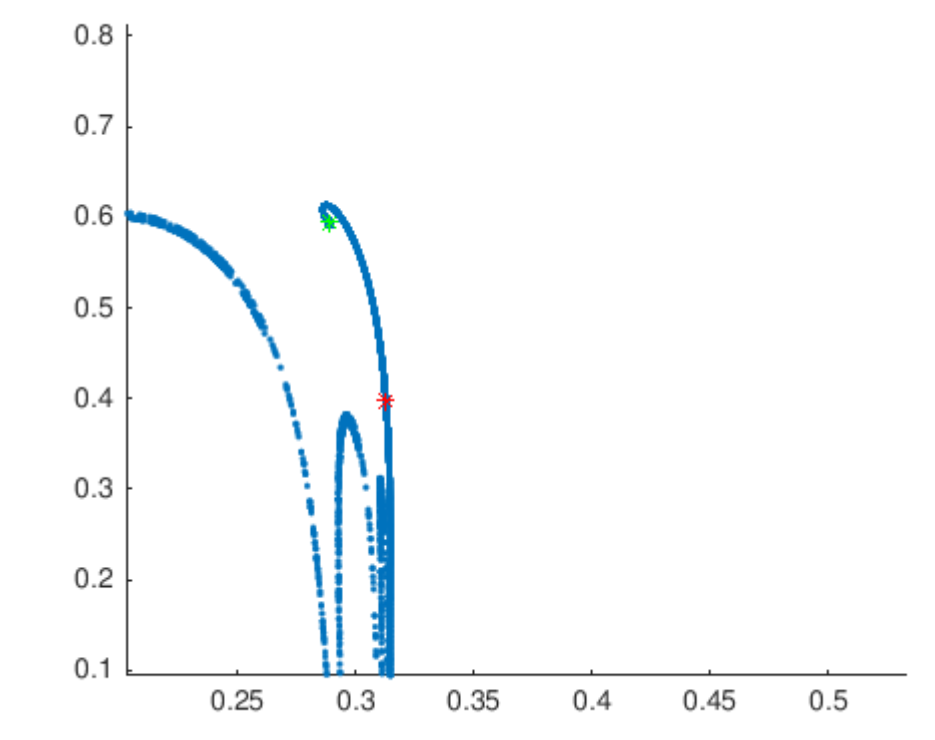}} 
\subfloat{\includegraphics[width = 2.1in]{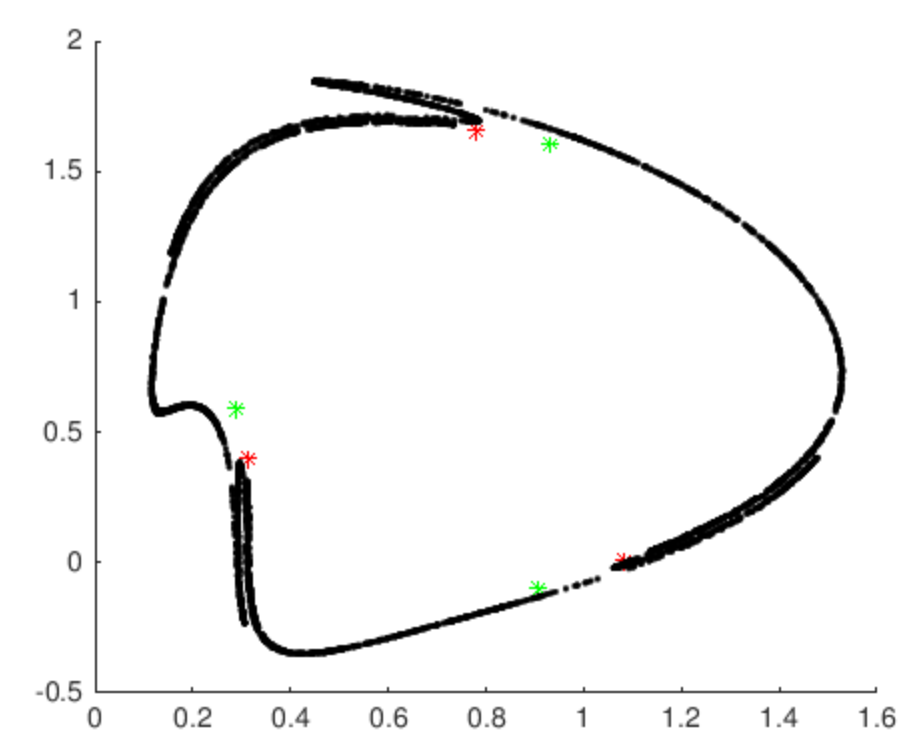}}
\caption{All three figures were plotted using $\alpha=0.95$. (Left) unstable manifold of the saddle period three cycle in the Poincar\'{e} section, colors have the same meaning as in Figure \ref{fig:preResTorus}. (Center) the fractal structure of $W^u(\mathcal{Q}_2)$ is no longer present going into the stable cycle. (Right) the torus $\mathcal{T}$ in phase space is destroyed and trajectories fall in a set with fractal dimension $2+d$, $d<1$. Both period three cycles have moved inside the invariant set (black). 
 }
\label{fig:attractor}
\end{figure}

Further numerical evidence for this claim is given in the three frames of Figure \ref{fig:attractor}.
 Here the bifurcation parameter is increased slightly further to $\alpha = 0.95$ so that the picture 
 opens up a little.  The frame on the right is obtained by iterating a large number of initial 
 conditions  until they numerically converge to the attractor, represented by the 
 black curve.  Note that the stable and saddle $3$-cycles (green and red collections of dots) appear to 
 have moved off the attractor as they do not touch the black curve.  
This indicates multi-stability in the system as the green orbit is itself attracting.  Moreover
the left frame shows the numerically computed unstable manifold $W^u(\mathcal{Q}_2)$, and it is clear by comparing 
the left with the right frame that -- while $W^u(\mathcal{Q}_2)$ is accumulating on the chaotic attractor 
union $\mathcal{Q}_1$, 
the unstable manifold of $\mathcal{Q}_2$  is no longer contained in the attractor.  This becomes even more clear
when we look again at the right frame of the figure and see no black line from the red dots to the green:
hence the attractor is no longer a resonant torus.  

Zooming in suggests in fact that the attractor is now a 
quite complicated shape, as illustrated in the middle frame of Figure \ref{fig:attractor}, and 
also in Figure \ref{fig:horseShoe}.
Indeed Figure \ref{fig:horseShoe} reveals the fractal structure of the unstable manifold of $\mathcal{Q}_2$.
The actual structure of the set is even more complicated than any picture can reveal, 
as results from topological dynamics imply that once there is a transverse homoclinic 
for $\mathcal{Q}_2$ the closure of the unstable manifold, which contains the attractor, 
is an indecomposable continuum \cite{MR908665,MR1277856}. 
We refer also the the work of 
\cite{MR968418} on the persistence of normally hyperbolic invariant manifolds in 
the absence of uniform rates and to the much more recent and constructive 
work of \cite{maciejContinum}.

\subsection{Visualizing the torus in phase space}

\begin{figure}[t!]
\centering
\includegraphics[width = 6in]{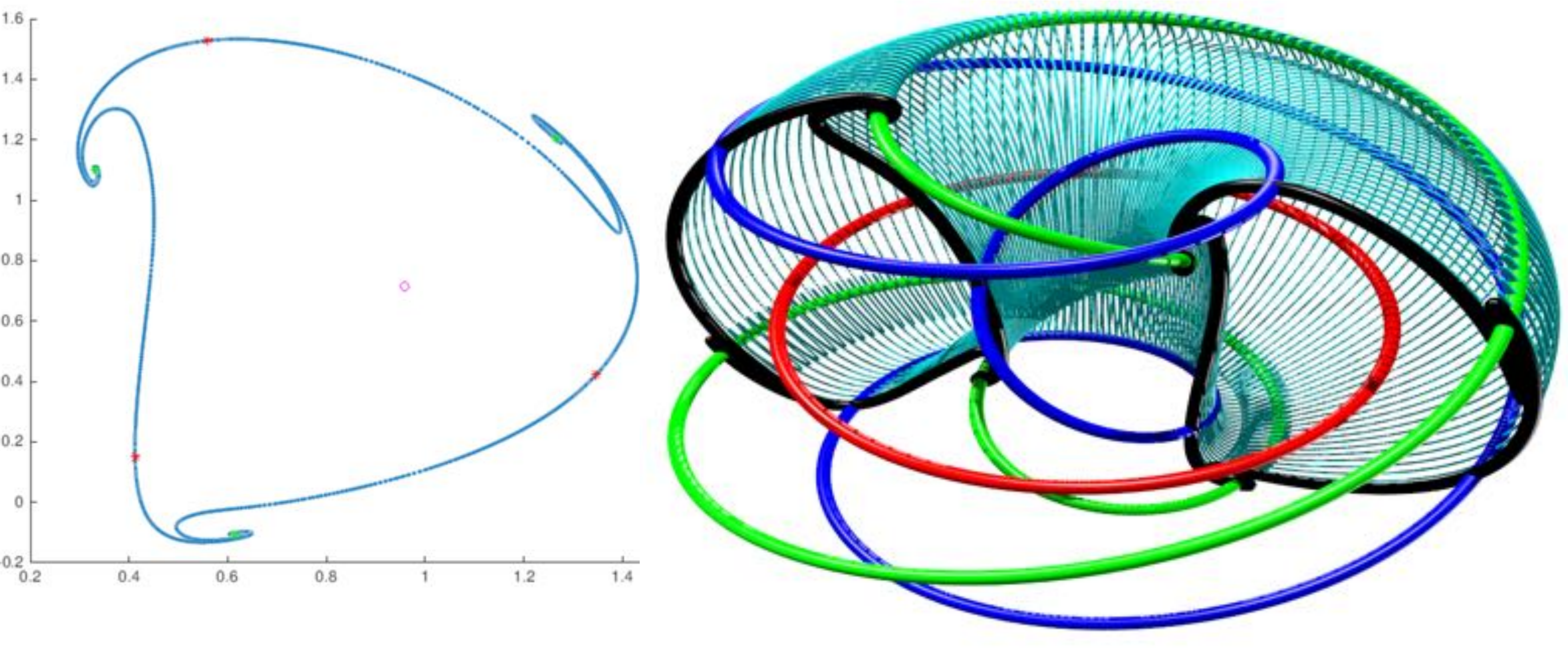}
\caption{\textbf{Cut-away at $\alpha = 0.85$:} The left frame recalls the 
invariant set in the Poincar\'{e} section when $\alpha = 0.85$, which is a resonant 
invariant circle formed by two period three cycles.  The right frame illustrates the corresponding 
invariant set in phase space.  The red curve is the repelling periodic orbit $\gamma$ which originally 
underwent the Neimark-Sacker bifurcation.  The green curve is the attracting periodic orbit $\gamma_1$
corresponding to the attracting period three cycle, while the blue curve is the saddle periodic orbit 
$\gamma_2$ corresponding  
to the saddle $3$-cycle in the Poincar\'{e} section. The half torus (colored in teal) is obtained
by advecting the section's invariant circle $\Gamma$ under the flow.  Since the invariant circle is formed by the 
unstable manifold of the saddle cycle, the teal surface represents $W^u(\gamma_2)$.
The resulting invariant set is a topological, but not smooth, invariant attracting 
resonant torus.} \label{fig:sideBySide085}
\end{figure}

\begin{figure}[t!]
\centering
\includegraphics[width = 6in]{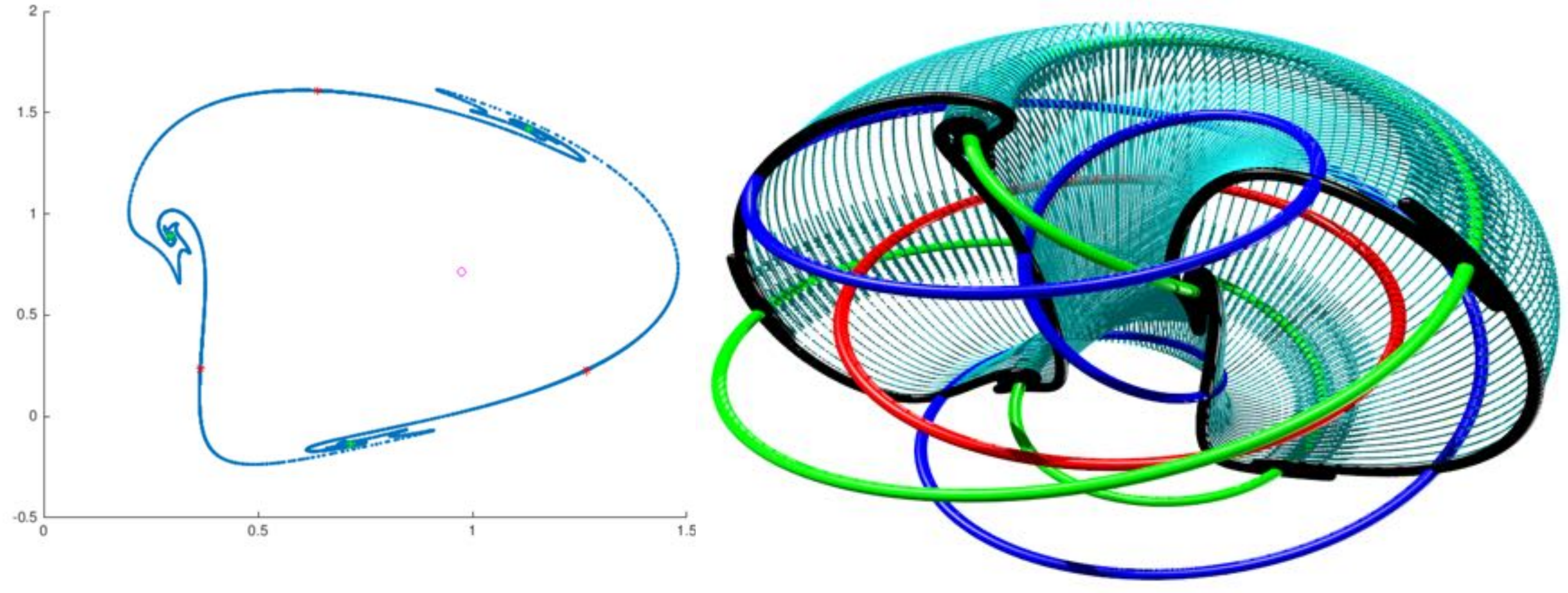}
\caption{\textbf{Cut-away at $\alpha = 0.9$:} same color scheme as in Figure 
\ref{fig:sideBySide085}. Complicated embedding of the
 attracting resonant invariant torus.} \label{fig:sideBySide09}
\end{figure}

\begin{figure}[t!]
\centering
\includegraphics[width = 6in]{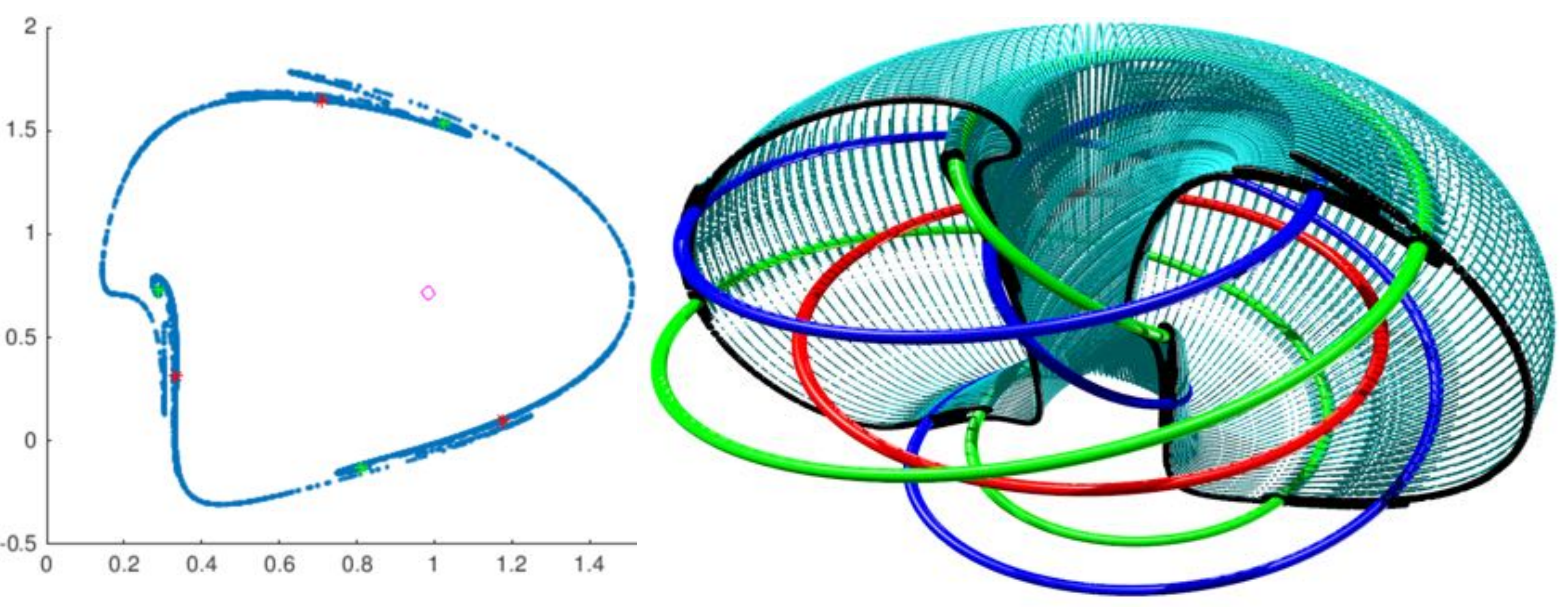}
\caption{\textbf{Cut-away at $\alpha = 0.929$:} 
same color scheme as in Figure 
\ref{fig:sideBySide085}.  The unstable manifold of the saddle type period 
three cycle in the Poincar\'{e} section side-by-side with the unstable 
manifold of $\gamma_2$ in phase space, just after the global bifurcation 
triggered by the homoclinic tangency at $\alpha_3$.  The unstable 
manifold accumulates on the union of the stable periodic orbit $\gamma_1$ and the 
chaotic attractor.} \label{fig:sideBySide09231}
\end{figure}

Studying the dynamics of a three dimensional system in a
two dimensional Poincar\'{e} section helps us to locate bifurcations 
of the invariant tori  by reducing them to invariant circles.
Nevertheless, it is still desirable to visualize dynamical structures
of the original system in the full phase space, 
and  to this end we provide several images which show side by side 
the invariant objects found in the Poincar\'{e} section and the corresponding 
invariant objects for the full Langford system.

See for example Figure \ref{fig:sideBySide085}.  The left frame illustrates the 
attracting invariant circle in the section
for the parameter value $\alpha = 0.85$.  The 
three red and green dots represent the saddle type and stable $3$-cycles
$\mathcal{Q}_2$ and $\mathcal{Q}_1$
respectively.  The magenta dot in the center of the frame represents
the repelling fixed point, while the blue curve is the unstable manifold of 
the saddle.  The blue curve is clearly absorbed into the basin of attraction 
of the stable $3$-cycle, forming a resonant invariant circle.  
The $3$-cycles $\mathcal{Q}_1, \mathcal{Q}_2$ give rise to periodic 
solutions in $\mathbb{R}^3$, which we denote by 
$\gamma_1$ and $\gamma_2$ respectively.

The right frame of the same figure illustrates the embedding in 
phase space of the same objects.  Here the red curve represents the repelling 
periodic orbit $\gamma$, green curve the attracting periodic orbit $\gamma_1$,
and the blue curve is the saddle periodic orbit $\gamma_2$.  
The unstable manifold of $\gamma_2$ accumulates at $\gamma_1$ forming the 
resonant torus.  Half of the torus is cut-away so that the skeleton given by 
the periodic orbits stands our clearly.   

Figures \ref{fig:sideBySide09} and 
\ref{fig:sideBySide09231} depict the same information at $\alpha = 0.9$ and 
$\alpha = 0.929$.  Taken together the three images provide much more insight 
into the structure of the invariant dynamics than can be gained by studying 
simulations of individual orbits like those illustrated in Figure 
\ref{fig:attractors}.

\section{Global dynamics}
\label{sec:globalD}
We now come to the second part of the present work, and study the dynamics not on 
but near the attracting invariant torus.  The remaining important features of the 
surrounding phase space are the equilibrium solutions on the $z-$ axis, and their invariant 
manifolds.  For this part of the study we abandon the Poincar\'{e} section and consider features 
of the three dimensional phase space.

\subsection{Accumulation of $W^u(p_0)$ on a component of the global attractor} \label{sec:accumulationOfWu}
One of the most important features in the phase space of the Langford system is the 
equilibrium point $p_0 \in \mathbb{R}^3$, which for all $\alpha \geq 0$ is located on the positive $z$ axis
and which has two dimensional unstable manifold and one dimensional stable manifold.
For all $\alpha \geq 0 $ the stable manifold of $p_0$ is a subset of the $z$-axis.
The two dimensional unstable manifold on the other hand is much more interesting.
All the calculations in this section are preformed using 
 the parameterization method/continuation scheme as discussed in 
the Section \ref{sec:parmMethod}.

\begin{figure}[t!]
\subfloat[$\alpha=0$]{\includegraphics[width = 3in]{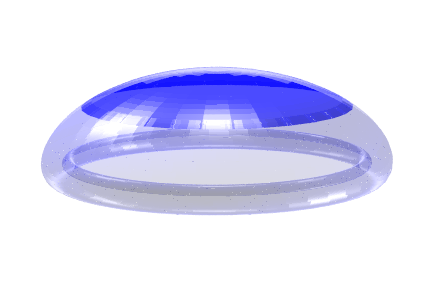}} \quad
\subfloat[$\alpha=0.6$]{\includegraphics[width = 3in]{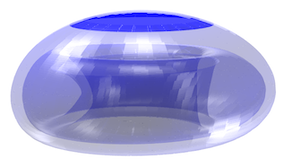}}\\
\subfloat[$\alpha=0.7$]{\includegraphics[width = 3in]{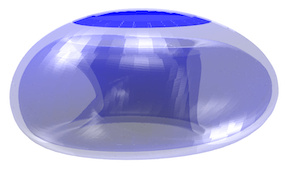}} \quad
\subfloat[$\alpha=0.8$]{\includegraphics[width = 3in]{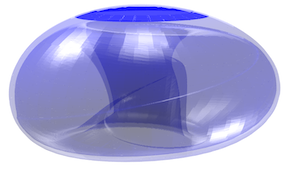}}\\
\subfloat[$\alpha=0.806$]{\includegraphics[width = 3in]{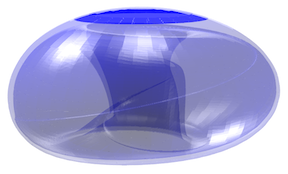}} \quad
\subfloat[$\alpha=0.9321$]{\includegraphics[width = 3in]{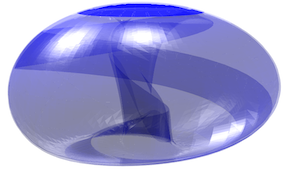}}
\caption{\textbf{Unstable manifold ``bubble'' for the saddle $p_0$}: 
for the indicated values of $\alpha$.    The computations suggest the 
existence of a periodic orbit which undergoes a Neimark-Sacker bifurcation.
The phase space is then dominated by the resulting 
smooth invariant torus.  The computations for higher 
$\alpha$ suggest that the smoothness 
of the torus may breakdown as 
$\alpha$ increases.} \label{fig:unstableBubbles1}
\end{figure}

\begin{figure}[t!]
\subfloat[$\alpha=0.95$]{\includegraphics[width = 3in]{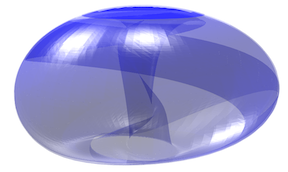}} \quad
\subfloat[$\alpha=1.1022$]{\includegraphics[width = 3in]{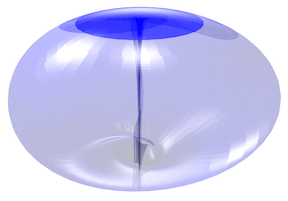}}
\caption{\textbf{Unstable manifold of $p_0$ after the saddle node bifurcation:}
(Left) $\alpha = 0.95$.  The bubble develops a ``stripe'' which is  due to 
the manifold folding over itself as it accumulates on the union of the 
chaotic attractor and the attracting periodic orbit.
(Right) $\alpha = 1.1022$, the quality of the bubble has changed dramatically.  It is 
more ``open'' and appears to accumulate 
on the $z$-axis.  
} 
\label{fig:unstableBubbles2}
\end{figure}

Recall that for $0 \leq \alpha \leq \alpha_4 \approx 0.9321697517861$, the point $p_0$ is the 
only equilibrium of the Langford system.  The manifold $W^u(p_0)$ is 
 illustrated in Figure \ref{fig:unstableBubbles1} for six such values of $\alpha$.
We see that for $\alpha < \alpha_1$ -- the value of the Neimark- Sacker bifurcation --
 $W^u(p_0)$ accumulates on the attracting periodic orbit $\gamma$ as seen in 
 frame (a) of Figure \ref{fig:unstableBubbles1}. The periodic orbit $\gamma$
appears to be the global attractor in this parameter range.

After the Neimark- Sacker bifurcation at $\alpha_1 \approx  0.697144898322973$ 
the periodic orbit $\gamma$ is repelling and $W^u(p_0)$ appears to accumulate
on the smooth attracting invariant torus $\mathcal{T}$ which was 
discussed at length in Section \ref{sec:global}.
This is seen in frames (b), (c) and (d) of Figure \ref{fig:unstableBubbles1}.  
Frames (d) and (e) illustrate the situation 
after the appearance of the attracting period three cycle in the Poincar\'{e} section
(see Figure \ref{fig:preResTorus}), and there is an 
attracting periodic orbit in phase space which we denote by 
$\gamma_1$. The system is bistable, with the 
attractor being the union of the invariant torus $\mathcal{T}$ and the
periodic orbit $\gamma_1$.  The manifold $W^u(p_0)$ accumulates on the 
union of these two objects -- a disjoint set -- and this is what introduces the 
rough folds in the embedding seen in Frames (d) and (e).

Frame (f) illustrates the embedding of $W^u(p_0)$ for $\alpha > \alpha_2$ but just before the 
global bifurcation at $\alpha = \alpha_4$ which destroys the torus.
Here the torus is resonant and only $C^0$, a fact which is clearly
reflected in the embedding of $W^u(p_0)$.

For  $\alpha > \alpha_4$ we are past the saddle node bifurcation, and $p_0$ is no longer
the unique equilibrium.  This has dramatic consequences for the global dynamics, and 
we illustrate the unstable manifold for two such parameter values in Figure 
\ref{fig:unstableBubbles2}.
Somewhere between $\alpha = 0.95$ (frame (a) of the figure) and $\alpha = 1.1022$
(frame (b)) something dramatic happens.  The change however is not easily understood
by looking only at $W^u(p_0)$, and we must consider the embedding of new invariant 
objects which appear only after the occurrence of the saddle node bifurcation.

\subsection{$W^s(p_1)$ as a separatrix} \label{sec:separatrix}
At $\alpha_4 \approx 0.9321697517861$ the system undergoes a saddle node bifurcation, 
resulting in the appearance of two new equilibrium solutions denoted 
$p_1$ and $p_2$.  For all $\alpha > \alpha_4$ the point $p_2$ is a stable 
equilibrium, making $p_2$ a new component of the global attractor.  
The equilibrium solution at $p_1$ on the other hand is a saddle, with 
one dimensional unstable manifold on the $z$-axis and a two dimensional 
stable manifold associated with a pair of complex conjugate eigenvalues.  

\begin{figure}[t!]
\centering
\subfloat{\includegraphics[width = 2.9in]{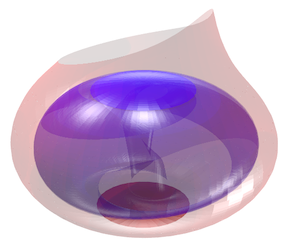}}\quad
\subfloat{\includegraphics[width = 2.9in]{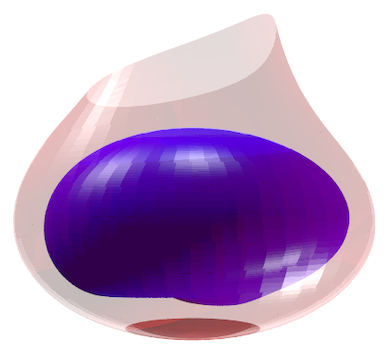}}
\caption{\textbf{2D stable and unstable manifolds of equilibria of $p_0$ and $p_1$:}.
At $\alpha = 0.95$ we note that 
 $W^u(p_0)$ (blue) and $W^s(p_1)$ (red) do not intersect at all.
We also remark that $p_2$ (not shown) 
 is below $p_1$ and is an attracting equilibrium point. 
 The attracting torus (or ``torus-like'' chaotic attractor) 
 is inside the ``bubble'' formed by these stable/unstable manifolds.} \label{fig:separatrix}
\end{figure}

Our simulations suggest that 
for some range of $ \alpha > \alpha_4 $,  
the two dimensional invariant manifold
$W^s(p_1)$ is a separatrix for the basin of attraction of $p_2$ and 
the attractor onto which $W^u(p_0)$ accumulates.  
In this sense, $W^s(p_1)$ forms a kind of ``bubble'', where inside the 
bubble we have an attractor comprised of either 
the resonant torus $\mathcal{T}$ or the chaotic 
set which appears after the destruction of $\mathcal{T}$.  The 
inside of the bubble is the basin of attraction of this 
attractor, and the outside of the bubble is the basin of the stable equilibrium
$p_2$.  The situation is illustrated in Figure $\ref{fig:separatrix}$.

\begin{figure}
\subfloat{\includegraphics[width = 3in]{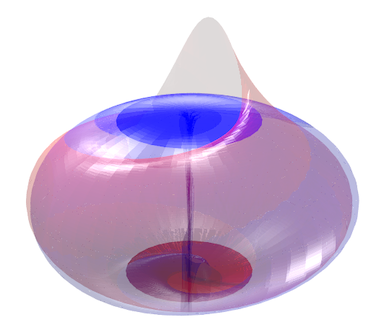}}\quad
\subfloat{\includegraphics[width = 3in]{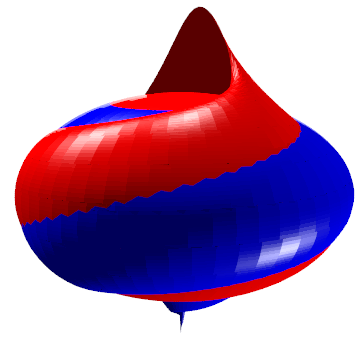}}
\caption{\textbf{2D stable and unstable manifolds of equilibria of $p_0$ and $p_1$
 for $\alpha = 1.1022$}. Here we see that 
 $W^u(p_0)$ (blue) and $W^s(p_1)$ (red) appear to intersect  
 transversely.  The intersection curves are then heteroclinic orbits from $p_0$ to 
 $p_1$.
 The unstable manifold accumulates on the $z$-axis, as seen in the transparency on the 
 left.  The frame on the right suggests that the unstable manifold enters the 
 basin of attraction of $p_2$.  In fact, for $\alpha = 1.1022$ it seems that 
 $p_2$ is the unique attractor.  
 }
\label{fig:heteroclinics}
\end{figure}

\begin{figure}
\includegraphics[width = 6in]{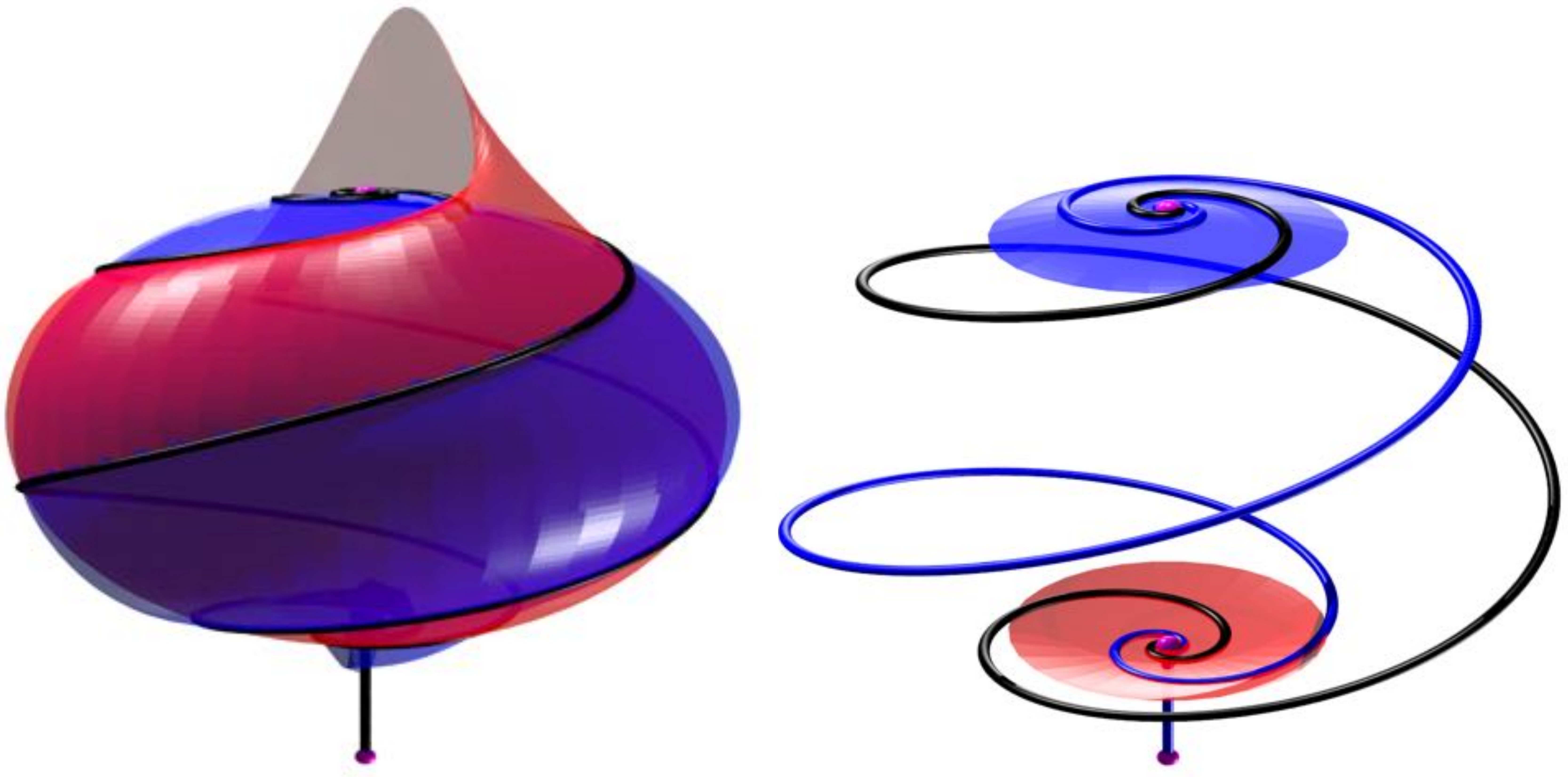}
\caption{\textbf{Heteroclinic connections from $p_0$ to $p_1$
 for $\alpha = 1.1022$}. Encouraged by the apparent transverse intersections
seen in Figure \ref{fig:heteroclinics}, we locate heteroclinic orbit segments 
starting on $W^u(p_0)$ and terminating on $W^s(p_1)$ by applying a 
Newton scheme to the boundary value problem describing the 
segments.  Observe that the heteroclinic orbit segments located are much smoother than 
the apparent intersection seen in Figure \ref{fig:heteroclinics}.  The apparent 
irregularity of the intersection is due to the fact that we compute piecewise 
linear triangulations of the fundamental domain and its iterates. 
 }
\label{fig:theConnections}
\end{figure}

\subsection{Heteroclinic intersections and the loss of bistability}
Studying $W^u(p_0)$ and $W^s(p_1)$ reveals yet another global 
bifurcation which dramatically alters the phase space dynamics
of the system.  
It appears that for some $\alpha_5 \approx 1.05$ these manifolds develop a tangency,
and that after this tangency there are transverse heteroclinic connections from 
$p_0$ to $p_1$.   
The situation is illustrated in Figures \ref{fig:heteroclinics}
and \ref{fig:theConnections}, where we see the transverse intersections of the manifolds and 
the resulting heteroclinic connections respectively.  

Once intersections appear between $W^u(p_0)$ and $W^s(p_1)$, the latter ceases to 
function as a separatrix, and orbits can pass from inside the bubble to outside.
The equilibrium $p_2$ remains attracting and its basin appears now to extend 
into the inside of the bubble.  
Indeed our numerical simulations suggest that for $\alpha > 1$, that is after the 
formation of intersections between the unstable/stable manifolds of $p_0$
and $p_1$, the equilibrium $p_2$ becomes the global attractor for the 
system.  That is, all orbits which start inside the bubble eventually 
accumulate there.  This finally explains the ``openness'' of 
Figure \ref{fig:unstableBubbles2} (b) remarked upon earlier:  
the occurrence of the heteroclinic  tangency 
between $W^u(p_0)$ and $W^s(p_1)$ appears to destroy the attractor which had 
previously dominated the dynamics inside the bubble.

\section{Conclusions and Discussion} \label{sec:conclusions}
We will summarize the results of the present work  
by sketching the main features of the global dynamics of the Langford system \eqref{eq:1}
as suggested by our analysis.  
First recall the main local and global bifurcations undergone by the system.
\begin{itemize}
\item At $\alpha_1 \approx  0.697144898322973$ the periodic orbit $\gamma$ 
undergoes a Neimark-Sacker
bifurcation.  This is a local bifurcation of $\gamma$ which results in 
the appearance of the invariant torus $\mathcal{T}$.
\item At $\alpha_2 \approx 0.823$ the invariant torus $\mathcal{T}$ develops a resonance.
After this $\mathcal{T}$  is the union of two periodic orbits $\gamma_2$ (saddle stability),
and $\gamma_1$ (attracting), and the unstable manifold of $\gamma_2$. The resonance
is triggered by the collision of a saddle periodic orbit with the invariant 
 torus.  Since the torus is an attractor before and after the collision, this bifurcation 
 involves no change in the stability of $\mathcal{T}$ and is hence a global bifurcation.
\item At $\alpha_3 \approx 0.925$ there is a global 
bifurcation triggered by the formation of a tangency between $W^s(\gamma_1)$
and $W^u(\gamma_1)$.  
\item At $\alpha_4 \approx 0.9321697517861$ there is a saddle node bifurcation resulting in the 
appearance of the equilibrium points $p_1$ (saddle-focus stability) and $p_2$ (attracting).
This is a local fold bifurcation for the equilibrium points.  
\item At $\alpha_5 \approx 1$ a global bifurcation is triggered by 
the development of a tangency between $W^u(p_0)$ and 
$W^s(p_1)$.    
\end{itemize}

Between the bifurcation values listed above, we conjecture
based on our numerical studies that the 
system has the following properties.
\begin{conjecture}[Sketch of the global dynamics] \label{conj:main}
The flow generated by the Langford vector field,
given in Equation \eqref{eq:1}, has that;
\begin{enumerate}
\item For $0 < \alpha < \alpha_1$  the periodic orbit $\gamma$ is the global attractor.  
\item For $\alpha_1 < \alpha < \alpha_2$ the global attractor is either $\mathcal{T}$
or $\mathcal{T} \cup \gamma_1$.  
\item For $\alpha_2 < \alpha < \alpha_3$ the global attractor is $\mathcal{T}$. 
\item For $\alpha_3 < \alpha < \alpha_4$ the global attractor is 
$\mathcal{T} \cup p_2$.
\item For $\alpha_4 < \alpha < \alpha_5$ there is multi-stability.  The global attractor is 
comprised of at least the components $\tilde{\mathcal{T}}$ - the chaotic attractor appearing 
after the break-up of the invariant torus, the attracting periodic orbit $\gamma_2$, and
the attracting equilibrium solution $p_2$.
\item For $\alpha > \alpha_5$ the equilibrium solution $p_2$ is the global attractor.
\item For $0 < \alpha < \alpha_3$ the unstable manifold $W^u(p_0)$ accumulates 
on the global attractor, which is either $\gamma$ (until $\alpha = \alpha_1$) or $\mathcal{T}$.
\item For $\alpha_3 < \alpha < \alpha_5$ the stable manifold $W^s(p_1)$ is 
a separatrix.  The basin of attraction of $p_2$ is outside the bubble formed by 
 $W^s(p_1)$.  
\end{enumerate}
\end{conjecture}
It is essential to stress that the eight points above are still just conjectures,
however well informed. The numerical work carried out in the present work is not 
sufficient to rule out other components of the global attractor, for example other attracting 
periodic orbits near the resonant torus or the chaotic attractor.  This point is elaborated on below.   

\begin{figure}[t!]
\centering
\includegraphics[width = 6.5in]{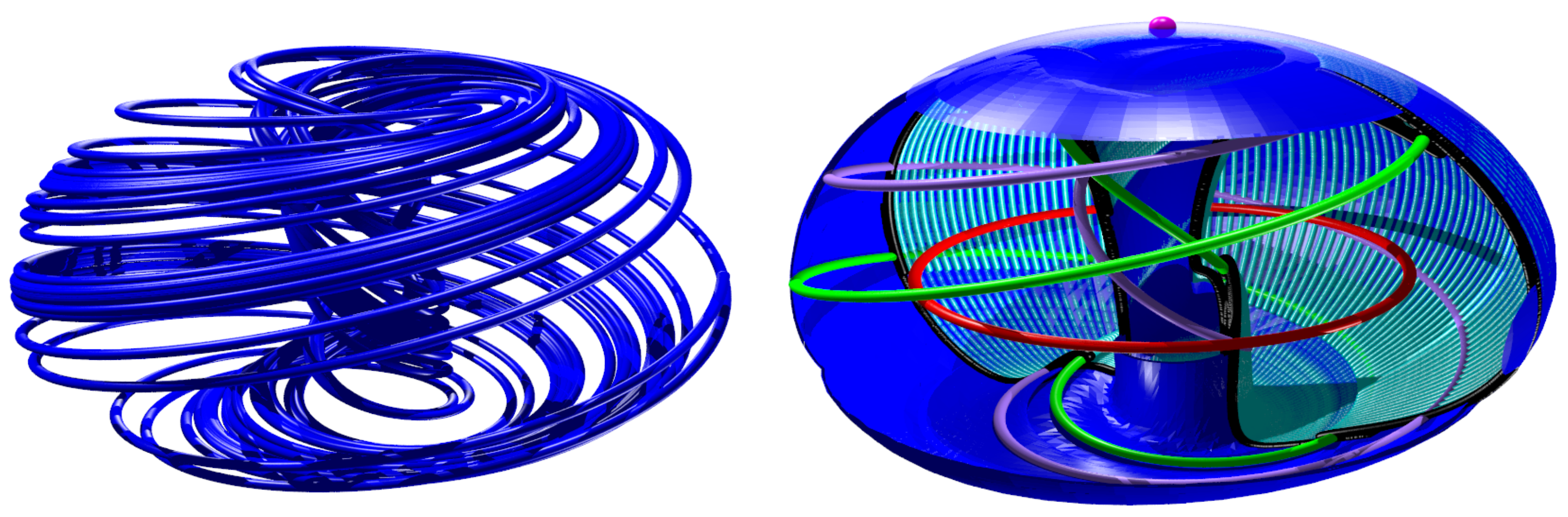}
\caption{\textbf{Visualization phase space structure: simulation versus
invariant manifolds at $\alpha = 0.9321$:} Left: direct simulation of an
initial condition for roughly one hundred time units. Right: the equilibrium 
solution $p_0$ (magenta dot), its unstable manifold (blue surface), 
the resonant torus comprised of a stable periodic orbit $\gamma_1$ 
(green curve), the saddle periodic orbit $\gamma_2$ (purple curve), 
and its unstable manifold (green torus).  Also shown is the 
repelling periodic orbit $\gamma$ (red curve).  Most of these objects have
unstable directions and are impossible to locate by direct simulation. 
Even the attracting resonant torus is very difficult to ``fill in'' by just 
simulating the system.
} \label{fig:finalPic}
\end{figure}

It is also worth remarking that softer sorts of conclusions are encapsulated in
the paper's many figures.  The deliberate calculations and three dimensional 
renderings of invariant manifolds throughout our work provide more
delicate insights into the dynamics of the system than are obtained by 
straight forward simulations of ensembles of initial conditions.    
As a final illustration of this point we give in Figure \ref{fig:finalPic}
a side by side comparison 
of the results obtained using the methodology of the present work with 
the results obtained by direct integration, for the parameter value $\alpha = 0.9321$.
Simulation results cannot illuminate the full attractor, as numerical integrations
will never reveal the unstable periodic orbit $\gamma_2$ which lies inside 
the invariant torus $\Gamma$.

More generally, 
since the Langford system is derived by truncating the normal form of a cusp-Hopf
singularity, we expect qualitatively similar behavior in systems exhibiting  
this bifurcation.   
An interesting topic for future research would be to repeat the numerical analysis performed
in the present work for other modifications of the normal form. 
For example one could perturb the system in such a way 
that the $z$-axis is no longer invariant.  Or, one could modify the system so that 
the Neimark-Sacker bfurcation is subcritical rather than supercritical.

 As remarked already in 
\cite{MR821035} (and in the introduction of the present work) complex 
dynamics are often generated by interactions between equilibrium and periodic 
solutions.  The fact that the Langford system is close to  
a simultaneous cusp-Hopf bifurcation is precisely what provides
the multiple equilibrium solutions in close proximity to a limit cycle.  This is the 
basic mechanism organizing many of the interesting dynamical
phenomena discussed 
in the present work.

The normal form unfolding a pitchfork-Hopf bifurcation was also 
studied by Langford in \cite{MR591900}, and it would be a nice 
project to apply the methods of the present work to this system, or 
to systems derived from the fold-Hopf bifurcation.  
Other important normal forms are discussed for example in 
the woks of \cite{MR508917,MR950168,MR633825}.

Another interesting topic of future work would be to prove -- possibly with computer
assistance -- as much of Conjecture \ref{conj:main} as possible.  
For example, the following Theorem is found in the Author's work with Maciej Capinski 
\cite{maciej_us_invariantTori}.
\begin{theorem}[Existence of a $C^0$ invariant torus]
For $\alpha = 0.85$ Equation \eqref{eq:1} has a $C^0$ resonant invariant torus, which is not even globally 
Lipschitz much less $C^1$.  More precisely, the torus contains exactly two periodic orbits 
-- one attracting and the other a saddle.   The Floquet exponents of the attracting periodic orbit are complex conjugates.
The saddle periodic orbit has one stable and one unstable Floquet exponent.  
The one dimensional unstable manifold of the saddle periodic orbit is completely 
captured in the basin of attraction of the stable periodic orbit, so that torus is the 
union of the stable periodic orbit, the saddle periodic orbit, and the unstable manifold of the saddle. 
\end{theorem}
The proof of this theorem is based on the  techniques developed in 
\cite{AMS, MR3281845,BDLM} for validating bounds on local manifold parameterizations and 
computer assisted proofs for 
heteroclinic connections, the methods developed in \cite{cnLohner, MR1961956, MR2262261, MR1816799}
for rigorous integration of vector fields and computer assisted proof in Poincar\'{e} sections, and the methods of 
\cite{MR3443692, MR3032848} for obtaining validated error bounds on stable/unstable manifolds in 
Poincar\'{e} sections. This one theorem provides a glimpse of what could be 
accomplished in this and similar systems by constructing computer assisted arguments.

For example, the 
techniques developed in \cite{MR1961956,MR2262261,MR2271217,MR2173545}
could be used to prove the existence of the global bifurcations observed above.  
Using the techniques developed in \cite{MR2533624,MR3504856}, it should be possible to 
study in a mathematically rigorous way the global attractor of the Langford
system over a large parameter range, and verify and/or clarify many of the claims
of Conjecture \ref{conj:main}.  For example parts (1), (2), (3), (4), and (5)
appear to us susceptible to this kind of analysis.  
Combining the techniques of the references just cited with the mathematically 
rigorous methods for computing stable/unstable invariant manifold 
atlases developed in \cite{manifoldPaper1} could provide means of verifying parts
(6) and (7) of the conjecture.  We also remark that the recent work of \cite{maciejContinum}
could be applied to give computer assisted proofs for the chaotic attractor after 
breakdown.  

Another project could be to combine the parameterization method for hyperbolic
invariant tori developed in 
\cite{MR3713933, MR3713932} with the methods of computer assisted proof 
developed in \cite{alexCAPKAM} to prove the existence of the invariant tori 
studied in the present work for $\alpha < \alpha_3$.  That is, to study the 
tori before the onset of resonance.  
As mentioned already in \cite{maciej_us_invariantTori}, the methods of computer assisted
proof developed there appear to struggle in the parameter range $\alpha_1 < \alpha < \alpha_2$
due to the apparent lack of uniform contraction rates near the torus.
If arguments like the ones suggested in this paragraph could
succeed for the Langford system, they then could also be 
applied to systems coming from other important normal forms.

\section{Acknowledgments}
The authors would like to thank Jordi-Llu\'is 
Figueras, Maciej Capi\'nski, and Vincent Naudot for many helpful suggestions and 
invaluable insights. We owe also a special thanks to Takahito Mitsui for bringing the 
paper of Langford \cite{MR821035} to our attention after reading a much rougher early 
version of this manuscript.  The published version of the manuscript benefitted 
greatly from the suggestions of two anonymous referees.
The second author was partially supported by NSF grant 
DMS-1813501.

{\footnotesize
\bibliographystyle{unsrt}
\bibliography{papers}

\begin{thebibliography}{100}

\bibitem{MR0132256}
Ju.~I. Ne\u{\i}mark.
\newblock Some cases of the dependence of periodic motions on parameters.
\newblock {\em Dokl. Akad. Nauk SSSR}, 129:736--739, 1959.

\bibitem{MR2615427}
Robert~John Sacker.
\newblock {\em On invariant surfaces and bifurcation of periodic solutions of
  ordinary differential equations}.
\newblock ProQuest LLC, Ann Arbor, MI, 1964.
\newblock Thesis (Ph.D.)--New York University.

\bibitem{MR1005055}
Seung-hwan Kim, R.~S. MacKay, and J.~Guckenheimer.
\newblock Resonance regions for families of torus maps.
\newblock {\em Nonlinearity}, 2(3):391--404, 1989.

\bibitem{MR1115870}
C.~Baesens, J.~Guckenheimer, S.~Kim, and R.~S. MacKay.
\newblock Three coupled oscillators: mode-locking, global bifurcations and
  toroidal chaos.
\newblock {\em Phys. D}, 49(3):387--475, 1991.

\bibitem{MR709899}
Kunihiko Kaneko.
\newblock Transition from torus to chaos accompanied by frequency lockings with
  symmetry breaking. {I}n connection with the coupled-logistic map.
\newblock {\em Progr. Theoret. Phys.}, 69(5):1427--1442, 1983.

\bibitem{MR834186}
V.~S. Afraimovich and L.~P. Shil\'{n}ikov.
\newblock Invariant two-dimensional tori, their breakdown and stochasticity.
\newblock In {\em Methods of the qualitative theory of differential equations},
  pages 3--26, 164. Gor\'{k}ov. Gos. Univ., Gorki, 1983.

\bibitem{MR906312}
D.~V. Turaev and L.~P. Shil\'{n}ikov.
\newblock Bifurcations of quasi-attractors of torus-chaos.
\newblock In {\em Mathematical mechanisms of turbulence ({R}ussian)}, pages
  113--121, iii. Akad. Nauk Ukrain. SSR, Inst. Mat., Kiev, 1986.

\bibitem{MR3095277}
Renato~C. Calleja, Alessandra Celletti, and Rafael de~la Llave.
\newblock Local behavior near quasi-periodic solutions of conformally
  symplectic systems.
\newblock {\em J. Dynam. Differential Equations}, 25(3):821--841, 2013.

\bibitem{MR3062760}
Renato~C. Calleja, Alessandra Celletti, and Rafael de~la Llave.
\newblock A {KAM} theory for conformally symplectic systems: efficient
  algorithms and their validation.
\newblock {\em J. Differential Equations}, 255(5):978--1049, 2013.

\bibitem{MR3713933}
Marta Canadell and \`Alex Haro.
\newblock Computation of quasiperiodic normally hyperbolic invariant tori:
  rigorous results.
\newblock {\em J. Nonlinear Sci.}, 27(6):1869--1904, 2017.

\bibitem{MR783349}
Alain Chenciner.
\newblock Bifurcations de points fixes elliptiques. {I}. {C}ourbes invariantes.
\newblock {\em Inst. Hautes Etudes Sci. Publ. Math.}, (61):67--127, 1985.

\bibitem{MR784530}
A.~Chenciner.
\newblock Bifurcations de points fixes elliptiques. {II}. {O}rbites
  p\'eriodiques et ensembles de {C}antor invariants.
\newblock {\em Invent. Math.}, 80(1):81--106, 1985.

\bibitem{MR932134}
Alain Chenciner.
\newblock Bifurcations de points fixes elliptiques. {III}. {O}rbites
  p\'{e}riodiques de ``petites'' p\'{e}riodes et \'{e}limination r\'{e}sonnante
  des couples de courbes invariantes.
\newblock {\em Inst. Hautes Etudes Sci. Publ. Math.}, (66):5--91, 1988.

\bibitem{MR1285950}
R.~S. MacKay.
\newblock Transport in {$3$}{D} volume-preserving flows.
\newblock {\em J. Nonlinear Sci.}, 4(4):329--354, 1994.

\bibitem{MR701669}
Kunihiko Kaneko.
\newblock Similarity structure and scaling property of the period-adding
  phenomena.
\newblock {\em Progr. Theoret. Phys.}, 69(2):403--414, 1983.

\bibitem{MR845031}
Kunihiko Kaneko.
\newblock {\em Collapse of tori and genesis of chaos in dissipative systems}.
\newblock World Scientific Publishing Co., Singapore, 1986.

\bibitem{MR2241302}
Frank Schilder, Werner Vogt, Stephan Schreiber, and Hinke~M. Osinga.
\newblock Fourier methods for quasi-periodic oscillations.
\newblock {\em Internat. J. Numer. Methods Engrg.}, 67(5):629--671, 2006.

\bibitem{MR3309008}
Marta Canadell and \`Alex Haro.
\newblock Parameterization method for computing quasi-periodic reducible
  normally hyperbolic invariant tori.
\newblock In {\em Advances in differential equations and applications},
  volume~4 of {\em SEMA SIMAI Springer Ser.}, pages 85--94. Springer, Cham,
  2014.

\bibitem{MR3713932}
Marta Canadell and \`Alex Haro.
\newblock Computation of quasi-periodic normally hyperbolic invariant tori:
  algorithms, numerical explorations and mechanisms of breakdown.
\newblock {\em J. Nonlinear Sci.}, 27(6):1829--1868, 2017.

\bibitem{MR821035}
W.~F. Langford.
\newblock Numerical studies of torus bifurcations.
\newblock In {\em Numerical methods for bifurcation problems ({D}ortmund,
  1983)}, volume~70 of {\em Internat. Schriftenreihe Numer. Math.}, pages
  285--295. Birkh\"{a}user, Basel, 1984.

\bibitem{MR0462175}
V.~S. Afrauimovic, V.~V. Bykov, and L.~P. Silnikov.
\newblock The origin and structure of the {L}orenz attractor.
\newblock {\em Dokl. Akad. Nauk SSSR}, 234(2):336--339, 1977.

\bibitem{MR1733750}
V.~I. Arnold, V.~S. Afrajmovich, Yu.~S. Ilyashenko, and L.~P. Shilnikov.
\newblock {\em Bifurcation theory and catastrophe theory}.
\newblock Springer-Verlag, Berlin, 1999.
\newblock Translated from the 1986 Russian original by N. D. Kazarinoff,
  Reprint of the 1994 English edition from the series Encyclopaedia of
  Mathematical Sciences [{{\i}t Dynamical systems. V}, Encyclopaedia Math.
  Sci., 5, Springer, Berlin, 1994; MR1287421 (95c:58058)].

\bibitem{MR1237641}
Jacob Palis and Floris Takens.
\newblock {\em Hyperbolicity and sensitive chaotic dynamics at homoclinic
  bifurcations}, volume~35 of {\em Cambridge Studies in Advanced Mathematics}.
\newblock Cambridge University Press, Cambridge, 1993.
\newblock Fractal dimensions and infinitely many attractors.

\bibitem{MR2471925}
V.~Araujo, M.~J. Pacifico, E.~R. Pujals, and M.~Viana.
\newblock Singular-hyperbolic attractors are chaotic.
\newblock {\em Trans. Amer. Math. Soc.}, 361(5):2431--2485, 2009.

\bibitem{Cabre1}
X.~Cabr{\'e}, E.~Fontich, and R.~de~la Llave.
\newblock The parameterization method for invariant manifolds. {I}. {M}anifolds
  associated to non-resonant subspaces.
\newblock {\em Indiana Univ. Math. J.}, 52(2):283--328, 2003.

\bibitem{Cabre2}
X.~Cabr{\'e}, E.~Fontich, and R.~de~la Llave.
\newblock The parameterization method for invariant manifolds. {II}.
  {R}egularity with respect to parameters.
\newblock {\em Indiana Univ. Math. J.}, 52(2):329--360, 2003.

\bibitem{Cabre3}
X.~Cabr{\'e}, E.~Fontich, and R.~de~la Llave.
\newblock The parameterization method for invariant manifolds. {III}.
  {O}verview and applications.
\newblock {\em J. Differential Equations}, 218(2):444--515, 2005.

\bibitem{Stone05}
Z.B. Stone and H.A. Stone.
\newblock Imaging and quantifying mixing in a model droplet micromixer.
\newblock {\em Phys. Fluids}, 17:063103, 2005.
\newblock \url{https://doi.org/10.1063/1.1929547}.

\bibitem{MR1704974}
K.~E. Lenz, H.~E. Lomel\'{i}, and J.~D. Meiss.
\newblock Quadratic volume preserving maps: an extension of a result of
  {M}oser.
\newblock {\em Regul. Chaotic Dyn.}, 3(3):122--131, 1998.
\newblock J. Moser at 70 (Russian).

\bibitem{MR2481277}
H.~R. Dullin and J.~D. Meiss.
\newblock Quadratic volume-preserving maps: invariant circles and bifurcations.
\newblock {\em SIAM J. Appl. Dyn. Syst.}, 8(1):76--128, 2009.

\bibitem{MR2259296}
Shawn~C. Shadden, John~O. Dabiri, and Jerrold~E. Marsden.
\newblock Lagrangian analysis of fluid transport in empirical vortex ring
  flows.
\newblock {\em Phys. Fluids}, 18(4):047105, 11, 2006.

\bibitem{MR880159}
Takashi Matsumoto, Leon~O. Chua, and Ryuji Tokunaga.
\newblock Chaos via torus breakdown.
\newblock {\em IEEE Trans. Circuits and Systems}, 34(3):240--253, 1987.

\bibitem{MR1488520}
O.~Sosnovtseva and E.~Mosekilde.
\newblock Torus destruction and chaos-chaos intermittency in a commodity
  distribution chain.
\newblock {\em Internat. J. Bifur. Chaos Appl. Sci. Engrg.}, 7(6):1225--1242,
  1997.

\bibitem{MR3435117}
Taoufik Bakri, Yuri~A. Kuznetsov, and Ferdinand Verhulst.
\newblock Torus bifurcations in a mechanical system.
\newblock {\em J. Dynam. Differential Equations}, 27(3-4):371--403, 2015.

\bibitem{MR3279518}
Taoufik Bakri and Ferdinand Verhulst.
\newblock Bifurcations of quasi-periodic dynamics: torus breakdown.
\newblock {\em Z. Angew. Math. Phys.}, 65(6):1053--1076, 2014.

\bibitem{Vadim}
Vadim~S. Anishchenko, Vladimir Astakhov, Alexander Neiman, Tatjana Vadivasova,
  and Lutz Schimansky-Geier.
\newblock {\em Nonlinear dynamics of chaotic and stochastic systems}.
\newblock Springer Series in Synergetics. Springer, Berlin, second edition,
  2007.
\newblock Tutorial and modern developments.

\bibitem{aizawaVideo}
Arash Mohammadi.
\newblock {T}he {A}izawa {A}ttractor.
\newblock \url{https://www.youtube.com/watch?v=RBqbQUu-p00}, November 2017.

\bibitem{bridges2018:491}
Michael Gagliardo.
\newblock 3d printing chaos.
\newblock In Carlo~S\'equin Eve~Torrence, Bruce~Torrence and Krist\'of
  Fenyvesi, editors, {\em Proceedings of Bridges 2018: Mathematics, Art, Music,
  Architecture, Education, Culture}, pages 491--494, Phoenix, Arizona, 2018.
  Tessellations Publishing.
\newblock Available online at
  \url{http://archive.bridgesmathart.org/2018/bridges2018-491.pdf}.

\bibitem{chaoticAtmsopheres}
\url{http://chaoticatmospheres.com/mathrules-strange-attractors}.
\newblock ``{S}trange {A}ttractors.'' {C}haotic {A}tmospheres.

\bibitem{Haro1}
A.~Haro and R.~de~la Llave.
\newblock A parameterization method for the computation of invariant tori and
  their whiskers in quasi-periodic maps: rigorous results.
\newblock {\em J. Differential Equations}, 228(2):530--579, 2006.

\bibitem{Haro2}
{\`A}.~Haro and R.~de~la Llave.
\newblock A parameterization method for the computation of invariant tori and
  their whiskers in quasi-periodic maps: numerical algorithms.
\newblock {\em Discrete Contin. Dyn. Syst. Ser. B}, 6(6):1261--1300
  (electronic), 2006.

\bibitem{MR2299977}
A.~Haro and R.~de~la Llave.
\newblock A parameterization method for the computation of invariant tori and
  their whiskers in quasi-periodic maps: explorations and mechanisms for the
  breakdown of hyperbolicity.
\newblock {\em SIAM J. Appl. Dyn. Syst.}, 6(1):142--207 (electronic), 2007.

\bibitem{MR3118249}
Gemma Huguet and Rafael de~la Llave.
\newblock Computation of limit cycles and their isochrons: fast algorithms and
  their convergence.
\newblock {\em SIAM J. Appl. Dyn. Syst.}, 12(4):1763--1802, 2013.

\bibitem{MR2551254}
Antoni Guillamon and Gemma Huguet.
\newblock A computational and geometric approach to phase resetting curves and
  surfaces.
\newblock {\em SIAM J. Appl. Dyn. Syst.}, 8(3):1005--1042, 2009.

\bibitem{maximePOman}
J.~D. Mireles~James and Maxime Murray.
\newblock Chebyshev-{T}aylor parameterization of stable/unstable manifolds for
  periodic orbits: implementation and applications.
\newblock {\em Internat. J. Bifur. Chaos Appl. Sci. Engrg.}, 27(14):1730050,
  32, 2017.

\bibitem{maximeJPMe}
Maxime Breden, Jean-Philippe Lessard, and Jason~D. Mireles~James.
\newblock Computation of maximal local (un)stable manifold patches by the
  parameterization method.
\newblock {\em Indag. Math. (N.S.)}, 27(1):340--367, 2016.

\bibitem{parmChristian}
Jan~Bouwe van~den Berg, Jason~D. Mireles~James, and Christian Reinhardt.
\newblock Computing (un)stable manifolds with validated error bounds:
  non-resonant and resonant spectra.
\newblock {\em J. Nonlinear Sci.}, 26(4):1055--1095, 2016.

\bibitem{fastSlow}
J.~B. van~den Berg and J.~D. Mireles~James.
\newblock Parameterization of slow-stable manifolds and their invariant vector
  bundles: theory and numerical implementation.
\newblock {\em Discrete Contin. Dyn. Syst.}, 36(9):4637--4664, 2016.

\bibitem{manifoldPaper1}
William~D. Kalies, Shane Kepley, and J.~D. Mireles~James.
\newblock Analytic continuation of local (un)stable manifolds with rigorous
  computer assisted error bounds.
\newblock {\em SIAM J. Appl. Dyn. Syst.}, 17(1):157--202, 2018.

\bibitem{jorgeMePerParm}
Jorge Gonzalez and J.~D. Mireles~James.
\newblock High-order parameterization of stable/unstable manifolds for long
  periodic orbits of maps.
\newblock {\em SIAM Journal on Applied Dynamical Systems}, 16(3):1748--1795,
  2017. \url{https://doi.org/10.1137/16M1090041}.

\bibitem{jayChrisParmDDE}
Chris~M. Groothedde and J.~D. Mireles~James.
\newblock Parameterization method for unstable manifolds of delay differential
  equations.
\newblock {\em Journal of Computational Dynamics}, pages 1--52, (First online
  September 2017). \url{doi:10.3934/jcd.2017002}.

\bibitem{MR3783519}
Lei Zhang and Rafael de~la Llave.
\newblock Transition state theory with quasi-periodic forcing.
\newblock {\em Commun. Nonlinear Sci. Numer. Simul.}, 62:229--243, 2018.

\bibitem{MR3797119}
Stavros Anastassiou, Anastasios Bountis, and Arnd B\"{a}cker.
\newblock Recent results on the dynamics of higher-dimensional {H}\'{e}non
  maps.
\newblock {\em Regul. Chaotic Dyn.}, 23(2):161--177, 2018.

\bibitem{MR3705136}
Stavros Anastassiou, Tassos Bountis, and Arnd B\"{a}cker.
\newblock Homoclinic points of 2{D} and 4{D} maps via the parametrization
  method.
\newblock {\em Nonlinearity}, 30(10):3799--3820, 2017.

\bibitem{Cana}
\`Alex Haro, Marta Canadell, Jordi-Llu\'\i~s Figueras, Alejandro Luque, and
  Josep-Maria Mondelo.
\newblock {\em The parameterization method for invariant manifolds}, volume 195
  of {\em Applied Mathematical Sciences}.
\newblock Springer, [Cham], 2016.
\newblock From rigorous results to effective computations.

\bibitem{AMS}
J.~D. Mireles~James.
\newblock Validated numerics for equilibria of analytic vector fields:
  invariant manifolds and connecting orbits.
\newblock {\em Proceedings of Symposia in Applied Mathematics}, 74:1--55, 2018.

\bibitem{MR2821596}
Jan~Bouwe van~den Berg, J.~D. Mireles~James, Jean-Philippe Lessard, and
  Konstantin Mischaikow.
\newblock Rigorous numerics for symmetric connecting orbits: even homoclinics
  of the {G}ray-{S}cott equation.
\newblock {\em SIAM J. Math. Anal.}, 43(4):1557--1594, 2011.

\bibitem{MR3022075}
D.~Ambrosi, G.~Arioli, and H.~Koch.
\newblock A homoclinic solution for excitation waves on a contractile
  substratum.
\newblock {\em SIAM J. Appl. Dyn. Syst.}, 11(4):1533--1542, 2012.

\bibitem{MR3281845}
Gianni Arioli and Hans Koch.
\newblock Existence and stability of traveling pulse solutions of the
  {F}itz{H}ugh-{N}agumo equation.
\newblock {\em Nonlinear Anal.}, 113:51--70, 2015.

\bibitem{MR2644324}
A.~Wittig, M.~Berz, J.~Grote, K.~Makino, and S.~Newhouse.
\newblock Rigorous and accurate enclosure of invariant manifolds on surfaces.
\newblock {\em Regul. Chaotic Dyn.}, 15(2-3):107--126, 2010.

\bibitem{1990mmcm.conf..285S}
C.~{Simo}.
\newblock {On the Analytical and Numerical Approximation of Invariant
  Manifolds}.
\newblock In D.~{Benest} and C.~{Froeschle}, editors, {\em Modern Methods in
  Celestial Mechanics, Comptes Rendus de la 13ieme Ecole Printemps
  d'Astrophysique de Goutelas (France), 24-29 Avril, 1989. Edited by Daniel
  Benest and Claude Froeschle. Gif-sur-Yvette: Editions Frontieres, 1990.,
  p.285}, page 285, 1990.

\bibitem{MR1713086}
Bernd Krauskopf and Hinke Osinga.
\newblock Two-dimensional global manifolds of vector fields.
\newblock {\em Chaos}, 9(3):768--774, 1999.

\bibitem{MR1870261}
Hinke Osinga.
\newblock Non-orientable manifolds of periodic orbits.
\newblock In {\em International {C}onference on {D}ifferential {E}quations,
  {V}ol. 1, 2 ({B}erlin, 1999)}, pages 922--924. World Sci. Publ., River Edge,
  NJ, 2000.

\bibitem{MR2114735}
John Guckenheimer and Alexander Vladimirsky.
\newblock A fast method for approximating invariant manifolds.
\newblock {\em SIAM J. Appl. Dyn. Syst.}, 3(3):232--260, 2004.

\bibitem{MR2834454}
A.~Zanzottera, G.~Mingotti, R.~Castelli, and M.~Dellnitz.
\newblock Intersecting invariant manifolds in spatial restricted three-body
  problems: design and optimization of {E}arth-to-halo transfers in the
  {S}un-{E}arth-{M}oon scenario.
\newblock {\em Commun. Nonlinear Sci. Numer. Simul.}, 17(2):832--843, 2012.

\bibitem{MR1391509}
Michael Dellnitz and Andreas Hohmann.
\newblock The computation of unstable manifolds using subdivision and
  continuation.
\newblock In {\em Nonlinear dynamical systems and chaos ({G}roningen, 1995)},
  volume~19 of {\em Progr. Nonlinear Differential Equations Appl.}, pages
  449--459. Birkh\"{a}user, Basel, 1996.

\bibitem{MR2338026}
M.~E. Henderson.
\newblock Covering an invariant manifold with fat trajectories.
\newblock In {\em Model reduction and coarse-graining approaches for multiscale
  phenomena}, pages 39--54. Springer, Berlin, 2006.

\bibitem{MR2179490}
Michael~E. Henderson.
\newblock Computing invariant manifolds by integrating fat trajectories.
\newblock {\em SIAM J. Appl. Dyn. Syst.}, 4(4):832--882, 2005.

\bibitem{MR2989589}
R.~C. Calleja, E.~J. Doedel, A.~R. Humphries, A.~Lemus-Rodr\'{\i}guez, and
  E.~B. Oldeman.
\newblock Boundary-value problem formulations for computing invariant manifolds
  and connecting orbits in the circular restricted three body problem.
\newblock {\em Celestial Mech. Dynam. Astronom.}, 114(1-2):77--106, 2012.

\bibitem{MR2136745}
B.~Krauskopf, H.~M. Osinga, E.~J. Doedel, M.~E. Henderson, J.~Guckenheimer,
  A.~Vladimirsky, M.~Dellnitz, and O.~Junge.
\newblock A survey of methods for computing (un)stable manifolds of vector
  fields.
\newblock {\em Internat. J. Bifur. Chaos Appl. Sci. Engrg.}, 15(3):763--791,
  2005.

\bibitem{MR2835474}
Roy~H. Goodman and Jacek~K. Wr\'{o}bel.
\newblock High-order bisection method for computing invariant manifolds of
  two-dimensional maps.
\newblock {\em Internat. J. Bifur. Chaos Appl. Sci. Engrg.}, 21(7):2017--2042,
  2011.

\bibitem{MR3021639}
Jacek~K. Wr\'{o}bel and Roy~H. Goodman.
\newblock High-order adaptive method for computing two-dimensional invariant
  manifolds of three-dimensional maps.
\newblock {\em Commun. Nonlinear Sci. Numer. Simul.}, 18(7):1734--1745, 2013.

\bibitem{shaneNumericalPaper}
Shane Kepley and J.~D. Mireles~James.
\newblock Homoclinic dynamics in a restricted four body problem: a
  multi-parameter study of transverse connections for the saddle-focus
  equilibrium solutions.
\newblock {\em (submitted to Celestial Mechanics and Dynamical Astronomy)},
  2018.

\bibitem{MR699057}
S.~Newhouse, J.~Palis, and F.~Takens.
\newblock Bifurcations and stability of families of diffeomorphisms.
\newblock {\em Inst. Hautes \'{E}tudes Sci. Publ. Math.}, (57):5--71, 1983.

\bibitem{doedel_paris}
E.~J. Doedel.
\newblock Lecture notes on numerical analysis of bifurcation problems.
\newblock In {\em International Course on Bifurcations and Stability in
  Structural Engineering}. Universit\'e Pierre et Marie Curie (Paris VI),
  November 2000.

\bibitem{MR910499}
H.~B. Keller.
\newblock {\em Lectures on numerical methods in bifurcation problems},
  volume~79 of {\em Tata Institute of Fundamental Research Lectures on
  Mathematics and Physics}.
\newblock Published for the Tata Institute of Fundamental Research, Bombay,
  1987.
\newblock With notes by A. K. Nandakumaran and Mythily Ramaswamy.

\bibitem{MR1314079}
A.~R. Champneys and Yu.~A. Kuznetsov.
\newblock Numerical detection and continuation of codimension-two homoclinic
  bifurcations.
\newblock {\em Internat. J. Bifur. Chaos Appl. Sci. Engrg.}, 4(4):785--822,
  1994.

\bibitem{MR0228014}
S.~Smale.
\newblock Differentiable dynamical systems.
\newblock {\em Bull. Amer. Math. Soc.}, 73:747--817, 1967.

\bibitem{MR516994}
S.~Newhouse, D.~Ruelle, and F.~Takens.
\newblock Occurrence of strange {A}xiom {A} attractors near quasiperiodic flows
  on {$T^{m}$},{$\,m\geq 3$}.
\newblock {\em Comm. Math. Phys.}, 64(1):35--40, 1978/79.

\bibitem{MR908665}
Marcy Barge.
\newblock Homoclinic intersections and indecomposability.
\newblock {\em Proc. Amer. Math. Soc.}, 101(3):541--544, 1987.

\bibitem{MR1277856}
Judy Kennedy.
\newblock How indecomposable continua arise in dynamical systems.
\newblock In {\em Papers on general topology and applications ({M}adison, {WI},
  1991)}, volume 704 of {\em Ann. New York Acad. Sci.}, pages 180--201. New
  York Acad. Sci., New York, 1993.

\bibitem{MR968418}
Andreas Floer.
\newblock A topological persistence theorem for normally hyperbolic manifolds
  via the {C}onley index.
\newblock {\em Trans. Amer. Math. Soc.}, 321(2):647--657, 1990.

\bibitem{maciejContinum}
Maciej~J. Capi\'nski and Hieronim Kubica.
\newblock Persistence of normally hyperbolic invariant manifolds in the absence
  of rate conditions.
\newblock {\em (submitted)}, 2019.

\bibitem{MR591900}
W.~F. Langford and G.~Iooss.
\newblock Interactions of {H}opf and pitchfork bifurcations.
\newblock In {\em Bifurcation problems and their numerical solution ({P}roc.
  {W}orkshop, {U}niv. {D}ortmund, {D}ortmund, 1980)}, volume~54 of {\em
  Internat. Ser. Numer. Math.}, pages 103--134. Birkh\"{a}user, Basel-Boston,
  Mass., 1980.

\bibitem{MR508917}
M.~Golubitsky and D.~Schaeffer.
\newblock A theory for imperfect bifurcation via singularity theory.
\newblock {\em Comm. Pure Appl. Math.}, 32(1):21--98, 1979.

\bibitem{MR950168}
Martin Golubitsky, Ian Stewart, and David~G. Schaeffer.
\newblock {\em Singularities and groups in bifurcation theory. {V}ol. {II}},
  volume~69 of {\em Applied Mathematical Sciences}.
\newblock Springer-Verlag, New York, 1988.

\bibitem{MR633825}
Martin Golubitsky and William~F. Langford.
\newblock Classification and unfoldings of degenerate {H}opf bifurcations.
\newblock {\em J. Differential Equations}, 41(3):375--415, 1981.

\bibitem{maciej_us_invariantTori}
Maciej~J Capinski, Emmanuel Fleurantin, and J.~D. Mireles~James.
\newblock Computer assisted proofs of two-dimensional attracting invariant tori
  for {O}{D}{E}s.
\newblock arXiv:1905.08116.

\bibitem{BDLM}
Jan~Bouwe van~den Berg, Andr\'ea Desch\^enes, Jean-Philippe Lessard, and
  Jason~D. Mireles~James.
\newblock Stationary coexistence of hexagons and rolls via rigorous
  computations.
\newblock {\em SIAM J. Appl. Dyn. Syst.}, 14(2):942--979, 2015.

\bibitem{cnLohner}
Daniel Wilczak and Piotr Zgliczynski.
\newblock $c^n$-lohner algorithm.
\newblock {\em Scheade Informaticae}, 20:9--46, 2011.

\bibitem{MR1961956}
Daniel Wilczak and Piotr Zgliczynski.
\newblock Heteroclinic connections between periodic orbits in planar restricted
  circular three-body problem---a computer assisted proof.
\newblock {\em Comm. Math. Phys.}, 234(1):37--75, 2003.

\bibitem{MR2262261}
Daniel Wilczak.
\newblock Symmetric homoclinic solutions to the periodic orbits in the
  {M}ichelson system.
\newblock {\em Topol. Methods Nonlinear Anal.}, 28(1):155--170, 2006.

\bibitem{MR1816799}
Gianni Arioli and Piotr Zgliczy{\'n}ski.
\newblock Symbolic dynamics for the {H}\'enon-{H}eiles {H}amiltonian on the
  critical level.
\newblock {\em J. Differential Equations}, 171(1):173--202, 2001.

\bibitem{MR3443692}
Maciej~J. Capi\'nski and Anna Wasieczko-Zajac.
\newblock Geometric proof of strong stable/unstable manifolds with application
  to the restricted three body problem.
\newblock {\em Topol. Methods Nonlinear Anal.}, 46(1):363--399, 2015.

\bibitem{MR3032848}
Maciej~J. Capi{\'n}ski.
\newblock Computer assisted existence proofs of {L}yapunov orbits at {$L_2$}
  and transversal intersections of invariant manifolds in the {J}upiter-{S}un
  {PCR}3{BP}.
\newblock {\em SIAM J. Appl. Dyn. Syst.}, 11(4):1723--1753, 2012.

\bibitem{MR2271217}
Daniel Wilczak.
\newblock The existence of {S}hilnikov homoclinic orbits in the {M}ichelson
  system: a computer assisted proof.
\newblock {\em Found. Comput. Math.}, 6(4):495--535, 2006.

\bibitem{MR2173545}
Daniel Wilczak.
\newblock Symmetric heteroclinic connections in the {M}ichelson system: a
  computer assisted proof.
\newblock {\em SIAM J. Appl. Dyn. Syst.}, 4(3):489--514 (electronic), 2005.

\bibitem{MR2533624}
Zin Arai, William Kalies, Hiroshi Kokubu, Konstantin Mischaikow, Hiroe Oka, and
  Pawe{\l} Pilarczyk.
\newblock A database schema for the analysis of global dynamics of
  multiparameter systems.
\newblock {\em SIAM J. Appl. Dyn. Syst.}, 8(3):757--789, 2009.

\bibitem{MR3504856}
Tomoyuki Miyaji, Pawe\l Pilarczyk, Marcio Gameiro, Hiroshi Kokubu, and
  Konstantin Mischaikow.
\newblock A study of rigorous {ODE} integrators for multi-scale set-oriented
  computations.
\newblock {\em Appl. Numer. Math.}, 107:34--47, 2016.

\bibitem{alexCAPKAM}
Jordi-Llu{\'{\i}}s Figueras, Alex Haro, and Alejandro Luque.
\newblock Rigorous computer assisted application of kam theory: a modern
  approach.
\newblock {\em (Submitted) arXiv:1601.00084 [math.DS]}, 2016.

\end{thebibliography}
}

\end{document}